\documentclass[a4paper,12pt,twoside]{amsart}                                   %
\usepackage{amsmath,amssymb,amsfonts}
\usepackage[french]{babel}
\usepackage[latin1]{inputenc}
\usepackage[T1]{fontenc}
\usepackage{pslatex}

\usepackage{hyperref}
\newtheorem*{prop}{Proposition}
\newtheorem*{thm}{Th\'eor\`eme}
\newtheorem*{lem}{Lemme}
\newtheorem*{cor}{Corollaire}
\theoremstyle{definition}
\newtheorem*{defn}{D\'efinition}
\theoremstyle{remark}

\def \dem {\begin{proof}[Démonstration]}

\def \id { {\rm id\, }}
\def \dim { {\rm dim\, }}
\def \Ad { {\rm Ad\, }}
\def \Vect { {\rm Vect\, }}
\def \en{ {\rm end\, }}

 \title[C*-groupoides quantiques et  inclusions de facteurs]
 {C*-groupoides quantiques et  inclusions de facteurs~:\\
 Structure sym\'etrique et autodualit\'e, \\ action sur le facteur hyperfini de type $\mathrm{II}_{1}$ } 
\author{Marie-Claude DAVID}   
\address{Mathématique, Bât. 425,  Université Paris-Sud, F-91405 Orsay 
Cedex.\hfil \break  
\indent mcld@math.u-psud.fr}
\begin{document}        

\maketitle

\textsc{Résumé :}
Etant données une inclusion $N_{0} \subset N_{1}$ de facteurs de type $\mathrm{II}_{1}$ de profondeur 
$2$ et d'indice fini et 
$$N_0 \,\subset \, N_1 \, \subset  \, N_2 \, \subset 
\,  N_3 \dots $$
 la tour de Jones correspondante, D. Nikshych et L. Vainerman 
ont muni les commutants relatifs $N'_{0} \cap N_{2}$ 
et $N'_{1} \cap N_{3}$  de structures duales de C*-groupoïde 
quantique.

 Je  modifie ici   la dualité 
et j'obtiens ainsi une construction symétrique qui n'exige pas une 
nouvelle définition des involutions. Alors
les algèbres de Temperley-Lieb sont des C*-groupoïdes 
quantiques autoduaux ; plus généralement on peut associer à une inclusion de 
profondeur finie et d'indice fini un  C*-groupoïde 
quantique autodual.

Je montre que tout C*- groupoïde quantique connexe de dimension finie 
agit extérieu\-rement sur le facteur hyperfini de type $\mathrm{II}_{1}$.
A la lumière de ce cas particulier, je propose une déformation de tout 
C*-groupoïde quantique fini en un C*-groupoïde quantique fini régulier.
\vskip  1cm

\textsc{Abstract :}
Let  $N_{0} \subset N_{1}$ a depth $2$, finite index inclusion of type 
$\mathrm{II}_{1}$ factors and 
$$N_0 \,\subset \, N_1 \, \subset  \, N_2 \, \subset \,  N_3 \dots $$
the corresponding Jones tower. D. Nikshych and L. Vainerman built dual structures of quantum C*-groupoid on 
the relative commutants $N'_{0} \cap N_{2}$ et $N'_{1} \cap N_{3}$. 

Here I define a new duality which allows a symmetric construction 
without changing the involution. So the Temperley-Lieb algebras are 
selfdual quantum C*-groupoids and the quantum C*-groupoids 
associated to a finite depth finite index inclusion can be chosen 
selfdual.

I show that every finite-dimensional connected quantum C*-groupoid acts 
outerly on the type $\mathrm{II}_{1}$ hyperfinite factor. In the 
light of this particular case, I propose a deformation of any finite quantum 
C*-groupoid to a regular finite quantum C*-groupoid.

\vskip  1cm

\textsc{Code Matière AMS :} 46L37, 16W30, 57T05, 22D35.
\vskip  5mm
\textsc{Mots Clefs :} Subfactors, quantum groupoids, Temperley-Lieb 
algebras, crossed product, action.
\vskip 1cm

\newpage
\tableofcontents

\section{Introduction}

Soient une inclusion de facteurs de type $\mathrm{II}_{1}$, $N_{0} \subset 
N_{1}$,  de profondeur 
$2$ d'indice fini~et 
$$N_0 \,\subset \, N_1 \, \raisebox{1.3 ex}{$\begin{matrix} f_{1}\\ \subset \end{matrix}$} \, 
N_2 \,\raisebox{1.3 ex}{$\begin{matrix} f_{2}\\ \subset \end{matrix}$} 
\,  N_3 \dots $$
la tour obtenue par construction de base. Si l'inclusion est irréductible, les commutants relatifs $A=N'_{0} \cap N_{2}$ 
et $B =N'_{1} \cap N_{3}$ peuvent être munis de structures duales d'algèbres 
de Kac de dimension finie ([Da1], [L], [Szy]).
Généralisant les méthodes de W. Szymanski aux inclusions réductibles, 
D. Nikshych et L. Vainerman 
définissent dans [NV 1] une dualité entre les commutants relatifs $A$ 
et $B$ à l'aide de la trace $tr$ de $N_{1}$ : 
$$\langle a,b\rangle = [N_{1}:N_{0}]^{2} tr(af_{2}f_{1}Hb)\quad (a\in A, \;b\in 
B)$$
(L'élément qui  rend compte du fait que l'inclusion n'est pas 
irréductible est l'indice $H$ de 
la restriction à $N'_{1} \cap N_{2}$ de la trace $tr$.)
A l'aide de cette dualité, ils  définissent des structures duales de C*-groupoïde 
quantique  sur $A$ et $B$. Les coproduits sont définis par dualité, aussi pour 
qu'ils soient compatibles avec l'involution, ils ont dû définir de  nouvelles 
involutions sur $A$ et $B$ différentes de celles héritées du facteur $N_{4}$.
L'étude de ces structures à l'aide de formules g\'en\'eralisant 
celles obtenues dans [Da] pour une inclusion irréductible fait 
apparaître un autre  inconvénient : 
Si on note $B(N_{1}\subset N_{2})$ la structure définie sur $B$, la 
structure duale sur $A$ n'est pas isomorphe, mais anti-isomorphe à $B(N_{0}\subset N_{1})$.

Je propose ici une  autre dualité qui permet une construction symétrique 
conservant l'involution :
$$\langle a,b\rangle = [N_{1}:N_{0}]^{2} tr(aH^{1/2}f_{2}f_{1}H^{1/2}b)\quad (a\in A, \;b\in 
B)$$
Avec cette nouvelle définition, si on note $B(N_{1}\subset N_{2})$ la structure définie sur $B$, la 
structure duale sur $A$ est $B(N_{0}\subset N_{1})$.

Les propriétés des structures de C*-groupoïde quantique construites à 
partir de cette dualité pourraient s'obtenir à partir des 
résultats de D. Nikshych et L. Vainerman.  Certaines démonstrations sont 
d'ailleurs fortement inspirées des leurs. Pourtant, comme 
je dispose maintenant de formules 
pour les co-produits et les antipodes, je donne souvent des 
démonstrations directes. 

L'intérêt de cette construction symétrique apparaît dans la partie 4. 
Dans le cas d'une inclusion de profondeur finie, on peut alors obtenir 
des  C*-groupoïdes quantiques  autoduaux. 

Dans la partie 5, je précise la structure de C*-groupoïde 
quantique des algèbres de Temper\-ley-Lieb. Grâce à la symétrie de la 
construction, ces C*-groupoïdes 
quantiques sont autoduaux. J'étudie en particulier le 
C*-groupoïde quantique de dimension 13 associé au graphe linéaire $A_{4}$ et 
montre qu' il est isomorphe à celui décrit
par G. Böhm et K. Szlachanyi dans [BSz - 5].

Dans la partie 6, j'étends aux C*-groupoïdes 
quantiques connexes de dimension finie un résultat de   
D. Nikshych [N] sur les algèbres de Kac faibles en 
les faisant agir extérieurement sur le facteur hyperfini de type $\mathrm{II}_{1}$.

Dans la partie 7, je montre qu'on peut 
déformer toute paire de C*-groupoïdes quantiques finis duaux 
en une paire de C*-groupoïdes quantiques finis réguliers sans modifier la structure de C*-algèbre.
    
L'essentiel de ce texte a été écrit au printemps 2001, 
la partie 6 l'a complété au printemps 2003 et la 7 à l'automne de la 
même année. Léonid Vainerman et Jean-Michel Vallin  ont été à 
l'origine de cet travail, je les en remercie vivement et plus particulièrement 
Léonid pour de nombreux échanges par courrier électronique à propos 
de la construction originale du C*-groupoide quantique. 
Mes remerciements vont aussi à Kornel Szlach\'anyi  pour ses réponses précises. Les calculs 
concernant le groupoïde quantique de dimension 13 ont été grandement 
facilités par les conseils Maple de Jacques Peyrière, je lui en suis 
reconnaissante.

   \section{$C^*$-groupoïdes quantiques finis} On rappelle ici les 
définitions de C*-groupoïde quantique fini et de C*-groupoïde quantique 
fini dual ([BNSz][N][NV 2] [NV 3]) ainsi que celles d'une action et du 
produit croisé.

\subsection{$C^*$-groupoïde quantique fini}\label{def}
\begin{defn}
Un {\bf C*-groupoïde quantique fini}  est une C*-algèbre $G$ de dimension 
finie (on note  $m$ la multiplication, $1$ l'unité , $^{*}$ 
l'involution) munie 
	d'une structure de co-algèbre associative avec  un coproduit 
	$\Delta$, une 
	co-unité $\varepsilon$ et une antipode $S$ tels que 
	
	i) $\Delta$ soit un $^{*}$-homomorphisme d'algèbres de $G$ dans $G 
	\otimes G$ vérifiant :
	$$(\Delta \otimes \id)\Delta(1)=(1\otimes \Delta(1))(\Delta(1)\otimes 1)$$
	
	ii) La co-unité soit une application linéaire de $G$ dans $G$ vérifiant :
	$$\varepsilon(fgh) 
	=\varepsilon(fg_{(1)})\varepsilon(g_{(2)}h)\qquad ((f,g,h)\in G^3)$$
(propriété est équivalente à 
		$$\varepsilon(fgh)=\varepsilon(fg_{(2)})\varepsilon(g_{(1)}h) \qquad ((f,g,h)\in 
		G^3))$$
	
	iii) L'antipode $S$ soit  un 
	anti-homomorphisme d'algèbre et de co-algèbre de $G$ dans $G$ 
	vérifiant pour tout $g$ de $G$ :
$$	m(\id \otimes S)\Delta(g)=(\varepsilon \otimes 
	\id)(\Delta(1)(g \otimes 1))$$
(propriété équivalente à
$$m(S \otimes \id  )\Delta(g)=(\id \otimes \varepsilon)((1 \otimes 
		g)\Delta(1))$$
\end{defn}

On appelle {\bf co-unité but} et {\bf co-unité source} les 
applications $\varepsilon_t$ et $\varepsilon_s$ définies pour tout $g$ de $G$ 
par :
$$\varepsilon_t(g)=(\varepsilon \otimes \id)(\Delta(1)(g \otimes 1))
	\qquad
	\varepsilon_s(g)=(\id \otimes \varepsilon)((1 \otimes g)\Delta(1))$$

\subsection{$C^*$-groupoïde quantique dual}
\begin{defn}
On définit sur $\hat G=Hom_{\mathbb C}(G,\mathbb C)$ une structure 
de   C*-groupoïde quantique dual de celle de $G$ grâce aux formules suivantes :
\begin{align*}
\langle h,\phi \psi\rangle&=\langle \Delta(h),\phi \otimes 
\psi\rangle\\
\langle h \otimes g,\hat{\Delta}(\phi) \rangle&=\langle hg,\phi \rangle\\
\langle h ,\hat{S}(\phi) \rangle &=\langle S(h),\phi \rangle\\
\langle h,\phi^{*} \rangle&=\overline{\langle S(h)^{*} ,\phi 
\rangle}
\end{align*}
pour  tous $\phi$ et  $\psi$ de  $\hat G$ et tous $h$ et $g$ de 
$G$.

L'unité de  $\hat G$ est $\varepsilon$ et la co-unité $\hat{\varepsilon}$ 
 est $\phi \mapsto  \langle 1, \phi \rangle$.
\end{defn}

\subsection{Projection de Haar, mesure de Haar}\label{Haar}
D'après [BNSz 4.5] et [NV 1-7.3.1], il existe une  unique projection $p$ de $G$ invariante par 
l'antipode, appelée 
{\bf projection de Haar normalisée} telle que pour tout $g$ de $G$, on 
ait  les 
 propriétés équivalentes suivantes :
\begin{align*}
\quad \text{(i)} \quad \varepsilon_{t}(g)p&=gp  &\varepsilon_{t}(p)=1\\
\quad \text{(ii)} \quad p\varepsilon_{s}(g)&=pg  &\varepsilon_{s}(p)=1
\end{align*}

La forme linéaire duale $\hat \phi$ de la projection de Haar normalisée est appelée {\bf mesure de 
Haar normalisée de $\hat G$}. Elle est fidèle, invariante par l'antipode et vérifie 
les 
 propriétés équivalentes suivantes :
\begin{align*}
\text{(i)} \quad (\id \otimes \hat \phi ) \hat \Delta&=(\hat 
\varepsilon_t \otimes \hat \phi ) \hat \Delta 
&\hat \phi \circ \hat \varepsilon_{t}=\hat \varepsilon\\
\text{(ii)} \quad (\hat \phi \otimes \id) \hat \Delta&=(\hat \phi \otimes \hat 
\varepsilon_{s}) \hat \Delta &
  \hat \phi \circ \hat \varepsilon_{s}=\hat \varepsilon
\end{align*}

   \subsection{Les sous-algèbres co-unitales}\label{cartan} [NV3 - 
   2.2][BNSz - 2.5 et 2.9]
   
   L'algèbre $G_{s}=\varepsilon_s(G)$ (resp. $G_{t}=\varepsilon_t(G)$) est 
appelée sous-algèbre co-unitale source (resp. sous-algèbre co-unitale 
but). 

Les co-unités but et source sont des homomorphismes idempotents de $G_{t}$ 
(resp. $G_{s}$)   et vérifient pour 
tout $g$ de $G$ :
$$(\id \otimes \varepsilon_t)\Delta(g)=1_{(1)}g \otimes 1_{(2)} \qquad 
(\varepsilon_s \otimes \id )\Delta(g)=1_{(1)}\otimes g 1_{(2)}$$
On a aussi les formules suivantes :
$$\varepsilon_t \circ S=\varepsilon_t \circ\varepsilon_s=S 
\circ \varepsilon_s\qquad \qquad\varepsilon_s \circ S 
=\varepsilon_s\circ\varepsilon_t=S \circ \varepsilon_t$$
Les sous-algèbres co-unitales commutent entre elles et vérifient :
\begin{align*}
    G_t & = & \{g\in G, \Delta(g)=1_{(1)}g \otimes 1_{(2)}=g1_{(1)} \otimes 1_{(2)}\} 
    & = & \{(\omega \otimes \id)\Delta(1), \omega \in \hat{G}\} \\
    G_s & = & \{g\in G, \Delta(g)=1_{(1)} \otimes g1_{(2)}=1_{(1)} \otimes 1_{(2)}g\} 
    & = & \{( \id \otimes \omega )\Delta(1), \omega \in \hat{G}\} 
\end{align*}

\subsection{C*-groupoïde quantique fini régulier} On dit que le C*-groupoïde quantique fini
 $G$ est  {\bf régulier} si son antipode est involutive sur les algèbres 
 co-unitales.
    \subsection{Action d'un groupoïde quantique}[NSzW - def.1.2.2] 
    [N - 2.2]
    \subsubsection{} Soit M une algèbre involutive unitaire. On dit qu'un groupoïde quantique 
    $G$
    fini agit à gauche (resp. à droite) sur $M$ s'il existe une 
    application linéaire $g\otimes m \mapsto g\triangleright m$ de $G\otimes M$ 
    dans $M$ (resp. $m\otimes g \mapsto m\triangleleft g$ de $M\otimes 
    G$ 
    dans $M$) 
    définissant une structure de $G$-module à gauche (resp. à droite) 
    sur $M$  et  vérifiant pour $g$ dans $G$ et $x$ et $y$ dans $M$
    
    \begin{enumerate}
        \item  $g \triangleright(xy)= (g _{(1)}\triangleright x)(g 
        _{(2)}\triangleright y)$\;\;(resp. $(xy)\triangleleft 
        g=(x\triangleleft g _{(1)}) (y\triangleleft g _{(2)})$)
    
        \item  $(g \triangleright x)^{\ast} =S(g)^{\ast} \triangleright x^{\ast}$
	\;\;(resp. $(x   \triangleleft g)^{\ast}= x^{\ast}\triangleleft S(g)^{\ast}$)
    
        \item $ g\triangleright 1 = \varepsilon_t(g) \triangleright 
        1$ \;\;(resp. $1\triangleleft g= 1 \triangleleft 
        \varepsilon_s(g)$)

    \end{enumerate}
    
    Si $M$ est une $C^*$-algèbre ou une algèbre de von Neumann, 
    l'application \break $g\otimes m \mapsto g\triangleright m$ 
    (resp. $m\otimes g \mapsto m\triangleleft g$) doit être continue en norme ou faiblement pour 
    tout $g$ de $G$.
    
     D'après [NSzW - def 1.2.4], une action à gauche est dite {\bf standard} 
     si l'application $x\otimes 1_M \mapsto x\triangleright 
    1_{M}$  est un isomorphisme de $A_{t}$ sur une sous-algèbre de $M$.
    
    Une action à gauche est  standard si et seulement si elle 
    vérifie :
    $$g\triangleright 1 =  0 \; \Leftrightarrow \; \varepsilon_t(g) = 0$$ 
    
    \subsubsection{}\label{dual} Rappelons les définitions des actions duales des 
  groupoïdes l'un sur l'autre (voir  [NSzW]  ou [N]). Posons $A=G$ 
  et $B=\hat{G}$.
  Le groupoïde $A$ agit à droite sur $B$ :
  $$b \triangleleft a=\langle a,b_{(1)}\rangle b_{(2)}\qquad (a\in A, b\in B)$$
  Le groupoïde $B$ agit à gauche sur $A$ :
  $$b \triangleright a=\langle a_{(2)},b\rangle a_{(1)}\qquad (a\in A, b\in B)$$
   De façon équivalente, pour tous $x$ et $a$ de $A$ et $y$ et $b$ de $B$, on a 
   $$\langle x,b \triangleleft a\rangle= \langle ax,b\rangle\quad 
   \text{et}\quad \langle b \triangleright a,y\rangle= \langle a,yb\rangle.$$
 
    \subsubsection{}\label{AsBt} Les actions que nous venons de définir 
   sont standard.
  \begin{prop}[BNSz-lemme 2.6] L'application $x \mapsto 
  1_{b}\triangleleft x$ est un isomorphisme 
  de l'algèbre $A_{s}$ sur l'algèbre $B_{t}$. Sa 
  réciproque est donnée par $y \mapsto y \triangleright 1_{a}$
  ($y \in B_{t}$).
  \end{prop}
  
  \subsubsection{} \label{rem27} [BNSz - 2.7] Des propriétés des 
  sous-algèbres co-unitales, on déduit les formules suivantes pour $b$ 
  dans $B$, $x$ dans $A_{t}$ et $y$ dans $A_{s}$ :
  \begin{align*}
 & x \triangleright  b	 =( x\triangleright 1_{b})\ b 
  &\quad & y\triangleright b = b\ (y \triangleright  1_{b})\\
  	&b\triangleleft  x = (1_{b}\triangleleft x)\ b   &\quad &b\triangleleft y = b\ (1_{b}\triangleleft y) 
  	\end{align*}

   \subsection{Produit croisé d'une algèbre par un groupoïde quantique}
   \subsubsection{Définition}[N - 2.2]
   Le produit croisé à gauche (resp. à droite) $M \rtimes G$ (resp. $G 
   \ltimes M$) est le $\mathbb{C}$-espace 
   vectoriel $M \otimes_{G_t} G$ (resp. $G \otimes_{G_s} M$) où on 
   identifie $m(z \triangleright 1)\otimes  g$ et $m\otimes zg$ 
   (resp. $gz\otimes  m$ et $g\otimes (1\triangleleft z)m$) pour $m$ 
   dans $M$, $g$ dans $G$ et $z$ dans $G_t$ (resp. $G_s$). 
   Soit $[m \otimes g]$ (resp. $[g\otimes m]$) 
   la classe de $m \otimes g$ (resp. $g\otimes m$).
   
   On munit le produit 
   croisé d'une structure d'une algèbre involutive en posant pour 
   tous $g$ et $h$ dans $G$ et $x$ et $y$ dans $M$ :
   \begin{align*}
        &[x \otimes g][y \otimes h]& = [x (g_{(1)} \triangleright y) \otimes 
   g_{(2)}h] & \qquad  & [x \otimes g]^{\ast} & =[(g_{(1)}^{\ast} \triangleright x^{\ast}) \otimes g^{\ast}_{(2)}] &   \\
      \text{(resp. }&[g \otimes x][h \otimes y]& = [g h_{(1)}  \otimes 
     (x \triangleleft h_{(2)})y] &  \qquad & [g \otimes x]^{\ast}  & 
     =[g_{(1)}^{\ast} \otimes  (x^{\ast} \triangleleft g^{\ast}_{(2)}) ] & \text{)}
   \end{align*}
  
   De plus si $M$ est une $C^*$-algèbre ou une algèbre de von Neumann, 
   le produit croisé devient une $C^*$-algèbre ou une algèbre de von Neumann.
   
   Les applications $i_G : g \mapsto [1\otimes g]$ (resp. $g \mapsto 
   [g \otimes 1]$) et $i_M : m \mapsto 
   [m \otimes 1]$ (resp. $m \mapsto 
   [1 \otimes m]$) sont des homomorphismes injectifs d'algèbres 
   involutives de $G$ et $M$ dans le produit 
   croisé telles que :
   $$M \rtimes G = i_M(M)i_G(G)\quad \text{(resp. } G \ltimes M = 
   i_G(G)i_M(M)\; \text{ ) }$$
   
   \subsubsection{Action duale sur le produit croisé} \label{dualprod}[N - 2.2].
   On définit l'action duale à gauche (resp. à droite) de $\hat{G}$ 
   sur $M\rtimes G$ (resp. $G\ltimes M$) par  :
   $$\qquad h \triangleright [m\otimes g]=[m\otimes h \triangleright 
   g]\qquad (g \in G, h\in \hat{G}, m\in M)$$
   $$\text{(resp. }  [g\otimes m]\triangleleft h=[g\triangleleft h\otimes 
   m]\qquad (g \in G, h\in \hat{G}, m\in M))$$
   
   \subsection{Inclusions de profondeur 2} \label{NSzW} Nous rappelons ici 
   quelques résultats de [NSzW] qui motivent cet article. Nous 
   considèrons un $C^*$-groupoïde quantique $A$ de dimension finie agissant sur 
    une algèbre de von Neumann $M$.
   
   \subsubsection{Produit croisé et tour de Jones}\label{415}
   \begin{cor}[NSzW - 4.1.5] Soient $N$ et $M$ des algèbres de von 
   Neumann et $A$ un $C^*$-groupoïde quantique de dimension finie tel que 
   $M$ soit $N\rtimes \hat{A}$. La tour 
   $$N \subset M \subset M\rtimes A \subset M\rtimes A\rtimes 
   \hat{A}\dots$$
   est une tour de Jones de profondeur 2.
   \end{cor}
   
   En particulier, si on prend $N=A_{t}$ et $M=A_{t}\rtimes A=A$, 
   alors 
   $$A_{t} \subset A \subset A\rtimes \hat{A} \subset  A\rtimes 
   \hat{A}\rtimes A\dots$$
    est une tour de Jones.
   \subsubsection{Action extérieure}\label{ext} Nous prendrons le résultat du théorème suivant comme définition 
   pour une action extérieure.
   \begin{thm}[NSzW - 4.2.3]
   L'action de $A$ sur $M$ est extérieure si et seulement si on a 
   l'égalité :
   $$M' \cap (M \rtimes A)=Z(M) \rtimes A_{s}$$	
   \end{thm}
   
   \subsubsection{$C^*$-groupoïde quantique connexe}\label{connexe}
   \begin{defn}
   Un $C^*$-groupoïde quantique $A$ est dit connexe (terminologie de [N] que 
   nous gardons car elle fait référence à l'inclusion $A_{s} \subset 
   A$) ou pur (terminologie de [NSzW]) si $A_{s} \cap Z(A)$ est réduit 
   à $\mathbb{C}$.
   \end{defn}
   
   \begin{prop}[NSzW - 2.4.6] Les conditions suivantes sont équivalentes :
   \begin{enumerate}
   	\item  $A$ est connexe
   
   	\item  $\hat{A}_{s} \cap \hat{A}_{t}=\mathbb{C}$
   
   	\item  $A_{t} \cap Z(A)=\mathbb{C}$
   \end{enumerate}
   \end{prop}
   
   \noindent{\it Remarques}
   \begin{enumerate}
   	\item  Si $A$ est connexe, toute action est standard ([NSzW - 2.2.1]).
   
   	\item D'après le théorème 3.1.1 de [NSzW], le centre de $M \rtimes A$ 
   contient nécessairement $1_{M} \rtimes (A_{t}\cap Z(A))$. Donc la 
   connexité de $A$ est nécessaire pour obtenir un facteur comme 
   produit croisé.
  \item D'après le corollaire 2.4.4 de [NSzW], si l'action de $A$ 
  sur  $M$ est standard, le centre de $M$ contient nécessairement une 
sous-algèbre isomorphe à $A_{t}\cap A_{s}$ donc la connexité de 
$\hat{A}$ est nécessaire à l'action de $A$ sur un facteur.
 \end{enumerate}
   
 \subsubsection{Action extérieure et facteur}\label{facteur}
 \begin{thm}[NSzW - 4.2.4]
   Si $A$ agit extérieurement  et de façon standard sur un facteur $M$,
   $M \rtimes A$ est un facteur si et seulement si $A$ est connexe.
   \end{thm}
   
   \begin{thm}[NSzW -  4.2.5]
   Si $A$ est connexe et agit extérieurement sur un facteur $M$ alors 
  la tour 
   $$M^{A} \subset M \subset M\rtimes A \subset M\rtimes A\rtimes 
   \hat{A}\dots$$
   est une tour de Jones de facteurs. De plus le groupoïde dual 
   $\hat{A}$ est aussi connexe et son action canonique sur $M\rtimes 
   A$ est extérieure.
   \end{thm}
   \subsubsection{Tour dérivée}
   Le résultat suivant précise la tour dérivée de 
   l'inclusion obtenue par l'action de $A$.
   
   \begin{cor} [NSzW -  4.3.5]
   Soient $A$ un $C^*$-groupoïde quantique fini agissant exté\-rieurement sur un facteur $M$
   et $N$ la sous-algèbre des points fixes de $M$ sous $A$. On a les 
   égalités suivantes :
   \begin{eqnarray*}
   		N' \cap M & = & 1_{M}\rtimes A_{t}  \\
   		M' \cap M \rtimes A& = & 1_{M}\rtimes A_{s}  \\
   		N' \cap M \rtimes A& = & 1_{M}\rtimes A
   	\end{eqnarray*}	
   \end{cor}

\section{$C^*$-groupoïdes quantiques associés à une inclusion 
d'indice fini de 
profondeur 2 de facteurs de type $\mathrm{II}_{1}$.}

 Soit $N_{0} \subset N_{1}$ 
une inclusion d'indice fini $\tau^{-1}$ de facteurs de 
type $\mathrm{II}_{1}$. On note 
$$N_0 \,\subset \, N_1 \, \raisebox{1.3 ex}{$\begin{matrix} f_{1}\\ \subset \end{matrix}$} \, 
N_2 \,\raisebox{1.3 ex}{$\begin{matrix} f_{2}\\ \subset \end{matrix}$} 
\,  N_3 \, \raisebox{1.3 ex}{$\begin{matrix} f_{3}\\ \subset \end{matrix}$} \, 
N_4 \dots \subset N_{n}  \, \raisebox{1.3 ex}{$\begin{matrix} f_{n}\\ \subset \end{matrix}$} \, 
N_{n+1} \dots$$
 la tour de Jones obtenue par construction de base [G.H.J. 3] 
et  $tr$ la trace normale finie normalisée sur les 
facteurs considérés. 

On suppose que l'inclusion $N_{0} \subset N_{1}$ est de profondeur $2$ c'est-à-dire qu'elle 
vérifie l'une des conditions équivalentes suivantes :
\begin{enumerate}
	\item  
le commutant relatif $N'_{0} \cap N_{3}$ est obtenu par construction de 
base à partir de $N'_{0} \cap N_{1} \subset N'_{0} \cap N_{2}$
	\item  
l'algèbre $N'_{0} \cap N_{3}$ est linéairement engendrée  par $(N'_{0} \cap 
N_{2})f_{2}(N'_{0} \cap N_{2})$
	\item  $\dim Z(N'_{0} \cap N_{1})=\dim Z(N'_{0} \cap N_{3})$
\end{enumerate}

On remarque que $N_{0} \subset N_{1}$ est de profondeur $2$ si et 
seulement si $N_{1} \subset N_{2}$ est de profondeur~$2$.

\subsection{Anti-automorphismes associés à la tour dérivée.}
Soit $J_{n}$  l'isomé\-trie bijective 
anti-linéaire canonique de l'espace standard $L^2(N_{n},tr)$ de 
$N_{n}$ ($n \in \mathbb N$). A partir de cette isométrie, A.Ocneanu définit un anti-automorphisme $j_{n}$ de $N'_{0} \cap N_{2n}$  en posant :
$$ j_{n}(x)=J_{n}x^*J_{n}\qquad( x\in N'_{0} \cap N_{2n}).$$
L'anti-automorphisme $j_{n}$ envoie $N'_{0} \cap N_{n}$ sur $N'_{n} \cap N_{2n}$.

\subsubsection{} On rappelle ici les principales propriétés de ces 
anti-automorphismes.
\begin{thm}[Da1- 2.2.1,2.2.2]\label{thj}
 Pour tout entier naturel $n$, les 
anti-automorphismes $j_{n}$ sont involutifs et satisfont les relations suivantes :

a) La restriction de $j_{n+2}j_{n+1}$ à $N'_{0} \cap N_{2n}$ 
	coïncide avec $j_{n+1}j_{n}$.
	
b) Si on note $ F_{n}$ le projecteur de Jones de l'inclusion $N_{0} 
\subset N_{n}$, pour tout $x$ de $N'_{0} \cap N_{n}$, on a l'égalité 
$F_{n}x = F_{n}j_{n}(x)$.

c) $j_{n}(f_{p})=f_{2n-p}\quad (1\leq p \leq n)$.

d) Si l'inclusion $N_0 \subset  N_1$ est de profondeur finie, l'anti-automorphisme $j_{n}$
conserve la trace de $N'_{0} \cap N_{2n}$ pour tout entier $n$.
\end{thm}

\subsubsection{}\label{gamma} Comme les isomorphismes $j_{n+1}j_{n}$ se prolongent 
les uns les autres, on peut définir un isomorphisme $\gamma$ de la tour dérivée 
par $\gamma_{/N'_{0} \cap N_{2n}}=j_{n+1}j_{n}.$
Cet isomorphisme $\gamma$ ajoute $2$ aux   indices, par exemple :
$$\gamma(f_{n})=f_{n+2}, \quad \gamma(N'_{0}\cap N_{2n})=N'_{2}\cap 
N_{2n+2} \dots $$

\subsubsection{Remarque} \label{remj}
	D'après [Da2. 2.2.3.v et 2.3], on peut affirmer que les 
	applications $S_{A}$ et $S_{B}$ utilisées par D. Nikshych et L. Vainerman 
	dans [NV 3-8.2  68]	sont  respectivement $j_{1}$ et $j_{2}$, nous les noterons ainsi, 
	gardant les notations $S_{A}$ et $S_{B}$ pour les antipodes. On a 
	donc  pour tous $a$ de $N'_{0} \cap N_{2}$ et $b$ de $N'_{1} \cap N_{3}$ :
	$$E_{N'_{1}}(j_{1}(a)f_{2}f_{1})= E_{N'_{1}}(f_{1}f_{2}a)\quad \text{et}\quad
	E_{N_{2}}(bf_{1}f_{2})=E_{N_{2}}(f_{2}f_{1}j_{2}(b))$$
	
	Nous utiliserons sans cesse et pour différentes constructions de base le résultat suivant de [PiPo1] :
	$$\forall x \in N_2 \qquad xf_1=\tau^{-1}E_{N_1}(xf_1)f_1$$

\subsection{L'opérateur $h$}\label{H1}
L'élément qui va rendre compte du fait que l'inclusion n'est pas 
irréductible est l'indice de 
la restriction à $N'_{1} \cap N_{2}$ de la trace $tr$  défini par 
Watatani [W]. Cet opérateur noté $H$ est un élément 
inversible et autoadjoint  du 
centre de $N'_{1} \cap N_{2}$.
Si l'algèbre $N'_{1} \cap N_{2}$ se décompose sur son centre comme 
$\oplus_{j\in J}\; M_{\nu_{j}}(\mathbb{C})q_{j}$ et que, pour tout $j$ 
de $J$, $t_{j}$ soit 
la valeur de la restriction à $N'_{1} \cap N_{2}$ de la trace $tr$ sur 
les projecteurs minimaux de $M_{\nu_{j}}(\mathbb{C})q_{j}$, alors $H$ 
est donné par la formule :
$$H=\sum_{j\in J} t_{j}^{-1}\nu_{j}q_{j}.$$
On vérifie facilement que la trace de $H$ est égale à la dimension 
de $N'_{1}\cap N_{2}$.
On notera  $h$ l'opérateur $H^{1/2}$.

 \begin{lem} Pour tout élément $x$ de $N'_{1} \cap N_{2}$, en 
 particulier pour $x=h$, on a :
 
 (a) $f_{2}j_{2}(x)=f_{2}x$ et $f_{1}x=f_{1}j_{1}(x)$

(b) $j_{2}(x)f_{1}f_{2}=f_{1}f_{2}j_{1}(x)$
\end{lem}

\dem

 (a)  D'après \ref{thj} (c).

(b) On utilise (a) et les propriétés 
de commutation des projecteurs $f_{1}$ et $f_{2}$.
\end{proof}

 \subsection{Unités matricielles et quasi-bases}

\subsubsection{Notations} Les notations et les résultats utilisés ici 
se trouvent dans [GHJ] aux paragraphes 2.3.11, 2.4.6 et 2.6.5. On peut 
voir en \ref{chemins} un exemple.
Soit une inclusion $K \subset L$ de $C^*$-algèbres unitaires de dimension 
finie :
$$K = \oplus_{j\in J} M_{\nu_{j}}(\mathbb{C})q_{j}  \subset  L= \oplus_{i\in I} 
M_{\mu_{i}}(\mathbb{C})p_{i} $$
On suppose qu'il existe sur $L$ une trace fidèle $tr$ et on note  
  $\vec t$ et $\vec s$ les vecteurs
définis par :
$$t_{j}=\nu^{-1}_{j}tr(q_{j})\;\;(j\in J), \qquad 
s_{i}=\mu^{-1}_{i}tr(p_{i})\;\;(i\in I).$$
 et $h$ l'élément $\sum_{j\in J}\sqrt{\nu_{j}t_{j}^{-1}}q_{j}$ de $K$.
 
 On considère le diagramme de Bratelli (augmenté d'un sommet *) de $K 
 \subset L$. On appelle étage $1$ du diagramme celui  des sommets 
 représentant les facteurs de $K$, étage $2$ l'étage de ceux de $L$.
 Soit $P$ l'ensemble des couples de chemins $p=(\xi, \eta)$ joignant 
le sommet * à un même sommet du deuxième étage noté 
$\en(p)$ et on note $s_{p}$ la valeur $s_{\en(p)}$. 
L'ensemble des opérateurs 
$\{T_{p}, p\in P\}$ est une famille d'unités matricielles de $L$.

\subsubsection{Quasi-base d'une espérance conditionnelle}\label{qb}
\begin{defn}[W - 1.2.2]

Soit $E$ une espérance conditionnelle fidèle de $L$ sur $K$. Une 
famille finie $\{(u_{1}, v_{1}) \dots (u_{n},v_{n})\}$ de $L\times L$ est dite une 
quasi-base si pour tout $x$ de $L$, on a :
$$\sum_{i=1}^nu_{i}E(v_{i}x)=x=\sum_{i=1}^nE(xu_{i})v_{i}$$
\end{defn}

\begin{prop}[W - 2.4.1] Soit $E$ l'espérance conditionnelle définie 
par la trace $tr$ de $L$ sur $K$. Posons, pour $p \in P$,
$u_{p}=\frac{1}{\sqrt{s_{p}}}T_{p}h^{-1}$ 
alors $\{(u_{p},u_{p}^{*}), p \in P \}$ est une quasi-base pour $E$.	
\end{prop}
On dira dans ce cas que $\{u_{p}, p \in P \}$ est une quasi-base pour 
$E$.

\subsubsection{Cas particulier de la tour dérivée}\label{umn}
On considère le diagramme de Bratelli de la tour 
dérivée de l'inclusion $N_{1}\subset N_{2}$ :
$$N'_{1}\cap N_{1}=\mathbb{C} \subset N'_{1}\cap N_{2} = \oplus_{j\in J} 
M_{\nu_{j}}(\mathbb{C})q_{j}  \subset  N'_{1}\cap N_{3}= \oplus_{i\in I} 
M_{\mu_{i}}(\mathbb{C})p_{i} 
\subset N'_{1}\cap N_{4} \dots$$
La restriction de la trace $tr$ aux algèbres de la tour dérivée est 
une trace de Markov caractérisée par les vecteurs $\vec t$ et $\vec s$ 
définis par :
$$t_{j}=\nu^{-1}_{j}tr(q_{j})\;\;(j\in J), \qquad 
s_{i}=\mu^{-1}_{i}tr(p_{i})\;\;(i\in I).$$

Dans le cas d'une inclusion irréductible de profondeur 2, une famille d'unités 
matricielles normalisées (par 
$tr(b^*_{p}b_{p})=1$) de $B$  est une base orthonormale de $B$ pour le 
produit scalaire issu de la trace $tr$ mais aussi une base de Pimsner-Popa de 
$N_{3}$ sur $N_{2}$. Dans le cas d'une inclusion réductible de profondeur 2, 
ces deux propriétés ne coïncident plus. Si la normalisation est modifiée par $h^{-1}$, nous 
obtenons des quasi-bases.
\begin{prop}\ 

\begin{itemize}
	\item [(i)] 
L'ensemble  $\{b_{p}=\frac{1}{\sqrt{s_{p}}}T_{p}, p\in P\}$ est une 
famille d'unités matricielles normalisées (ou une base orthonormale) de $B$. 
	\item  [(ii)] L'ensemble $\{b_{p}h^{-1}, p\in P\}$ est une 
quasi-base de $N_{3}\cap N'_{1}$ sur $N_{2}\cap N'_{1}$.
	\item  [(iii)] Si l'inclusion $N_{1} \subset N_{2}$ est de 
	profondeur $2$, l'ensemble $\{b_{p}h^{-1}, p\in P\}$ est une 
quasi-base de $N_{3}$ sur $N_{2}$, c'est-à-dire que pour tout $x$ de 
$N_{3}$, on a :
$$x=\sum_{p\in P} b_{p}h^{-1}E_{N_{2}}(h^{-1}b^{*}_{p}x)=
\sum_{p\in P} E_{N_{2}}( xb_{p}h^{-1})h^{-1}b^{*}_{p}.$$
\end{itemize}
\end{prop}
\dem

La première affirmation est évidente puisque la trace de 
$T^{*}_{p'}T_{p}$ vaut $\delta(p,p')s_{p}$ pour $(p,p')$ dans $P\times 
P$. La deuxième résulte du lemme 
2.4.1 de [W] rappelé en \ref{qb}. 

Si l'inclusion $N_{1} \subset N_{2}$ est de 
profondeur $2$, le commutant relatif $N_{4}\cap N'_{1}$ est obtenu 
par contruction de base à partir de l'inclusion 
$N'_{1}\cap N_{2} \subset  N'_{1}\cap N_{3}$ donc il existe une 
famille finie $\{(u_{\mu}, v_{\mu}), \mu \in M\} $ de couples de $N'_{1}\cap 
N_{3}$ telle que $\sum_{\mu \in M} u_{\mu}f_{3}v_{\mu}=1$. Alors  $\{(u_{\mu}, v_{\mu}), \mu \in 
M\} $ est une quasi-base  de $N_{3}$ sur $N_{2}$ en effet $\sum_{\mu \in M} u_{\mu}f_{3}v_{\mu}=1$
implique $\sum_{\mu \in M} u_{\mu}E_{N_{2}}(v_{\mu}x)f_{3}=xf_{3}$ 
pour tout $x$ de $N_{3}$ et donc grâce à [GHJ - 2.6.7(iii)], on peut écrire :
$$\sum_{\mu \in M} u_{\mu}E_{N_{2}}(v_{\mu}x)=x$$
En utilisant (ii) pour chaque $u_{\mu}$, on obtient :
\begin{align*}
x&=\sum_{\mu \in M,p\in P} 
b_{p}h^{-1}E_{N_{2}}(h^{-1}b^{*}_{p} u_{\mu})E_{N_{2}}(v_{\mu}x)\\
&=\sum_{\mu \in M,p\in P} 
b_{p}h^{-1}E_{N_{2}}(h^{-1}b^{*}_{p} u_{\mu}E_{N_{2}}(v_{\mu}x))\\
&=\sum_{p\in P} b_{p}h^{-1}E_{N_{2}}(h^{-1}b^{*}_{p}x)\end{align*}
\end{proof}

\subsubsection{Remarque} Dans le cas d'une inclusion de profondeur 
$2$, on a donc une quasi-base formée d'éléments du commutant relatif 
$N_{3}\cap N'_{1}$, de plus on peut la choisir très proche d'une base 
orthonormale de cette algèbre. Nous verrons dans les calculs que ces 
propriétés sont très précieuses.

\subsection{Dualité}\label{dualite}
Contrairement à D. Nikshych et L. Vainerman, nous définissons une 
dualité  entre $A$ et $B$ en utilisant une formule symétrique :
$$\langle a,b\rangle = \tau^{-2} tr(ahf_{2}f_{1}hb)\quad (a\in A, \;b\in 
B)$$
On remarque que, grâce \ref{H1}, la dualité s'écrit aussi :
$$\langle a,b\rangle = \tau^{-2} 
tr(aj_{2}(h)f_{2}f_{1}j_{1}(h)b)\quad (a\in A, \;b\in 
B)$$
et on a pour tout $a$ de $A$, tout $b$ de  $B$ et tout $x$ de $N'_{1} \cap N_{2}$. :
$$
\langle a,bx \rangle = \langle xa,b \rangle
\qquad \text{et} \qquad
 \langle a,j_{2}(x)b \rangle = \langle aj_{1}(x),b \rangle $$

\noindent {\it Remarque :} Ce crochet définit une dualité car 
l'inclusion est de profondeur $2$ (voir [Szy] ou [NV1 -3.2]).
	
\subsection{Les co-algèbres $A$ et $B$ } 
Les algèbres $A$ et $B$ sont des C*-algèbres. 
Grâce à la dualité entre $A$ et $B$, nous définissons sur $A$ et $B$ 
des structures de co-algèbres co-associatives. Nous allons voir 
que les co-produits et les co-unités de ces co-algèbres sont définis 
par des formules analogues pour  $A$ et $B$ et montrer qu'ils 
vérifient les propriétés  (i) et (ii) des C* groupoïdes quantiques 
(\ref {def}).

\subsubsection{Formule pour les co-produits}\label{formcop}

Le co-produit $\Delta_{B}$ de $B$ est défini comme dual de la 
multiplication de l'algèbre $A$. De même pour le co-produit $\Delta_{A}$ 
de $A$.
\begin{prop} Soient $\{b_{p},p \in P\}$ une famille d'unités 
matricielles normalisées de $B$ et $\{a_{s},s \in S\}$ une famille d'unités 
matricielles normalisées de $A$.

Premières formules :

\noindent Le co-produit $\Delta_{B}$ de $B$ est donn\'e pour $x$ dans $B$ par :
$$\Delta_{B}(x)=\tau ^{-2}\sum_{p\in P} 
E_{N_{3}}(f_{3}xE_{N'_{2}}(b^{*}_{p}f_{3}h^{-1}f_{2}h^{-1}))\otimes b_{p}$$
et le co-produit $\Delta_{A}$ de $A$ vérifie pour $x$ dans $A$ une 
formule analogue :
$$\Delta_{A}(x)=\tau ^{-2}\sum_{s\in S} 
E_{N_{2}}(f_{2}xE_{N'_{1}}(a^{*}_{s}f_{2}j_{1}(h^{-1})f_{1}j_{1}(h^{-1}))) \otimes a_{s}.$$
formule qui s'écrit aussi :
$$\Delta_{A}(x)=\tau ^{-2}\sum_{s\in S} 
E_{N_{2}}(f_{2}xE_{N'_{1}}(a^{*}_{s}f_{2}h^{-1}f_{1}h^{-1})) \otimes a_{s}.$$

Deuxième  formule pour $\Delta_{B}$ :

Si $\{\alpha_{r},r \in R\}$ est une base de Pimsner-Popa de $A$ sur 
$N'_{0}\cap N_{1}$, on a aussi :
$$\Delta_{B}(x)= \sum_{p\in P} 
\sum_{r\in R}E_{A}(x\alpha_{r}b^{*}_{p})h^{-1}f_{2}h^{-1}\alpha^*_{r} \otimes b_{p}$$

\end{prop}

\dem

Soient $a_{1}$ et $a_{2}$ deux éléments de $A$, on a :
\begin{align*}
\langle a_{1}\otimes a_{2},\Delta_{B}(x) \rangle
&=\tau 
^{-6}\sum_{p\in P}tr(a_{1}hf_{2}f_{1}hf_{3}xE_{N'_{2}}(b^{*}_{p}h^{-1}f_{3}f_{2}h^{-1}))
tr(a_{2}hf_{2}f_{1}hb_{p})\\
&=\tau 
 ^{-6}tr(f_{1}hf_{3}xa_{1}E_{N'_{2}}\left(\sum_{p\in 
 P}tr(a_{2}hf_{2}f_{1}hb_{p})b^{*}_{p}h^{-1}f_{3}f_{2}\right)f_{2})\\
 &=\tau  
 ^{-6}tr(f_{1}hf_{3}xa_{1}E_{N'_{2}}[E_{N'_{1}}(a_{2}hf_{2}f_{1})f_{3}f_{2}]f_{2})
\end{align*}

Comme les tours $N'_3 \,\subset \, N'_2 \, \raisebox{1.3 ex}{$\begin{matrix} 
f_{2}\\ \subset \end{matrix}$} \, 
N'_1 $ et $N'_2 \,\subset \, N'_1 \, \raisebox{1.3 ex}{$\begin{matrix} 
f_{1}\\ \subset \end{matrix}$} \, 
N'_0 $  sont standard, on obtient :
\begin{align*}
\langle a_{1}\otimes a_{2},\Delta_{B}(x) \rangle
&=\tau 
^{-5}tr(f_{1}hxa_{1}E_{N'_{1}}(a_{2}hf_{2}f_{1})f_{3}f_{2}f_{3})\\
&=\tau 
^{-4}tr(f_{1}hxa_{1}E_{N'_{1}}(a_{2}hf_{2}f_{1})f_{3})\\
&=\tau 
^{-3}tr(a_{1}E_{N'_{1}}(a_{2}hf_{2}f_{1})f_{1}hx)\\
&=\tau 
^{-2}tr(a_{1}a_{2}hf_{2}f_{1}hx)\\
&=\langle a_{1} a_{2},x \rangle
\end{align*}

La première formule pour $\Delta_{A}$ se montre de manière analogue.

Nous pouvons établir une expression du co-produit $\Delta_{B}$ semblable 
à celle de 
[Da1-5.3.1], c'est la deuxième formule.
Soit $\{\alpha_{r},r \in R\}$ une base de Pimsner-Popa de $A$ sur 
$N'_{0}\cap N_{1}$, on vérifie facilement que $\{\alpha_{r},r \in R\}$ est 
une base de Pimsner-Popa de $N_{2}$ sur 
$N_{1}$. On sait d'après [Bi 2.7] que pour tout $y$ de 
$N'_{1}\cap N_{4}$, on a :
$$\tau \sum_{r\in R} \alpha_{r}y\alpha^{*}_{r}=E_{N'_{2}}(y).$$
En appliquant ce résultat à la première formule de $\Delta_{B}$ , on 
obtient la seconde.
\end{proof}

On utilisera la notation habituelle de Sweedler : 
$\Delta(x)=x_{(1)}\otimes x_{(2)}$.

\subsubsection{}\label{coppa}
Des définitions des co-produits par 
dualité et des propriétés de la dualité, on déduit : 
\begin{prop}
Soit $x$ un élément de $N'_{1} \cap N_{2}$. Pour tous $a$ de $A$ et 
$b $ de $B$, on a :
\begin{align*}
&\Delta_{B}(bx)= \Delta_{B}(b)(x\otimes1)
&\Delta_{B}(j_{2}(x)b) = (1\otimes j_{2}(x))\Delta_{B}(b)\\
&\Delta_{A}(aj_{1}(x)) = \Delta_{A}(a) (j_{1}(x) \otimes 1) 
&\Delta_{A}(xa)=(1 \otimes  x)\Delta_{A}(a)
\end{align*}
\end{prop}

\subsubsection{Les projecteurs $\Delta_{A}(1)$ et 
$\Delta_{B}(1)$}\label{delta(1)}
Nous calculons maintenant les images de l'unité par les co-produits et 
nous faisons le lien avec les projecteurs relatifs aux produits fibrés 
en dimension finie définis par Jean-Michel Vallin dans [V1].

\begin{prop} Si $\{\lambda_{k},k \in K\}$ est une famille d'unités 
matricielles  de $N'_{1} \cap N_{2}$ telle que l'élément 
$\lambda_{k}$ appartienne au facteur $M_{\nu_{j_{k}}}(\mathbb 
C)q_{j_{k}}$, 
l'élément $\Delta_{B}(1)$ est le  projecteur 
$ {\displaystyle \sum_{k}\frac{1}{\nu_{j_{k}}} j_{2}(\lambda^{*}_{k}) \otimes 
\lambda_{k}}$ 
de $N'_{2} \cap N_{3}\otimes N'_{1} \cap N_{2}$.
 \end{prop}

\dem
 
Ecrivons la deuxième formule du co-produit pour $x=1$, comme 
$\sum_{r\in R} \alpha_{r}f_{2}\alpha^{*}_{r}$ vaut $1$, la propriété 
\ref{H1} (a)  donne :
$$\Delta_{B}(1) = \sum_{p} j_{2}(h^{-1}E_{N_{2}}(b^{*}_{p})h^{-1}) \otimes  
b_{p}.$$
Si $\{\lambda_{k},k \in K\}$ est une famille d'unités 
matricielles  de $N'_{1} \cap N_{2}$, 
$\{\sqrt{\frac{1}{t_{j_{k}}}}\lambda_{k},t \in I\}$ est une base 
orthonormale de $N'_{1} \cap N_{2}$
et comme $E_{N_{2}}(b^{*}_{p})h^{-1}$ appartient à $N'_{1} \cap N_{2}$, on a :
\begin{align*}
\Delta_{B}(1) &= \sum_{p,k}\frac{1}{t_{j_{k}}} j_{2}(h^{-1}\lambda^{*}_{k}) 
\otimes tr(E_{N_{2}}(b^{*}_{p})h^{-1}\lambda_{k})
b_{p}\\
&= \sum_{k}\frac{1}{t_{j_{k}}} j_{2}(h^{-1}\lambda^{*}_{k}) \otimes h^{-1}\lambda_{k}
\end{align*}
On obtient donc la formule annoncée et l'appartenance de 
$\Delta_{B}(1)$ à $N'_{2} \cap N_{3}\otimes N'_{1} \cap N_{2}$~;
 $\Delta_{B}(1)$ est un projecteur car $\Delta_{B}$ est un homomorphisme 
 d'algèbres, comme on le montre 
 au paragraphe suivant.
\end{proof}
\indent Considérons l'identité comme représentation de $N'_{1} \cap N_{2}$ dans 
$A$ et $B$, $j_{1}$ (resp. $j_{2}$) comme antireprésentation de $N'_{1} \cap N_{2}$ dans 
$A$ (resp. $B$).
Alors $\Delta_{A}(1)$ (resp. $\Delta_{B}(1)$) est le projecteur 
$e_{Id,j_{1}}$ (resp. $e_{j_{2},Id}$) défini par Jean-Michel Vallin [V1].

\subsubsection{Propriétés des co-produits}\label{cophom}
\begin{prop}
Les co-produits sont des homomorphismes d'algèbres involutives \\
vérifiant :
$$(\Delta \otimes \id)\Delta(1)=(1\otimes \Delta(1))(\Delta(1)\otimes 1)$$
\end{prop}

\dem

Montrons la proposition pour $\Delta_{B}$ par exemple. 
D'après \ref{formcop}, pour tous $x$ et $y$ de $B$, on a :
\begin{align*}
\Delta_{B}(x)\Delta_{B}(y)
&=\tau ^{-2}\sum_{p,p'\in P} 
E_{N_{3}}(f_{3}xE_{N'_{2}}(b^{*}_{p}f_{3}h^{-1}f_{2}h^{-1}))
\sum_{r\in R}E_{A}(y\alpha_{r}b^{*}_{p'})h^{-1}f_{2}h^{-1}\alpha^*_{r}
\otimes b_{p} b_{p'}\\
&=\tau ^{-2}\sum_{p,p'\in P}\sum_{r\in R} 
E_{N_{3}}[f_{3}xE_{A}(y\alpha_{r}b^{*}_{p'})h^{-1}E_{N'_{2}}(b^{*}_{p}f_{3}h^{-1}f_{2}j_2(h^{-1}))f_{2}h^{-1}\alpha^*_{r}]
\otimes b_{p} b_{p'}\\
&=\tau ^{-2}\sum_{p,p'\in P}\sum_{r\in R} 
E_{N_{3}}[f_{3}xE_{A}(y\alpha_{r}b^{*}_{p'})H^{-1}E_{N'_{2}}(b^{*}_{p}f_{3}h^{-1}f_{2})f_{2}h^{-1}\alpha^*_{r}]
\otimes b_{p} b_{p'}\\
&=\tau ^{-1}\sum_{p,p'\in P}\sum_{r\in R} 
E_{N_{3}}[f_{3}xE_{A}(y\alpha_{r}b^{*}_{p'})H^{-1}b^{*}_{p}f_{3}h^{-1}f_{2}h^{-1}\alpha^*_{r}] 
\otimes b_{p} b_{p'}\\
&=\sum_{p,p'\in P}\sum_{r\in R}  
E_{N_{2}}[xE_{A}(y\alpha_{r}b^{*}_{p'})H^{-1}b^{*}_{p}]h^{-1}f_{2}h^{-1}\alpha^*_{r} 
\otimes b_{p} b_{p'}\\
&=\sum_{p,p',q\in P}\sum_{r\in R}  
E_{N_{2}}[xE_{A}(y\alpha_{r}b^{*}_{p'}tr(b_{p'}b^{*}_{q}b_{p} ))H^{-1}b^{*}_{p}]h^{-1}f_{2}h^{-1}\alpha^*_{r} 
\otimes b_{q} \\
&=\sum_{p,q\in P}\sum_{r\in R}  
E_{N_{2}}[xE_{A}(y\alpha_{r}b^{*}_{q}b_{p} h^{-1})h^{-1}b^{*}_{p}]h^{-1}f_{2}h^{-1}\alpha^*_{r}
\otimes b_{q} \\
&=\sum_{p,q\in P}\sum_{r\in R}  
E_{N_{2}}(xy\alpha_{r}b^{*}_{q})h^{-1}f_{2}h^{-1}\alpha^*_{r} 
\otimes b_{q}\quad \text{(d'après \ref{umn}})\\
&=\Delta_{B}(xy)
\end{align*}

Traitons maintenant le cas de l'involution :
\begin{align*}
\Delta_{B}(x)^*&= \sum_{p\in P} 
\sum_{r\in R} \alpha_{r}h^{-1}f_{2}h^{-1}E_{A}(b_{p}\alpha^{*}_{r}x^{*}) 
\otimes b^{*}_{p}\\
&= \sum_{p\in P} 
\sum_{r,s\in R} \alpha_{r}h^{-1}f_{2}h^{-1}E_{N_{1}}(b_{p}\alpha^{*}_{r}x^{*}\alpha_{s})\alpha^{*}_{s} 
\otimes b^{*}_{p}\\
&= \sum_{p\in P} 
\sum_{r,s\in R} \alpha_{r}E_{N_{1}}(b_{p}\alpha^{*}_{r}x^{*}\alpha_{s})h^{-1}f_{2}h^{-1}\alpha^{*}_{s} 
\otimes b^{*}_{p}
 \end{align*}
 Comme $b_{p}$ commute avec $N_{1}$, on a :
 $$E_{N_{1}}(b_{p}\alpha^{*}_{r}x^{*}\alpha_{s})
 =E_{N_{1}}(\alpha^{*}_{r}x^{*}\alpha_{s}b_{p})$$
Puisque $\{\alpha_{r},r \in R\}$ est une base de Pimsner-Popa de $A$ sur 
$N'_{0}\cap N_{1}$ et que \hfill \break $\{b^{*}_{p}, p\in P\}$ est une famille 
d'unités matricielles normalisées, on obtient  :
\begin{align*}
\Delta_{B}(x)^*&=\sum_{p\in P} 
\sum_{s\in R}E_{A}(x^{*}\alpha_{s}b_{p})h^{-1}f_{2}h^{-1}\alpha^{*}_{s} 
\otimes b^{*}_{p}
= \Delta_{B}(x^{*})
\end{align*}
Montrons la relation pour $\Delta_{B}(1)$. 
Comme $1_{(1)}$ appartient à $N'_{2}\cap N_{3}$, d'après \ref{coppa}, 
on a :
$$(\Delta_{B} \otimes \id)\Delta_{B}(1)= \Delta_{B}(1_{(1)})\otimes 
1_{(2)}=(1\otimes 1_{(1)})\Delta_{B}(1)\otimes 1_{(2)}=(1\otimes 
\Delta_{B}(1))(\Delta_{B}(1)\otimes 1)$$
\end{proof}
\subsubsection{Etude des co-unités}\label{coun}
\begin{prop} Pour tout élément $b$ de $B$, 
la co-unité $\varepsilon_{B}$ de $B$ est donn\'ee par : 
$$\varepsilon_{B}(b)=\tau^{-1}tr(hf_{2}hb).$$
Pour tout élément $a$ de $A$, 
la co-unité $\varepsilon_{A}$ de $A$, vérifie une formule analogue :
$$\varepsilon_{A}(a)=\tau^{-1}tr(hf_{1}ha)=\tau^{-1}tr(j_{1}(h)f_{1}j_{1}(h)a).$$
Chaque co-unité vérifie :
$\varepsilon(xyz) =\varepsilon(xy_{(2)})\varepsilon(y_{(1)}z)$.
\end{prop}
\dem

Les formules sont évidentes compte-tenu des propriétés des 
projecteurs de Jones.
Montrons la relation pour la co-unité de $B$.
Soient $x$, $y$ et $z$ trois éléments de $B$. On pose :
$$\Delta_{A}(1)=1_{(1)}\otimes 1_{(2)} \quad \text{et} \quad 
\Delta_{B}(y)=y_{(1)}\otimes y_{(2)}.$$
On peut alors écrire :
\begin{align*}
\varepsilon(xy_{(2)})\varepsilon(y_{(1)}z)&=\langle 1_{(1)},x\rangle \langle 
1_{(2)},y_{(2)} 
\rangle  \langle1_{(1)},y_{(1)}\rangle \langle 1_{(2)},z \rangle\\
&=\langle 1_{(1)},x\rangle \langle1_{(2)},\varepsilon_{B}(y_{(1)}1_{(1)}) y_{(2)}\rangle \langle 1_{(2)},z \rangle
\end{align*}
Or comme $1_{(1)}$ est un élément de $N'_{1} \cap N_{2}$ et $B$ une 
co-algèbre, d'après 
\ref{coppa}, on a :
$$\varepsilon_{B}(y_{(1)}1_{(1)}) y_{(2)}=(\varepsilon_{B}\otimes 
\id)(\Delta_{B}(y)(1_{(1)}\otimes 1))=(\varepsilon_{B}\otimes 
\id)\Delta_{B}(y1_{(1)})= y1_{(1)}$$

On obtient donc
\begin{align*}
\varepsilon(xy_{(2)})\varepsilon(y_{(1)}z)
&=\langle 1_{(1)},x\rangle \langle1_{(2)},y1_{(1)}\rangle \langle 1_{(2)},z \rangle\\
&=\langle 1,xy1_{(1)}\rangle \langle 1_{(2)},z \rangle\\
&=\langle 1_{(1)},xy\rangle \langle 1_{(2)},z \rangle\\
 &=\varepsilon_{B}(xyz).
\end{align*}
\end{proof}

Nous donnons maintenant les expressions des co-unités but et source 
de $B$, les formules pour $A$ sont analogues :

\begin{prop} Pour tout $x$ de $B$, on a :
\begin{align*}
\varepsilon^t_{B}(x)&=\tau^{-1}E_{N'_{1}\cap N_{2}}(xhf_{2}h^{-1}) \\
\varepsilon^s_{B}(x)&=\tau^{-1}j_{2}(E_{N'_{1}\cap N_{2}}(xhf_{2}h^{-1}))
=\tau^{-1}E_{N'_{2}\cap N_{3}}(j_{2}(x)h^{-1}f_{2}h)
\end{align*}

La sous-algèbre co-unitale but $A_{t}$ de $A$ est $N'_{0} \cap N_{1}$, la 
sous-algèbre co-unitale source $B_{s}$ de $B$ est $N'_{2} \cap 
N_{3}$, les sous-algèbres co-unitales but $B_{t}$ et source $A_{s}$ 
coïncident avec $N'_{1}\cap N_{2}$.
\end{prop}

\dem

Avec les notations de \ref{delta(1)} et grâce à la proposition 
précédente, on a : 

\begin{align*}
\varepsilon^t_{B}(x)&=\tau^{-1}\sum_{k}\frac{1}{\nu_{j_{k}}} 
tr(xhf_{2}hj_{2}(\lambda^{*}_{k})) \lambda_{k}\\
&=\tau^{-1}\sum_{k}\frac{1}{\nu_{j_{k}}} 
tr(xhf_{2}\lambda^{*}_{k}h) \lambda_{k}\\
&=\tau^{-1}\sum_{k}
tr(xhf_{2}\frac{1}{\sqrt{t_{j_{k}}}} \lambda^{*}_{k}) 
\frac{1}{\sqrt{t_{j_{k}}}} \lambda_{k}h^{-1}
\end{align*}
On obtient la formule annoncée en remarquant que 
$\{\frac{1}{\sqrt{t_{j_{k}}}}\lambda_{k},k \in K\}$ est une famille d'unités 
matricielles  normalisées de $N'_{1} \cap N_{2}$. La deuxième formule 
se montre de même.
\end{proof}

\subsection{Antipodes sur $A$ et $B$}\label{SB}
Pour obtenir des structures de C*-groupoïdes quantiques duaux sur $A$ et 
$B$, nous complétons nos données par les antipodes $S_{A}$ et $S_{B}$ 
définies pour tous $a$ de $A$ et $b$ de $B$ par :
$$\langle S_{A}(a)^{*},b \rangle=\overline{\langle a,b^{*} \rangle} 
\qquad \langle a,S_{B}(b)^{*} \rangle= \overline{\langle a^{*},b 
\rangle}$$
\subsubsection{Formules pour les antipodes} Nous obtenons alors les formules suivantes :
\begin{align*}
S_{B}(b) &= hj_{2}(h^{-1})j_{2}(b)h^{-1}j_{2}(h)\\
S_{A}(a) &= j_{1}(h)h^{-1}j_{1}(a)j_{1}(h^{-1})h
\end{align*}
Il est évident que $S_{A}$ et $S_{B}$ sont des anti-automorphismes 
d'algèbres conservant l'unité.
On vérifie facilement la formule de dualité : 
$$\langle S_{A}(a),b \rangle=\langle a,S_{B}(b) \rangle.$$
Et par dualité, on obtient que $S_{A}$ et $S_{B}$ sont des anti-automorphismes 
de co-algèbres conservant les co-unités.

\subsubsection{} Il nous reste à vérifier la formule liant l'antipode, la co-unité et 
le co-produit :
\begin{prop} Les antipodes $S_{A}$ et $S_{B}$ vérifient la relation :
$$	m(\id \otimes S)\Delta(x)=(\varepsilon \otimes 
	\id)(\Delta(1)(x \otimes 1))=\varepsilon_{t}(x)$$
\end{prop}
\dem

Soit $x$ un élément de $B$, avec les notations de \ref{formcop}, on 
a :
\begin{align*}
m(\id \otimes S_{B})\Delta_{B}(x)&=x_{(1)}S_{B}(x_{(2)})\\
&=\tau ^{-2}\sum_{p\in P} 
E_{N_{3}}[f_{3}xE_{N'_{2}}(b^{*}_{p}f_{3}h^{-1}f_{2}h^{-1})]hj_{2}(h^{-1})j_{2}(b_{p})h^{-1}j_{2}(h)\\
&=\tau ^{-2}\sum_{p\in P} 
E_{N_{3}}[f_{3}xE_{N'_{2}}(b^{*}_{p}f_{3}h^{-1}f_{2})hj_{2}(b_{p}h^{-2})h^{-1}j_{2}(h)]
\end{align*}
Or, si $\{\alpha_{r},r \in R\}$ est une base de Pimsner-Popa de $A$ sur 
$N'_{0}\cap N_{1}$, la formule de  [Da1-3.2.1] nous donne une expression de 
$j_{2}(b_{p}h^{-2})$ :
$$ j_{2}(b_{p}h^{-2})=\tau^{-1} \sum_{r\in R} 
E_{A}(f_{2}\alpha_{r}b_{p}h^{-2})f_{2}\alpha^{*}_{r}.$$
On peut donc écrire :
\begin{align*}
m(\id \otimes S_{B})\Delta_{B}(x)&=\tau ^{-3}\sum_{p\in P,r\in R} 
E_{N_{3}}[f_{3}xhE_{A}(f_{2}\alpha_{r}b_{p}h^{-2})E_{N'_{2}}(b^{*}_{p}f_{3}h^{-1}f_{2})f_{2}
\alpha^{*}_{r}  h^{-1}j_{2}(h)]\\
&=\tau ^{-2}\sum_{p\in P,r\in R} 
E_{N_{3}}[f_{3}xhE_{A}(f_{2}\alpha_{r}b_{p}h^{-1})h^{-1}b^{*}_{p}f_{3}h^{-1}f_{2}\alpha^{*}_{r}  h^{-1}j_{2}(h)]
\end{align*}
Comme l'ensemble $\{b_{p}h^{-1}, p\in P\}$ est une 
quasi-base de $N_{3}$ sur $N_{2}$ [\ref{umn}], on a :
$$\sum_{p\in P} 
E_{A}(f_{2}\alpha_{r}b_{p}h^{-1})h^{-1}b^{*}_{p}= f_{2}\alpha_{r}$$
d'autre part puisque $\{\alpha_{r},r \in R\}$ est une base de Pimsner-Popa, 
la somme
${\displaystyle \sum_{r\in R}\alpha_{r}f_{2}\alpha^{*}_{r}}$ vaut $1$,
on en déduit :
\begin{align*}
m(\id \otimes S_{B})\Delta_{B}(x)&=\tau ^{-2}\sum_{r\in R} 
E_{N_{3}}(f_{3}xhf_{2}\alpha_{r}f_{3}j_{2}(h^{-1})f_{2}\alpha^{*}_{r}  h^{-1}j_{2}(h))\\
&=\tau ^{-2} 
E_{N_{3}}(f_{3}xhf_{2}f_{3} h^{-1})\\
&=\tau ^{-2} 
E_{N_{3}}(E_{N_{2}}(xhf_{2}h^{-1})f_{3})\\
&=\tau ^{-1} E_{N_{2}}(xhf_{2}h^{-1})
\end{align*}
On obtient la formule annoncée  grâce à l'expression de 
$\varepsilon^t_{B}$ démontrée en \ref{coun}.
\end{proof}

\subsubsection{Remarques} Les formules définissant $S_{A}$ et $S_{B}$ sont 
analogues. Les automorphismes $S_{A}^2$ et $S_{B}^2$ sont 
intérieurs :
$$S_{A}^2 = \Ad (j_{1}(H) H^{-1}) \qquad S_{B}^2 = \Ad (j_{2}(H^{-1}) H)$$
Les C*-groupoïdes quantiques construits sont donc réguliers.

\subsection{Projection de Haar, mesure de Haar} On démontre facilement 
la proposition suivante.
\begin{prop} Soit $d$ la dimension de 
$N'_{1}\cap N_{2}$. La projection de Haar normalisée du C*-groupoïde 
quantique $B$ est $p_{B}=d^{-1}hf_{2}h$. La mesure de Haar normalisée du C*-groupoïde 
quantique $B$ est $\phi_{B}$ définie par :
$$\phi_{B}(x) = d^{-1}tr(Hj_{2}(H) x)\quad (x \in B).$$
On a des formules analogues pour le C*-groupoïde quantique $A$. 
\end{prop}

\subsection{Actions des C*-groupoïdes quantiques}\ 

D'après [NV 1 - 6.1], le C*-groupoïde quantique $B$
 agit à gauche sur l'algèbre $N_{2}$. Nous définissons ici une action 
 à gauche du C*-groupoïde quantique $A$ sur $N_{1}$ puis nous obtiendrons 
 par dualité une  action à gauche du C*-groupoïde quantique $B$ sur 
 $N_{2}$. On trouvera toutes les définitions concernant les actions et 
 les produits croisés dans la partie 2.
 \subsubsection{Action du C*-groupoïde quantique $A$ sur 
 $N_{1}$}\label{act}
 \begin{prop} L'application : $\begin{array}{ccc}A \otimes 
 N_{1}&\rightarrow& N_{1}\\ a \otimes x &\mapsto&
 a \triangleright x=\tau^{-1} 
 E_{N_{1}}(axhf_{1}h^{-1})\end{array}$ est une action 
 à gauche standard et extérieure du C*-groupoïde quantique $A$ sur l'algèbre $N_{1}$.
 \end{prop}
 
 \dem

 Tout au long de cette démonstration, nous utilisons l'égalité 
 : $f_{1}h=f_{1}j_{1}(h)$ ( voir \ref{H1} (a)).
  Montrons d'abord que l'application définit une structure de 
 $A$-module à gauche sur $N_{1}$ :
 \begin{align*}
 1\triangleright x&=\tau^{-1}  E_{N_{1}}(xhf_{1}h^{-1})=x\\
 a\triangleright(c \triangleright x)
 &=\tau^{-2} E_{N_{1}}(a E_{N_{1}}(cxhf_{1}h^{-1})hf_{1}h^{-1})\\
 &=\tau^{-2} E_{N_{1}}(a E_{N_{1}}(cxhf_{1})f_{1}h^{-1})\\
 &=\tau^{-1} E_{N_{1}}(acxhf_{1}h^{-1})=(ac)\triangleright x
 \end{align*}
 
 Nous étudions maintenant l'action de $A$ sur un produit.
 
 \begin{lem} Pour tout $a$ de $A$ et $x$ de $N_{1}$, on a :
 $$(a_{(1)}\triangleright x)a_{(2)}=ax$$
 \end{lem}
 
 \dem
 
 Soit $\{a_{s},s\in S\}$ une famille d'unités matricielles normalisées 
 de $A$. D'après \ref{formcop}, on peut écrire :
 \begin{align*}
  (a_{(1)} \triangleright x)a_{(2)}
  &=\tau ^{-3}\sum_{s\in S} 
E_{N_{1}}[E_{N_{2}}(f_{2}aE_{N'_{1}}(a^{*}_{s}f_{2}h^{-1}f_{1}h^{-1})) 
x hf_{1}h^{-1}]a_{s}\\
  &=\tau ^{-3}\sum_{s\in S} 
E_{N_{1}}[f_{2}a xE_{N'_{1}}(a^{*}_{s}f_{2}h^{-1}f_{1}) f_{1}h^{-1}] a_{s}\\
 &=\tau ^{-2}\sum_{s\in S} 
E_{N_{1}}(f_{2}a xa^{*}_{s}f_{2}h^{-1} f_{1}h^{-1}) 
 a_{s}\\
  &=\tau ^{-2}\sum_{s\in S} 
E_{N_{1}}(f_{2}E_{N_{1}}(a xa^{*}_{s})j_{1}(h^{-1}) 
f_{1}j_{1}(h^{-1}))
 a_{s}\\
 &=\sum_{s\in S} 
E_{N_{1}}(a xa^{*}_{s}j_{1}(h^{-1}))
 j_{1}(h^{-1}) a_{s}\\
  \end{align*}

  On conclut en remarquant que si 
 $\{a_{s},s\in S\}$ une famille d'unités matricielles normalisées 
 de $A$, $\{a^{*}_{s},s\in S\}$ l'est aussi et d'après \ref{umn}, 
 $\{a^{*}_{s}j_{1}(h^{-1}),s\in S\}$ est une 
 quasi-base de $N_{2}$ sur $N_{1}$.
 \end{proof}
 
 On reprend maintenant la démonstration des propriétés de l'action. En 
 appliquant le lemme précédent pour  $a$ dans $A$ et $x$ et $y$ dans 
 $N_{1}$, on obtient :
 $$(a_{(1)}\triangleright x)(a_{(2)}\triangleright y)
 =\tau^{-1}  E_{N_{1}}((a_{(1)}\triangleright 
 x)a_{(2)}yhf_{1}h^{-1})=a\triangleright xy$$
 
 Pour montrer la relation : $(a\triangleright x)^{*}=S_{A}(a)^{*}\triangleright 
 x^{*}$, commençons par un lemme :
 \begin{lem} Si $x$ et $y$ sont des éléments  de $N_{1}$ et que $a$ 
 appartienne à $A$, alors on a :
 $$tr(j_{1}(a)xf_{1}y)= tr(f_{1}xay)$$
 \end{lem}
 \dem
 
 Ici, en considérant l'inclusion des commutants, la formule [NV 1 - 4.5 (i)] s'écrit :
 $$j_{1}(a)=\tau^{-3}E_{N_{2}}(f_{2}f_{1}E_{N'_{1}}(af_{2}f_{1}))$$
 En utilisant les propriétés de commutation, on obtient donc :
  \begin{align*}
 tr(j_{1}(a)xf_{1}y)
 &=\tau^{-3} tr(E_{N_{2}}(f_{2}f_{1}E_{N'_{1}}(af_{2}f_{1}))xf_{1}y)\\
 &=\tau^{-3} tr(f_{2}f_{1}xE_{N'_{1}}(af_{2}f_{1})f_{1}y)\\
 &=\tau^{-2} tr(f_{2}f_{1}xaf_{2}f_{1}y)\\
 &= tr(E_{N_{1}}(f_{1}xa)y)\\
 &=tr(f_{1}xay)
 \end{align*}\end{proof}
 
 Grâce à ce lemme, pour $x$ et $y$  éléments  de $N_{1}$ et  $a$ 
 dans  $A$, on peut écrire :
 \begin{align*}
 tr((S_{A}(a)^{*}\triangleright x^{*})y)
 &=\tau^{-1} tr(hj_{1}(h^{-1})j_{1}(a^{*})h^{-1}j_{1}(h)x^{*}hf_{1}h^{-1}y)\\
  &=\tau^{-1} tr(j_{1}(ha^{*}h^{-1})x^{*}f_{1}y)\\
  &=\tau^{-1} tr(f_{1}x^{*}ha^{*}h^{-1}y)\\
  &=\tau^{-1} tr(h^{-1}f_{1}hx^{*}a^{*}y)\\
  &=tr((a\triangleright x)^{*}y)
 \end{align*}
 
 On étudie ensuite l'action de $A$ sur l'unité de $N_{1}$ :
 $$a\triangleright 1=\tau^{-1} 
 E_{N_{1}}(ahf_{1}h^{-1})=\varepsilon^t_{A}(a)=\varepsilon^t_{A}(a)\triangleright 1$$
 On vérifie facilement que $a\triangleright 1$ est nul si et 
 seulement si $ \varepsilon^t_{A}(a)$ l'est. L'action est donc 
 standard. La proposition \ref{N2} nous permet d'affirmer qu'elle est 
 extérieure (voir \ref{ext}) puisque d'après \ref{coun}, $N'_{1}\cap N_{2}$ est 
 la sous-algèbre co-unitale 
 $A_{s}$ .
 \end{proof} 
 
 \subsubsection{Points fixes sous l'action de $A$} \label{ptsfixes} Par définition, un 
 élément $x$ de $N_{1}$ est un point fixe sous l'action de $A$ si on a 
 l'égalité
 $$a\triangleright x=\varepsilon^t_{A}(a)\triangleright x \qquad (*)$$

\begin{prop} L'algèbre $N_{0}$ est l'algèbre des points fixes de 
$N_{1}$ sous l'action du C*-groupoïde quantique $A$.
\end{prop}

\dem

Si $x$ est dans $N_{0}$, alors $x$ commute à $A$ donc vérifie 
l'égalité (*).
Si $x$ est un point fixe, écrivons l'égalité (*) pour $a=f_{1}h^{-1}$, 
en utilisant le lemme \ref{H1}, nous obtenons alors la suite d'égalités équivalentes :
\begin{align*}
(f_{1}h^{-1})\triangleright x&=\varepsilon^t_{A}(f_{1}h^{-1})\triangleright 
x\\
\tau^{-1} E_{N_{1}}(f_{1}xf_{1}h^{-1})&=
\tau^{-2} E_{N_{1}}(E_{N_{1}}(f_{1}h^{-1})xhf_{1}h^{-1})\\
 E_{N_{0}}(x) j_{1}(h^{-1})&=xj_{1}(h^{-1})\\
E_{N_{0}}(x)&=x
\end{align*}
Donc $x$ appartient à $N_{0}$.
\end{proof}

 \subsubsection{Produit croisé de $N_{1}$ par $A$}\label{N2}
 
\begin{prop}
L'application $\Theta$ : $[x \otimes a] \longmapsto xa$ est un isomorphisme 
d'algèbres de von Neumann entre $N_{1}\rtimes A$ et $N_{2}$.
\end{prop}

\dem

La démonstration est semblable à celle donnée par D. Nikshych et L. Vainerman
 dans [NV 1 - 6.3], nous la donnons pour être complets.
 
Si $z$ appartient à $A_{t}=N'_{0} \cap N_{1}$, on a :
$$z \triangleright 1=\varepsilon^t_{A}(z)=z$$
donc  $\Theta$ définit une application linéaire de $N_{1}\otimes_{A_{t}} A$ 
dans $N_{2}$ qui est surjective puisque $A$ fournit une quasi-base de 
$N_{2}$ sur $N_{1}$.

En utilisant le premier lemme démontré en \ref{act}, on obtient pour $a$ et 
$c$ dans $A$, $x$ et $y$ dans $N_{1}$ :
 \begin{align*}
 \Theta([x \otimes a][y \otimes c])
 &=\Theta([x(a_{(1)}\triangleright y) \otimes a_{(2)}c])\\
 &=x(a_{(1)}\triangleright y)a_{(2)}c\\
 &=xayc\\
&=\Theta([x \otimes a])\Theta([y \otimes c])\\
 \Theta([x \otimes a]^*)&=\Theta([(a^{*}_{(1)}\triangleright 
 x^{*})\otimes a^{*}_{(2)}])\\
 &=(a^{*}_{(1)}\triangleright x^{*}) a^{*}_{(2)}\\
 &=a^{*}x^{*}\\
 &=\Theta([x \otimes a])
 \end{align*}
 L'application  est donc un homomorphisme surjectif d'algèbres involutives. 
Puisque ces algèbres sont des facteurs de type $\mathrm{II}_{1}$, $\Theta$ est 
injectif et la proposition est démontrée.
\end{proof}

\subsubsection{Action du C*-groupoïde quantique $B$ sur 
$N_{2}$}\label{actB}
On précise maintenant l'action à gauche de $B$ sur $A$.

\begin{lem} L'action à gauche de $B$ sur $A$ définie par dualité 
est :
$$b\triangleright a= \tau^{-1} E_{A}(bahf_{2}h^{-1})\qquad(a\in 
A,b\in B)$$
\end{lem}
\dem

D'après \ref{dual}, l'action duale de $B$ sur $A$ est définie par 
$$ b\triangleright a=a_{(1)}\langle a_{(2)},b\rangle  \qquad(a\in A,b\in B).$$
Pour tout $c$ dans $B$, on a donc :
 \begin{align*}
 \langle b\triangleright a,c\rangle
 &=\langle a_{(1)},c\rangle \langle a_{(2)},b\rangle\\
 &=\langle a,cb\rangle \\
 &=\tau^{-2} tr(bahf_{2}f_{1}hc)\\
  &=\tau^{-2} tr(\tau^{-1}E_{N_{2}}(bahf_{2})f_{2}f_{1}hc)\\
  &=\tau^{-2} tr(\tau^{-1}E_{N_{2}}(bahf_{2}h^{-1})hf_{2}f_{1}hc)\\
  &=\langle \tau^{-1}E_{A}(bahf_{2}h^{-1}),c\rangle
 \end{align*}
\end{proof}

\begin{prop} L'application : $\begin{array}{ccc}B \otimes N_{2}&\rightarrow& 
N_{2}\\ b \otimes x &\mapsto&b \triangleright x=\tau^{-1} 
 E_{N_{2}}(bxhf_{2}h^{-1})
 \end{array}$ 
 est une action 
 à gauche du C*-groupoïde quantique $B$ sur l'algèbre $N_{2}$ qui 
 prolonge l'action duale de $B$ sur $A$ et  dont l'algèbre de points 
fixes est $N_{1}$.
Le facteur $N_{3}$ est isomorphe au produit croisé de 
$N_{2}$ par  $B$.
\end{prop}
\dem

Il suffit de remarquer que la formule est analogue à celle qui définit 
l'action de $A$ sur $N_{1}$ et que si $x$ appartient à $A$, on 
retrouve l'action duale. C'est l'action duale de $B$ sur $N_{1 
}\rtimes A$ (voir \ref{dualprod}).
Comme l'algèbre $N_{2}$ est linéairement engendrée par les produits 
$xa$ ($x\in N_{1}, a\in A$), on a plus précisément  pour $x$ dans 
$N_{1}$ et $a$ dans $A$ l'\'egalit\'e $b\triangleright xa= x(b\triangleright a).$
\end{proof}

\subsection{Conclusion} Nous avons donc défini sur les  
commutants relatifs $N'_{0}\cap N_{2}$ et $N'_{1}\cap N_{3}$ des 
structures duales de C*-groupoïde quantique données par des 
formules analogues. Ces C*-groupoïdes quantiques agissent 
extérieurement de manière 
analogue et duale sur les facteurs $N_{1}$ et $N_{2}$. L'intérêt de la symétrie 
des définitions apparaîtra dans les parties suivantes.

Les C*-groupoïdes quantiques associés à une inclusion de profondeur 
$2$ de facteurs de type $II_1$ sont un peu particuliers (en 7, nous montrons que la régularité est possible par déformation), d'une part 
ils sont réguliers, 
d'autre part ils sont connexes en effet :
$$A_t \cap Z(A) \subset N'_{0}\cap N_{1} \cap \{f_{1}\}' = N'_{0}\cap 
N_{0}= \mathbb{C}$$
Cela tient au fait qu'on considère des inclusions de facteurs.

\section{C*-groupoïdes quantiques associés à une inclusion d'indice 
fini de 
profondeur finie de facteurs de type $\mathrm{II}_{1}$}

\subsection{Inclusion de profondeur finie}\label{proffinie}
  Soit $P_{0} \subset P_{1}$ 
une inclusion d'indice fini $\delta^{-1}$ et de profondeur finie de facteurs de 
type $\mathrm{II}_{1}$. On note 
$$P_0 \,\subset \, P_1 \, \raisebox{1.3 ex}{$\begin{matrix} e_{1}\\ \subset \end{matrix}$} \, 
P_2 \,\raisebox{1.3 ex}{$\begin{matrix} e_{2}\\ \subset \end{matrix}$} 
\,  P_3 \subset \dots P_{n}\, \raisebox{1.3
ex}{$\begin{matrix} e_{n}\\ \subset \end{matrix}$} \, P_{n+1}\dots $$
 la tour de Jones obtenue par construction de base [G.H.J. 3] 
et  $tr$ la trace normale finie normalisée sur les 
facteurs considérés.

\subsubsection{}\label{prof2} La proposition suivante montre 
qu'une inclusion de profondeur finie peut être vue comme 
intermédiaire d'une inclusion de profondeur $2$.
\begin{prop}[NV 2 - 4.1] 
Soit  $P_{0} \subset P_{1}$ 
une inclusion d'indice fini  et de profondeur finie $p$ de facteurs de 
type $\mathrm{II}_{1}$. Si l'entier $m$ est supérieur ou égal à $p-1$ alors 
l'inclusion $P_{0} \subset P_{m}$ est de profondeur $2$.
\end{prop}

\subsubsection{} On suppose donc que l'inclusion $P_{0} \subset P_{m}$ 
est de profondeur $2$ (\ref{prof2}). Pour cette inclusion, on prend 
les notations suivantes :
$$N_{0}=P_0 \,\subset \,N_{1}= P_m \, \raisebox{1.3 ex}{$\begin{matrix} f_{1}\\ \subset \end{matrix}$} \, 
N_2=P_{2m} \,\raisebox{1.3 ex}{$\begin{matrix} f_{2}\\ \subset \end{matrix}$} 
\,  N_3 =P_{3m}$$
D'après [PP2], les projecteurs $f_{j}$ s'expriment en fonction des 
projecteurs $e_{i}$, on a par exemple :
$$f_{1}=\delta^{-m(m-1)/2} (e_{m}e_{m-1}\dots e_{1})(e_{m+1}e_{m}\dots 
e_{2})(e_{2m-1}e_{2m-2}\dots e_{m}).$$

 L'anti-automorphisme $j_{n}$ est l'anti-automorphisme 
de $P'_{0} \cap P_{2n}$ défini en posant :
$$ j_{n}(x)=J_{n}x^*J_{n}\qquad( x\in P'_{0} \cap P_{2n}).$$
où $J_{n}$ l'isométrie bijective 
anti-linéaire canonique de l'espace standard $L^2(P_{n},tr)$ de 
$P_{n}$ ($n \in \mathbb N$) (voir \ref{thj}).

Les commutants relatifs $A=N'_{0} \cap N_{2}$ et $B=N'_{1} \cap 
N_{3}$ sont donc munis de structures duales  de C*-groupoïde quantique.

\subsection{Facteur intermédiaire et *-sous-algèbre co-idéale}

D'après le théorème 4.3 de [NV 2], le commutant relatif $P'_{m} \cap P_{2m+1}$
est une *-sous-algèbre co-idéale de $B$	et $P_{2m+1}$ est isomorphe 
au produit croisé de $P_{2m}$ par ce co-idéal. Précisons la 
structure de  $P'_{m} \cap P_{2m+1}$, son action sur $P_{2m}$ et les 
points fixes de $P_{2m}$ sous cette action.

\subsubsection{}\label{copideal}
Nous gardons les notations de \ref{formcop} pour la proposition 
suivante.
\begin{prop}
Le co-produit d'un élément $y$ du co-idéal à gauche $P'_{m} \cap 
P_{2m+1}$ est donné par :
$$\Delta_{B}(y)=\sum_{l \in L} 
\sum_{r\in 
R}E_{A}(y\alpha_{r}\mu^{*}_{l})h^{-1}f_{2}h^{-1}\alpha^*_{r} \otimes 
\mu_{l}$$
où $\{\mu_{l}, l\in L\}$ est une famille d'unités matricielles 
normalisées de $P'_{m} \cap P_{2m+1}$.

Les restrictions à $P'_{m} \cap P_{2m+1}$ de la co-unité et de la 
co-unité but de $B$ vérifient pour  tout élément $y$  de $P'_{m} \cap P_{2m+1}$ :

$$\varepsilon_{B}(y)= \delta^{-1}tr(yhe_{2m}h)\qquad
\varepsilon^t_{B}(y)= \delta^{-1}E_{N_{2}}(yhe_{2m}h^{-1})$$

La restriction à $P'_{m} \cap P_{2m+1}$ de l'action de $B$ sur $P_{2m}$ est donnée 
par 
$$y\triangleright x= \delta^{-1}E_{N_{2}}(yxhe_{2m}h^{-1})\qquad 
(y\in P'_{m} \cap P_{2m+1} ,x\in P_{2m})$$

Si l'on appelle algèbre des points fixes de $P_{2m}$ sous l'action du 
co-idéal $P'_{m} \cap P_{2m+1}$, l'algèbre $P_{2m}^f$ définie par :
$$P_{2m}^f=\{x\in P_{2m},y\triangleright x=\varepsilon^t_{B}(y)\triangleright 
x,\forall y\in P'_{m} \cap P_{2m+1}\}$$
alors $P_{2m}^f$ est l'image de $P_{2m-1}$ par l'automorphisme 
intérieur $\Ad(h)$.
\end{prop}

\dem

Pour le co-produit, on écrit la deuxième formule :
$$\Delta_{B}(y)= \sum_{p\in P} 
\sum_{r\in R}E_{A}(y\alpha_{r}b^{*}_{p})h^{-1}f_{2}h^{-1}\alpha^*_{r} \otimes b_{p}$$
Comme $y\alpha_{r}$ est un élément de $P'_{0} \cap P_{2m+1}$, on a :
$$E_{A}(y\alpha_{r}b^{*}_{p})=E_{A}(y\alpha_{r}E_{P_{2m+1}}(b^{*}_{p}))$$
Soit $\{\mu_{l}, l\in L\}$ une famille d'unités matricielles 
normalisées de $P'_{m} \cap P_{2m+1}$.
En décomposant $E_{P_{2m+1}}(b^{*}_{p})$ sur cette base, on obtient 
la formule annoncée.

Les formules pour la co-unité, la co-unité but et l'action résultent de leurs définitions 
et du lemme suivant :

\begin{lem} $E_{P_{2m+1}}(f_{2})$ vaut $\delta^{m-1} e_{2m}$.
\end{lem}

\dem

Comme $f_{2}$ est le projecteur de Jones de l'inclusion $P_{m} \subset 
P_{2m}$ qui est d'indice $\delta^{-m}$, on sait que $E_{P_{2m}}(f_{2})$ 
est le scalaire $\delta^{m}$. D'autre part, d'après [PP2], on a 
l'égalité : $f_{2}=f_{2}e_{2m}$. On en déduit, pour tout $x$ de 
$P_{2m+1}$ :
$$tr(E_{P_{2m+1}}(f_{2})x)=tr(f_{2}x)=tr(f_{2}e_{2m}x)=\delta^{-1} tr(f_{2}e_{2m}E_{P_{2m}}(e_{2m}x))
=\delta^{m-1}tr(e_{2m}x)$$
La formule annoncée en résulte.\end{proof}

La démonstration du troisième point  est analogue à celle de \ref{ptsfixes}.
\end{proof}

\subsubsection{} De cette proposition et du théorème 4.3 de [NV 2], on déduit le corollaire :
\begin{cor}
La tour $P_{1} \subset P_{2} \subset P_{3}$ est isomorphe 	à la tour  
$$hP_{2m-1}h^{-1} \subset P_{2m} \subset P_{2m}\ltimes (P'_{m} \cap P_{2m+1}).$$

Si $m$ est impair, la tour $P_{0} \subset P_{1} \subset P_{2}$ est isomorphe 	à la tour  
$$j_{m}(h)P_{m-1}j_{m}(h^{-1}) \subset P_{m} \subset P_{m}\ltimes (P'_{0} \cap P_{m+1}).$$
\end{cor}

La deuxième assertion résulte de la symétrie de la construction.

\subsubsection{Calcul de $\Delta_{B}(e_{2m})$}\label{e2m}

\begin{cor} Si $\{\mu_{l}, l\in L\}$ est une famille d'unités matricielles 
normalisées de \break $P'_{m} \cap P_{2m+1}$, on a la formule :
$$\Delta_{B}(e_{2m})= \delta \sum_{l \in L} j_{2m}(h^{-1}\mu^{*}_{l}) \otimes 
h^{-1}\mu_{l}$$
\end{cor}

\dem

De la proposition et du lemme \ref{copideal}, on déduit :
$$\Delta_{B}(e_{2m})= \delta^{-(m-1)}\sum_{l \in L} 
\sum_{r\in 
R}E_{A}(f_{2}\alpha_{r}\mu^{*}_{l})f_{2}\alpha^*_{r}j_{2m}(h^{-1}) \otimes 
h^{-1}\mu_{l}$$
Et la formule 3.2.1 de [Da1] permet d'écrire l'égalité annoncée.\end{proof}

\subsection{Autodualité} Pour la fin de cette partie, nous choisissons
  $m$ pair ($m=2k$) et nous adoptons les notations suivantes :
 $M_{q}$ est le facteur $P_{qk}$
 et $N_{q}$ le facteur $M_{2q}$ pour tout entier $q$.
 L'inclusion $N_{0} \subset N_{1}$ est donc de profondeur $2$  et  
  son indice est $\tau^{-1}=[M_{1}:M_{0}]^2$.
 Les projecteurs de  Jones sont indiqués sur les tours :
 $$\begin{matrix}
 &M_0 \,\subset \, &M_1 \, \raisebox{1.3 ex}{$\begin{matrix} e_{1}\\ \subset \end{matrix}$} \, 
&M_2 \,\raisebox{1.3 ex}{$\begin{matrix} e_{2}\\ \subset \end{matrix}$} \,  &M_3 \, \raisebox{1.3
ex}{$\begin{matrix} e_{3}\\ \subset \end{matrix}$} \, &M_4  \, \raisebox{1.3
ex}{$\begin{matrix} e_{4}\\ \subset \end{matrix}$} \, &M_5\, \raisebox{1.3
ex}{$\begin{matrix} e_{5}\\ \subset \end{matrix}$} \, &M_6\\
 &N_0 \,&\subset \, &N_1 \, &\raisebox{1.3 ex}{$\begin{matrix} f_{1}\\ \subset \end{matrix}$} \, 
&N_2 \,&\raisebox{1.3 ex}{$\begin{matrix} f_{2}\\ \subset \end{matrix}$} 
\,  &N_3 
\end{matrix}$$

Ici l'anti-automorphisme $j_{q}$ est l'anti-automorphisme 
de $M'_{0} \cap M_{2q}$ défini à partir de l'isométrie bijective 
anti-linéaire canonique de l'espace standard $L^2(M_{q},tr)$. Dans les diverses 
formules de la partie $3$, il faut donc remplacer $j_{2}$ par $j_{4}$ 
et $j_{1}$ par $j_{2}$.

 Les algèbres $A=N'_{0}\cap N_{2}$ 
et $B=N'_{1}\cap N_{3}$ sont munies de structures de C*-groupoïde 
quantique que nous allons comparer. Quand on choisit  
$m$ pair, on bénéficie de l'existence de l'isomorphisme 
$\gamma=j_{4}j_{3}=j_{3}j_{2}$ qui envoie l'algèbre  
involutive  $(A,j_{2})$ sur $(B,j_{4})$ (\ref{gamma}). L'isomorphisme $\gamma$ 
décale de 2 les indices de la tour des $M_{n}$
 mais de 1 ceux de la tour des $N_{q}$. De plus 
comme l'anti-isomorphisme $j_{3}$ de $N'_{0} \cap N_{3}$ conserve la 
trace et échange les algèbres 
$A$ et $B$, il conserve aussi l'indice $H$ de 
la restriction à $N'_{1} \cap N_{2}$ de la trace $tr$, on a donc :
$$\gamma(j_{2}(h))=h \qquad \text{et}\qquad \gamma(h)=j_{4}(h).$$

Ces remarques et la symétrie des structures de C*-groupoïde quantique  
permet\-tent d'affirmer que $\gamma$ est
 un isomorphisme du C*-groupoïde quantique 
$A$
sur le C*-groupoïde quantique $B$.

\begin{thm}
Soit  $P_{0} \subset P_{1}$ 
une inclusion d'indice fini  et de profondeur finie  de facteurs de 
type $\mathrm{II}_{1}$ telle que $P_{0} \subset P_{2k}$ soit de profondeur 
$2$. Alors  les C*-groupoïdes quantiques $P'_{0} \cap P_{4k}$ et $P'_{2k} \cap 
P_{6k}$ sont isomorphes, ils sont 
donc autoduaux.
\end{thm}

\section{Structure de C*-groupoïde quantique sur les algèbres de 
Temperley-Lieb}
Comme dans [NV 3 - 2.7], nous précisons la structure de C*-groupoïde quantique 
sur les algèbres de Temperley-Lieb dans le cas non générique. Ces algèbres sont apparues dès le 
début de l'étude des inclusions ([J]). Ce sont 
les commutants relatifs  des facteurs de Jones dans le facteur 
hyperfini.

\subsection{Facteurs de Jones et algèbres de Temperley-Lieb} [GHJ - 
2.1, 4.7b,II.7].

Soient $l$ un entier supérieur à 2  et $(e'_{i})_{i?0}$ une 
suite de projecteurs satisfaisant les relations suivantes :
$$e'_{i}e'_{i \pm 1}e'_{i}=\delta e'_{i}\quad \text{et} \quad 
e'_{i}e'_{j}=e'_{j}e'_{i}\quad \text{pour } |i-j| \geq 2$$
avec $\delta=(4\cos^2\frac{\pi}{l+1})^{-1}$.
Le facteur $P_{1}$ engendré par les projecteurs $(e'_{i})_{i?0}$ est 
le facteur hyperfini de type $\mathrm{II}_{1}$ et le sous-facteur  $P_{0}$ de $P_{1}$ 
engendré par les projecteurs $(e'_{i})_{i?1}$ est le sous-facteur de 
Jones d'indice $\delta^{-1}$. On note $tr$ la trace normalisée de $P_{1}$.
Le graphe principal de l'inclusion $P_{0} 
\subset P_{1}$ est le graphe 
linéaire $A_{l}$ à $l$ sommets (voir 
[GHJ - 1.4.3], [J - 4,5]), de même pour  l'inclusion $P_{1} \subset 
P_{2}$ qui lui est isomorphe. La profondeur de ces inclusions est 
donc $l-1$. Les commutants relatifs sont des algèbres de 
Temperley-Lieb de paramètre non-générique $\delta^{-1}$ :
$$P'_{q}\cap P_{n}=(1,e_{q+1},e_{q+2}\dots e_{n-2}, e_{n-1})"$$

D'après [NV 2 - 4.1] (voir \ref{prof2}), pour $m=l-2$, l'inclusion $P_{0} \subset P_{m}$ 
est de profondeur $2$, de plus elle est isomorphe \`a l'inclusion $P_{m} \subset P_{2m}$.
 L'algèbre de Temperley-Lieb $A=P'_{0} \cap P_{2m}  
=(1,e_{1},e_{2},\dots,e_{2m-1})"$ est donc munie d'une structure de C*-groupoïde 
quantique autodual qu'on va préciser.

\subsection{Structure de C*-groupoïde quantique des algèbres de Temperley-Lieb} \label{TL}

\begin{prop} Le co-produit de $A$ est donné par :
\begin{align*}
\Delta_{A}(1)&= \sum_{k}\frac{1}{\nu_{j_{k}}} j_{m}(\lambda^{*}_{k}) \otimes \lambda_{k}\\
\Delta_{A}(e_{p}) &= \Delta_{A}(1) (e_{p} \otimes 1)= (e_{p} \otimes 
1)\Delta_{A}(1)
\qquad &(1\leq p \leq m-1)
\\
 \Delta_{A}(e_{q})&=\Delta_{A}(1)(1 \otimes e_{q})=(1 \otimes  
 e_{q})\Delta_{A}(1)
\qquad &(m+1\leq q \leq 2m-1)\\
\Delta_{A}(e_{m}) &=\delta \sum_{l \in L} j_{m}(j_{m}(h^{-1})\mu^{*}_{l}) \otimes 
j_{m}(h^{-1})\mu_{l} .
\end{align*}
où $\{\lambda_{k}, k \in K\}$ est une famille d'unités matricielles 
 de $P'_{0} \cap P_{m}$ ($\nu_{j_{k}}$ est la dimension du facteur de 
 $P'_{0} \cap P_{m}$ auquel appartient $\lambda_{k}$), $\{\mu_{l}, l\in L\}$ une famille d'unités matricielles 
normalisées de $P'_{0} \cap P_{m+1}$ et $h$ la racine carrée de 
l'indice de la restriction à $P'_{m} \cap P_{2m}$ de $tr$.

La co-unité de $A$ est donnée par :
$$\varepsilon_{A}(x)=\delta^{-m}tr(hf_{1}hx)\quad (x \in A)$$
où $f_{1}$ est le projecteur de Jones de l'inclusion $P_{0} \subset 
P_{m}$ :
$$f_{1}=\delta^{-m(m-1)/2} (e_{m}e_{m-1}\dots e_{1})(e_{m+1}e_{m}\dots 
e_{2})(e_{2m-1}e_{2m-2}\dots e_{m}).$$

L'antipode  de $A$ est donnée par :
\begin{align*}
	S_{A}(e_{p})& =  e_{2m-p}
\qquad (1\leq p\leq 2m-1,p \neq m)\\
S_{A}(e_{m}) &= j_{m}(h)h^{-1}e_{m}j_{m}(h^{-1})h
\end{align*}
\end{prop}

\dem

Le co-produit est un homomorphisme d'algèbres, il suffit donc de le 
connaître sur les générateurs $1, e_{1},\dots e_{2m-1}$ de $A$. 
L'expression de $\Delta(1)$ résulte de \ref{delta(1)}.
Tous les projecteurs 
sauf $e_{m}$ sont soit dans $N'_{0} \cap N_{1}$ soit dans $N'_{1} \cap 
N_{2}$, leur co-produit est calculé grâce à \ref{coppa} et \ref{delta(1)}.

L'expression de $\Delta(e_{m})$ est donnée par le corollaire 
\ref{e2m} grâce à
la symétrie de la construction.

La formule de la co-unité découle de \ref{coun} et de [PiPo 2].

 L'antipode est un 
anti-automorphisme d'algèbres conservant l'unité, il suffit donc de la 
connaître sur les projecteurs de Jones. D'après [Da1-2.2.1], on sait 
que, pour $p=1\dots 2m-1$, $j_{m}(e_{p})$ est le projecteur $e_{2m-p}$ ; 
les formules résultent alors  des propriétés de commutation de $h$.
\end{proof}

\subsection{C*-groupoïde quantique de dimension $13$ associé au facteur de Jones de graphe 
$A_{4}$}\label{chemins}

 Dans cette partie, on suppose que $l$ vaut $4$, l'inclusion $P_{0} \subset 
 P_{1}$ est alors d'indice $\delta^{-1}=4\cos^2 \frac{\pi}{5}$, de graphe 
principal  $A_{4}$ et  l'inclusion $P_{0} \subset P_{2}$ 
 est  de profondeur 2.
La C*-algèbre  $A_{\delta}=P'_{0} \cap P_{4}$ est  un C*-groupoïde quantique autodual ;
  nous  le décrivons  
 et montrons qu'il  est isomorphe à celui, que nous nommerons $G$, décrit
    par G. Böhm et K. Szlachanyi dans [BSz - 5]. Dans [NV1], D. Nikshych 
    et L. Vainerman munissent cette même algèbre d'une structure de 
    groupoïde quantique pour laquelle l'involution est modifiée mais 
    on peut montrer par les méthodes employées ici qu'elle est 
    isomorphe aux deux autres.
    
Pour simplifier les calculs, 
nous utilisons le paramètre $z=\sqrt[4]{\delta}$ introduit  dans 
[BSz-5] 
et qui vérifie les relations suivantes et bien d'autres encore :
$$\begin{array}{lll}
	1-3\delta +\delta^2=0 \qquad
&\qquad	z^4+z^2-1=0 &\qquad	1+z^2=z^{-2}\\
	z^2=1-\delta=\displaystyle{\frac{\delta}{1-\delta}}
&\qquad	z^3=\sqrt{\delta (1-\delta)}
\end{array}$$

\subsubsection{Algèbre des chemins de $A_{\delta}$}\label{graphe}Comme dans [GHJ 2.3.11],
 nous représentons l'algèbre $A_{\delta}$ comme algèbre des 
chemins du graphe $A_{4}$ avec les notations suivantes pour les 
sommets du graphe et les chemins ;
la trace des projecteurs minimaux des algèbres correspondant aux 
sommets du graphe est donnée sur le graphe  de droite (d'après [J 5.2]) :

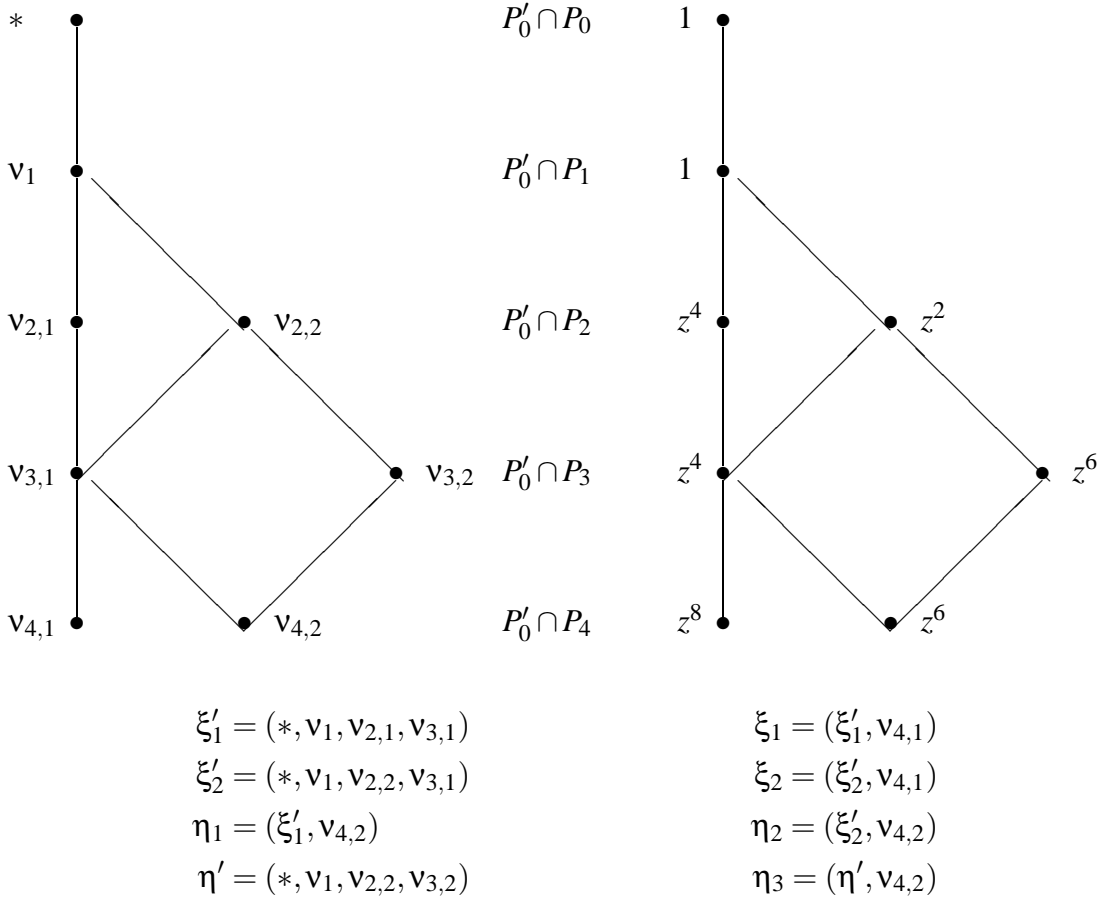
\unitlength=1cm
\begin{picture}(8,9)(0,0)
	%le nom des sommets
\put(0,8){{$*$}}
\put(0,6){{$\nu_{1}$}}
\put(0,4){{$\nu_{2,1}$}}
\put(0,2){{$\nu_{3,1}$}}
\put(0,0){{$\nu_{4,1}$}}
\put(3.5,0){{$\nu_{4,2}$}}
\put(3.5,4){{$\nu_{2,2}$}}
\put(5.5,2){{$\nu_{{3,2}}$}}
%les points
\put(0.8,0){{$\bullet$}}
\put(0.8,2){{$\bullet$}}
\put(0.8,4){{$\bullet$}}
\put(0.8,6){{$\bullet$}}
\put(0.8,8){{$\bullet$}}
\put(3,0){{$\bullet$}}
\put(3,4){{$\bullet$}}
\put(5,2){{$\bullet$}}
%les algèbres
\put(6.5,8){$P'_{0} \cap P_{0}$}
\put(6.5,6){$P'_{0} \cap P_{1}$}
\put(6.5,4){$P'_{0} \cap P_{2}$}
\put(6.5,2){$P'_{0} \cap P_{3}$}
\put(6.5,0){$P'_{0} \cap P_{4}$}
%les traits
\put(0.9,0.2){\line(0,1){1.8}}
\put(0.9,2.2){\line(0,1){1.8}}
\put(0.9,4.2){\line(0,1){1.8}}
\put(0.9,6.2){\line(0,1){1.8}}
\put(3.1,0){\line(-1,1){2}}
\put(3.1,0){\line(1,1){2}}
\put(0.9,2){\line(1,1){2}}
\put(5.2,2){\line(-1,1){2}}
\put(3.1,4){\line(-1,1){2}}
%deuxième diagramme
\put(8.5,0){\begin{picture}(6,9)
	%le nom des sommets
\put(0.3,8){{$1$}}
\put(0.3,6){{$1$}}
\put(0.3,4){{$z^4$}}
\put(0.3,2){{$z^4$}}
\put(0.3,0){{$z^8$}}
\put(3.5,0){{$z^6$}}
\put(3.5,4){{$z^2$}}
\put(5.5,2){{$z^6$}}
%les points
\put(0.8,0){{$\bullet$}}
\put(0.8,2){{$\bullet$}}
\put(0.8,4){{$\bullet$}}
\put(0.8,6){{$\bullet$}}
\put(0.8,8){{$\bullet$}}
\put(3,0){{$\bullet$}}
\put(3,4){{$\bullet$}}
\put(5,2){{$\bullet$}}
%les traits
\put(0.9,0.2){\line(0,1){1.8}}
\put(0.9,2.2){\line(0,1){1.8}}
\put(0.9,4.2){\line(0,1){1.8}}
\put(0.9,6.2){\line(0,1){1.8}}
\put(3.1,0){\line(-1,1){2}}
\put(3.1,0){\line(1,1){2}}
\put(0.9,2){\line(1,1){2}}
\put(5.2,2){\line(-1,1){2}}
\put(3.1,4){\line(-1,1){2}}
\end{picture}}
\end{picture}

\begin{align*}
     \xi'_{1}&=(*,\nu_{1}, \nu_{2,1}, \nu_{3,1}) \qquad
	&\xi_{1}=(\xi'_{1}, \nu_{4,1})\\
	\xi'_{2}&=(*,\nu_{1}, \nu_{2,2}, \nu_{3,1})\qquad
	&\xi_{2}=(\xi'_{2}, \nu_{4,1})\\
	\eta_{1}&=(\xi'_{1}, \nu_{4,2})\qquad
	&\eta_{2}=(\xi'_{2}, \nu_{4,2})\\
	\eta'&=(*,\nu_{1}, \nu_{2,2}, \nu_{3,2})\qquad
&	\eta_{3}=(\eta', \nu_{4,2})
\end{align*}

L'algèbre $A_{\delta}$ est somme directe de l'algèbre $C=\Vect\{c_{i,j},(i,j) \in \{1,2\}^2\}$ 
et l'algèbre $D=\Vect\{d_{h,k},(h,k) \in \{1,2,3\}^2\}$ avec 
$c_{i,j}=T_{\xi_i,\xi_j}$ et  $d_{h,k}=T_{\eta_h,\eta_k}$.

L'algèbre $P'_{0} \cap P_{3}$ est somme directe de $\Vect\{b_{i,j},(i,j) \in \{1,2\}^2\}$ 
et $\mathbb{C}b_{5}$ avec 
$b_{i,j}=T_{\xi'_i,\xi'_j}=c_{i,j}+d_{i,j}$ et  $b_{5}=T_{\eta',\eta'}=d_{3,3}$.

\subsubsection{Projecteurs de Jones} L'algèbre $A_{\delta}$ est 
engendrée, en tant qu'algèbre, par l'unité et les projecteurs $e_{1}, e_{2}$ et $ e_{3}$. La formule [GHJ 2.6.5.4] nous 
fournit l'expression de ces projecteurs dans la base d'unités 
matricielles de l'algèbre des chemins.

D'après [GHJ 2.6.5.4], on a donc :
$$\begin{array}{cccc}
	e_{1}&=\begin{pmatrix}1&0 \\ 0&0 \end{pmatrix} &+& \begin{pmatrix}1&0&0 
	\\ 0&0&0 \\0&0&0 \end{pmatrix}\\
	e_{2}&=\begin{pmatrix}z^4 &z^3 \\ z^3&z^2\end{pmatrix} &+& 
	\begin{pmatrix}z^4 &z^3&0 \\ z^3&z^2&0 \\0&0&0 \end{pmatrix}\\
	e_{3}&=\begin{pmatrix}1&0 \\ 0&0 \end{pmatrix} &+& \begin{pmatrix}0&0&0 
	\\ 0&z^2&z^3 \\ 0&z^3&z^4 \end{pmatrix}
\end{array}$$
Dans ce cas, d'après  \ref{H1} et \ref{thj}, on a  :
 $$h= z^{-2}e_{3}+z^{-1}(1-e_{3})\qquad j_{2}(h)= 
 z^{-2}e_{1}+z^{-1}(1-e_{1}).$$

\subsubsection{Expression des unités matricielles en fonction des projecteurs 
de Jones}\label{cd}
On vérifie facilement par le calcul les expressions suivantes pour 
les unités matricielles de $A_{\delta}$ :
\begin{align*}
	c_{1,1}&=e_{1}e_{3}\\
	c_{1,2}&=z^{-3}e_{3}e_{1}( e_{2}-\delta)\\
	d_{1,1}&=e_{1}(1-e_{3})\\
	d_{1,2}&=z^{-3}(1-e_{3})e_{1}( e_{2}-\delta)\\
	d_{1,3}&=z^{-6} e_{1}(1-e_{3})(e_{2}-\delta)(e_{3}-z^2)	
	\end{align*}
Et on obtient les autres expressions grâce aux relations entre les unités 
matricielles.

\subsubsection{Nouvelles unités matricielles}\label{num}  Nous définissons maintenant 
des nouvelles unités matricielles sur $D$ qui vont permettre d'obtenir des 
formules plus simples pour le co-produit et l'antipode et d'identifier 
$A_{\delta}$ et le C*-groupoïde quantique que nous appellerons $G$ décrit en [BSz 5].

\begin{prop}$(D,*)$ admet $\{e_{i,j}, i,j=1,2,3\}$ comme unités 
matricielles avec : 
\begin{align*}
	e_{1,2}&=z^{2} d_{1,2}-zd_{1,3}=z^{-3}e_{1}(1-e_{3})(e_{2}-\delta)(1-e_{3})\\
	e_{1,3}&=z d_{1,2}+z^{2}d_{1,3}=\delta^{-1}e_{1}(1-e_{3})(e_{2}-\delta)e_{3}
\end{align*}
 On a alors en particulier :
 \begin{align*} 
 e_{1}&=c_{1,1} + e_{1,1} \\
e_{3}&=c_{1,1} + e_{3,3}\\
e_{2}&=z^4 c_{1,1} +z^3 c_{1,2} + z^3 c_{2,1} + z^2 c_{2,2}
+ z^4 e_{1,1}+z^5 e_{1,2} +z^4 e_{1,3}\\
&\;+ z^5 e_{2,1}+z^6 e_{2,2}+z^5 e_{2,3}
 + z^4 e_{3,1}+z^5 e_{3,2}+ z^4 e_{3,3}
 \end{align*}
\end{prop}

Les calculs nécessaires à la vérification de cette proposition et des 
suivantes se font facilement à l'aide d'un logiciel de calcul formel.

\subsubsection{Expression du co-produit} 
\begin{prop}[BSz1 5] Si $\{e^0_{i,j},i,j=1,2\}$ (resp. 
$\{e^1_{i,j},i,j=1,2,3\}$) est une famille 
d'unités matricielles du facteur de dimension $4$ (resp. $9$) de $G$,
 le co-produit de $G$ est donné par :
\begin{align*}
\Delta_{G}(e^0_{1,1})&=e^0_{1,1} \otimes e^0_{1,1}+ e^1_{1,1} 
\otimes e^1_{3,3}\\
\Delta_{G}(e^0_{1,2})&=e^0_{1,2} \otimes e^0_{1,2}+ z^2e^1_{1,3} \otimes e^1_{3,1}
+ z e^1_{1,2} \otimes e^1_{3,2}\\
\Delta_{G}(e^0_{2,2})&=e^0_{2,2} \otimes e^0_{2,2}+ z^4e^1_{3,3} \otimes e^1_{1,1}
+ z^3 e^1_{3,2} \otimes e^1_{1,2}+ z^3e^1_{2,3} \otimes e^1_{2,1}
+ z^2 e^1_{2,2} \otimes e^1_{2,2}\\
\Delta_{G}(e^1_{1,1})&=e^0_{1,1} \otimes e^1_{1,1}+e^1_{1,1} \otimes e^0_{2,2}+ e^1_{1,1} 
\otimes e^1_{2,2}\\
\Delta_{G}(e^1_{1,2})&=e^0_{1,2} \otimes e^1_{1,2}+ e^1_{1,2} \otimes e^0_{2,2}
+ z e^1_{1,3} \otimes e^1_{2,1}- z^2 e^1_{1,2} \otimes e^1_{2,2}\\
\Delta_{G}( e^1_{1,3})&=e^0_{1,2}\otimes  e^1_{1,3} + e^1_{1,3} \otimes  e^0_{2,1}
+	 e^1_{1,2} \otimes e^1_{2,3}\\
\Delta_{G}(e^1_{2,2})&=e^0_{2,2} \otimes e^1_{2,2}+e^1_{2,2} \otimes e^0_{2,2}
+z^4e^1_{2,2} \otimes e^1_{2,2} \\
&\;+ z^2e^1_{3,3} \otimes e^1_{1,1}
- z^3 e^1_{3,2} \otimes e^1_{1,2}- z^3e^1_{2,3} \otimes e^1_{2,1}\\
\Delta_{G}(e^1_{2,3})&=e^0_{2,2} \otimes e^1_{2,3}+ e^1_{2,3} \otimes e^0_{2,1}
+ z e^1_{3,2} \otimes e^1_{1,3}- z^2 e^1_{2,2} \otimes e^1_{2,3}\\
\Delta_{G}( e^1_{3,3})&=e^0_{2,2}\otimes  e^1_{3,3} + e^1_{2,2} \otimes  e^1_{3,3} +	 e^1_{3,3} \otimes e^0_{1,1}
\end{align*}
\end{prop}

Dans la proposition suivante, nous donnons deux séries de formules pour le 
coproduit $\Delta_{A}$, l'une en fonction des nouvelles unités matricielles permet de le 
comparer avec $\Delta_{G}$, l'autre en fonction des projecteurs de 
Jones permet de le comparer avec celui de [NV1 - 7.3].

\begin{prop} Le co-produit $\Delta_{A}$ de $A_{\delta}$ est 
l'homomorphisme d'algèbres déterminé 
par les égalités suivantes :
\begin{align*}
\Delta_{A}(1)&= c_{1,1} \otimes c_{1,1} + c_{1,1} \otimes e_{1,1}
 + c_{2,2} \otimes c_{2,2} + c_{2,2} \otimes e_{2,2} +c_{2,2} \otimes e_{3,3}
 \\& \quad +   e_{1,1}\otimes c_{2,2}+  e_{1,1}\otimes e_{2,2} +    e_{1,1}\otimes e_{3,3}
 +e_{2,2} \otimes c_{2,2}+ e_{2,2} \otimes e_{2,2} 
 \\& \quad + e_{2,2} \otimes e_{3,3}
  + e_{3,3}\otimes c_{1,1} + e_{3,3} \otimes e_{1,1}\\
 \Delta_{A}(e_{1})&=c_{1,1} \otimes c_{1,1} + c_{1,1} \otimes e_{1,1}+   e_{1,1}\otimes c_{2,2}
 +  e_{1,1}\otimes e_{2,2} +    e_{1,1}\otimes e_{3,3}\\
 \Delta_{A}(e_{3})&=c_{1,1} \otimes c_{1,1} + c_{2,2} \otimes 
e_{3,3} +    e_{1,1}\otimes e_{3,3}  + e_{2,2} \otimes 
e_{3,3}  + e_{3,3}\otimes c_{1,1}
\end{align*}
$$\begin{array}{llllll}
	\Delta_{A}(e_{2})&=
		z^4 c_{1,1} \otimes c_{1,1} &+ z^4 c_{1,1} \otimes e_{1,1}\\
		&+z^3 c_{1,2} \otimes c_{1,2}& + z^5 c_{1,2} \otimes e_{1,2}&+z^4 
		c_{1,2} \otimes e_{1,3}\\
		&+z^3 c_{2,1} \otimes c_{2,1} &+ z^5 c_{2,1} \otimes e_{2,1}&+z^4 
		c_{2,1} \otimes e_{3,1}\\
		&+z^2 c_{2,2} \otimes c_{2,2} &+ z^6 c_{2,2} \otimes e_{2,2}&+z^5 c_{2,2} \otimes e_{2,3}&+ 
		z^5 c_{2,2} \otimes e_{3,2}&+z^4 c_{2,2} \otimes e_{3,3}\\
		&+z^4 e_{1,1} \otimes c_{2,2} &+ z^4 e_{1,1} \otimes e_{2,2}&+z^4 e_{1,1} \otimes e_{3,3}\\
		&+z^5 e_{1,2} \otimes c_{2,2} &- z^7 e_{1,2} \otimes e_{2,2}
		& + z^4 e_{1,2} \otimes e_{2,3}&+z^4 e_{1,2} \otimes e_{3,2}\\
		&+ z^4 e_{1,3} \otimes c_{2,1} &+ z^6 e_{1,3} \otimes e_{2,1} 
    	&+z^5 e_{1,3} \otimes e_{3,1}\\
		&+ z^5 e_{2,1} \otimes c_{2,2} &- z^7 e_{2,1} \otimes e_{2,2}
		&+z^4 e_{2,1} \otimes e_{2,3}&+z^4 e_{2,1} \otimes e_{3,2}\\
		&+z^6 e_{2,2} \otimes c_{2,2} &+ 2z^6 e_{2,2} \otimes e_{2,2}&-z^7 e_{2,2} \otimes e_{2,3}
    	&	- z^7 e_{2,2} \otimes e_{3,2}&+z^4 e_{2,2} \otimes e_{3,3}\\
		&+z^5 e_{2,3} \otimes c_{2,1} &+ z^7 e_{2,3} \otimes e_{2,1} &+z^6 e_{2,3} \otimes e_{3,1}\\
		&+z^4 e_{3,1} \otimes c_{1,2} &+ z^6 e_{3,1} \otimes e_{1,2}&+z^5 e_{3,1} \otimes e_{1,3}\\
		&+z^5 e_{3,2} \otimes c_{1,2} &+ z^7 e_{3,2} \otimes e_{1,2}&+z^6 e_{3,2} \otimes e_{1,3}\\
		&+z^4 e_{3,3} \otimes c_{1,1} &+ z^4 e_{3,3} \otimes e_{1,1}
\end{array}$$
On a aussi :
\begin{align*}
 \Delta_{A}(1)&= e_{3}\otimes e_{1} + (1-e_{3}) \otimes (1-e_{1})\\
 \Delta_{A}(e_{1})&= e_{1}e_{3}\otimes e_{1} + e_{1}(1-e_{3}) \otimes (1-e_{1})	\\
 \Delta_{A}(e_{3})&= e_{3}\otimes e_{1}e_{3} + (1-e_{3}) \otimes (1-e_{1})e_{3}
\end{align*}
\begin{align*}
\Delta_{A}(e_{2})&=\bigg(1-\frac{(e_{3}-e_{2})^2}{(1-\delta)}\bigg) \otimes
\bigg(1-\frac{(e_{1}-e_{2})^2}{(1-\delta)}\bigg)+\delta e_{3}\otimes e_{1}&\\
&+\;\frac{1}{\sqrt{\delta(1-\delta)}}\; e_{3}(e_{2}-\delta)\otimes e_{1}(e_{2}-\delta)
+\frac{1}{\sqrt{\delta(1-\delta)}}\;  (e_{2}-\delta)e_{3}\otimes 
(e_{2}-\delta)e_{1}\\
&+(1-\delta)\bigg(\frac{(e_{3}-e_{2})^2}{(1-\delta)}-e_{3}\bigg)\otimes 
\bigg(\frac{(e_{1}-e_{2})^2}{(1-\delta)}-e_{1}\bigg) 	
\end{align*}

\end{prop}

\dem

Comme l'algèbre $P'_{0} \cap  P_{2}$ égale $\mathbb{C}e_{1} \otimes 
\mathbb{C}(1-e_{1})$, on a :
$$\Delta_{A}(1)= e_{3}\otimes e_{1} + (1-e_{3}) \otimes (1-e_{1}).$$
En utilisant les expressions des projecteurs de Jones en fonction des 
unités matricielles, on obtient :
$$\Delta_{A}(1)= (c_{1,1} + e_{3,3})\otimes (c_{1,1} + e_{1,1}) + (c_{2,2} + 
 e_{1,1}+e_{2,2}) \otimes (c_{2,2}+ e_{2,2} +  e_{3,3}).$$
 
 De même, les autres formules sont la traduction à l'aide des unités 
 matricielles des égalités (\ref{TL}) :
 $$\Delta_{A}(e_{1})=(e_{1}\otimes 1)\Delta_{A}(1) \qquad \qquad 
 \Delta_{A}(e_{3})=\Delta_{A}(1)(1 \otimes e_{3}).$$
 
 Le calcul de $\Delta_{A}(e_{2})$ demande un peu plus de travail :
 D'après \ref{graphe} et  \ref{cd}, les unités matricielles normalisées de $P'_{0} \cap 
 P_{3}$ sont :
 $$\begin{array}{lll}
 \mu_{1}=z^{-2}e_{1} &\mu_{2}=z^{-2}(c_{1,2}+d_{1,2})\;\; &\mu_{3}=z^{-2}(c_{2,1}+d_{2,1})\\
 	\mu_{4}=z^{-2}(c_{2,2}+d_{2,2})\;\; & \mu_{5}=z^{-3}d_{3,3}
 \end{array}$$
 Comme $j_{2}(h^{-1})$ vaut $z^{2}e_{1}+z(1-e_{1})$, la formule
 \ref{TL} s'écrit ici :
\begin{align*}
	\Delta_{A}(e_{2})&= \delta (j_{2}(z^2 \mu_{1}) \otimes z^2 \mu_{1} 
 +  j_{2}(z \mu_{3}) \otimes z^2 \mu_{2}+  j_{2}(z^2 \mu_{2}) \otimes z \mu_{3} \\
&\qquad\qquad + j_{2}(z \mu_{4}) \otimes z \mu_{4}
+  j_{2}(z \mu_{5}) \otimes z \mu_{5})\\
&=j_{2}(z^2 e_{1}) \otimes z^2 e_{1} 
+ j_{2}(z b_{2,1}) \otimes z^2 b_{1,2}+ j_{2}(z^2 b_{1,2}) \otimes z b_{2,1}\\
&\qquad\qquad + j_{2}(z b_{2,2}) \otimes z b_{2,2} +j_{2}(d_{3,3}) \otimes d_{3,3}
\end{align*}
Or on a :
\begin{align*}
	b_{1,2}&=z^{-3}e_{1}(e_{2}-\delta)\\
	b_{2,2}&=z^{-6}(e_{2}-\delta)e_{1}(e_{2}-\delta)
	=z^{-2}[(e_{2}-e_{1})^2-(1-\delta)e_{1}]\\
	d_{3,3}&=z^{-12}(e_{3}-z^2)(e_{2}-\delta)(1-e_{3})e_{1}(e_{2}-\delta)(e_{3}-z^2)=
	1-\frac{(e_{1}-e_{2})^2}{(1-\delta)}
\end{align*}
Ces expressions des unités matricielles en fonction des projecteurs de 
Jones permettent de préciser les valeurs prises par $j_{2}$ puis on 
exprime le résultat en fonction des nouvelles unités matricielles et 
on obtient les formules annoncées.
 \end{proof}

\begin{cor}
Le co-produit de $A_{\delta}$ coïncide avec celui de $G$.
\end{cor}

\dem

Comme le projecteur $e_{3}$ est l'image de $e_{1}$ par l'antipode qui 
est un anti-automorphisme de co-algèbre (on verra plus loin que les 
antipodes coïncident), il suffit pour comparer 
$\Delta_{A}$ et le co-produit $\Delta_{G}$ de $G$ de considérer leurs 
valeurs en $1$, $e_{1}$ et $e_{2}$. On vérifie par le calcul qu'elles 
coïncident.
\end{proof}

\subsubsection{Expression de la co-unité} Comme $f_{1}$ coïncide avec $e_{2}$ dans $C$
et est nul dans $D$, la co-unité $\varepsilon_{A}$ est nulle sur 
$D$ et comme elle est linéaire et compatible avec l'involution, 
il suffit de la connaître sur $c_{1,1}$ et 
$c_{1,2}$. Comme la trace des projecteurs minimaux de $C$ est 
$\delta^2$, si on note $Tr_{0}$ la trace de $C$ qui vaut $1$ sur les 
projecteurs minimaux, on a :
\begin{align*}
\varepsilon_{A}(c_{1,2})&=Tr_{0} (hf_{1}hc_{1,2})=1\\
\varepsilon_{A}(c_{1,1})&=Tr_{0} (f_{1}hc_{1,1})=1
\end{align*}

La co-unité de $A$ coïncide avec celle de $G$.

\subsubsection{Expression de l'antipode}
\begin{prop} L'antipode de $A_{\delta}$ est entièrement déterminée par les formules 
sui\-vantes ; elle coïncide avec l'antipode de $G$.
\begin{align*}
S_{A}(c_{1,2})&=c_{2,1} \qquad
 &S_{A}(c_{2,1})&=c_{1,2}\\
S_{A}(e_{1,2})&=z^{-1}e_{2,3}\qquad
&S_{A}(e_{2,1})&=ze_{3,2}\\
S_{A}(e_{1,3})&=z^{-2}e_{1,3}\qquad
&S_{A}(e_{3,1})&=z^{2}e_{3,1}
\end{align*}
\end{prop}
	
\dem 

		Ces formules résultent  des expressions des nouvelles unités 
matricielles en fonction des projecteurs de 
Jones (\ref{num}) et de \ref{thj} (c).
Elles suffisent pour connaître $S_{A}$ qui est un 
anti-automorphisme d'algèbre.
\end{proof}

\section{Action d'un groupoïde quantique fini sur un facteur}
Dans cette partie, nous considèrons $(A, m_{a}, 1_{a},\Delta_{a}, \varepsilon_{a}, S_{a}, \phi_{a}, p_{a})$ 
et \break $(B, m_{b}, 1_{b},\Delta_{b}, \varepsilon_{b}, S_{b}, \phi_{b}, p_{b})$ deux 
groupoïdes quantiques finis en dualité (notée $\langle a,b\rangle$). 
Le but de cette partie est de construire une inclusion $M_1\subset M_2$ de facteurs 
(que nous obtiendrons hyperfinis de type $II_{1}$) telle que 
$M_2$ soit le produit croisé $M_1\rtimes A$. 
Nous généralisons ainsi les résultats de [N] qui  concernaient les algèbres de Kac faibles.

\subsection{Hypothèse et remarque importante}\label{hyp} 
Les résultats de [NSzW] rappelés en \ref{NSzW} et ceux de la partie 3 conduisent 
imposer l'hypothèse : $A$ et $B$ sont connexes (voir \ref{connexe}).
Par contre l'hypothèse de régularité est inutile à la construction de 
l'inclusion sur laquelle agissent  $A$ et $B$. Ce qui laisse penser 
que la structure obtenue dans la partie $3$ à partir de l'inclusion 
construite n'est pas nécessairement la structure originelle (voir \ref{conclusion}).

  \subsection{Produit croisé des groupoïdes en dualité : L'algèbre $A.B$} \label{AB} 
  Les produits croisés $A\ltimes B$ et
  $A \rtimes B$ sont isomorphes en effet grâce à \ref{AsBt}, 
  l'identification, pour  $y$ dans $A_{s}$, de $[ay\otimes b]$ et 
  $[a\otimes (1_{b} \triangleleft y)b]$ dans $A\ltimes B$  
  correspond à celle, pour $z$ dans $B_{t}$, de $[a(z \triangleright 1_{a})\otimes b]
  $ et $[a\otimes zb]$ dans $A \rtimes B$. On vérifie facilement 
  que les lois sont compatibles et que les injections $a \mapsto 
  [a\otimes  1_{b}]$ et $b \mapsto [1_{a}\otimes b]$ sont des 
  homomorphismes d'algèbres qui permettent d'écrire $A\ltimes B$ et
  $A \rtimes B$ comme $A.B$.
  
  \subsection{Mesures de Haar et espérances conditionnelles} 
  \subsubsection{Mesures de Haar sur $A_{s}=B_{t}$} \label{phiAsBt}
 \begin{prop}
 Les restrictions des mesures de Haar $\phi_{a}$ et $\phi_{b}$ coïncident sur les 
 algèbres identifiées $A_{s}$ et $B_{t}$.
 \end{prop}
 \dem

D'après \ref{Haar}, la restriction à $A_{s}$ (resp. $B_{t}$) de 
 $\phi_{a}$ (resp. $\phi_{b}$) vaut $\varepsilon_{a}$ (resp. 
 $\varepsilon_{b}$). Or, pour $a$ dans $A_{s}$, on a :
 $$\varepsilon_{b}(1_{b} \triangleleft 
 a)=\langle1_{a},1_{b(2)}\rangle\langle a,1_{b(1)}\rangle=\langle a,1_{b}\rangle=
 \varepsilon_{a}(a)$$
 donc $\phi_{a}$ et $\phi_{b}$ coïncident sur les 
 algèbres identifiées $A_{s}$ et $B_{t}$.\end{proof}
%%%%%%%%%%%%%%%%%%%%%%%%%%%%%%%%%%%%%%%%%%%%%%%%%%%%%%%%%%%%%%%%%%%%%%%%%%%%%%%%%
%																				%
%																				%
%	  Si les antipodes sont	involutives	sur	les	sous-algèbres co-unitales, la	%
%	 restriction des mesures de	Haar à ces algèbres	est	leur trace				%
%	 canonique.																	%
%	 D'après [BSz -	4.6], les restrictions des co-unités aux algèbres			%
%	 co-unitales sont liées	aux	traces canoniques de ces algèbres par une		%
%	 relation du type :															%
%	 $$Tr_{A_{t}}(x)=\varepsilon_{a}(x1_{(2)} S(1_{(1)})\quad (x\in	A_{t})$$	%
%	 Si	on suppose que $S^2$ est l'identité	sur	$A_{t}$	on peut	écrire :		%
%	 $$1_{(2)} S(1_{(1)})=S(1_{(1)}	S(1_{(2)})=S(\varepsilon^t_{a}(1))=1$$		%
%	 Ainsi si les antipodes	sont involutives sur les sous-algèbres				%
%	 co-unitales,																%
%	 la	restriction	 de	$\phi_{a}$ (ou $\phi_{b}$) aux sous-algèbres			%
%	 co-unitales coïncide														%
%	 avec la trace																%
%	 canonique sur ces algèbres.												%
%																				%
%%%%%%%%%%%%%%%%%%%%%%%%%%%%%%%%%%%%%%%%%%%%%%%%%%%%%%%%%%%%%%%%%%%%%%%%%%%%%%%%%
 
  \subsubsection{Des espérances conditionnelles}\label{escond}
  \begin{prop}
  On pose pour $a$ dans $A$ :
  $$F_{A_{t}}(a)=(\id \otimes \phi_{a})\Delta_{a}(a)=p_{b}\triangleright a\quad 
  \text{et}\quad F_{A_{s}}(a)=(\phi_{a} \otimes 
  \id)\Delta_{a}(a)=a\triangleleft p_{b}$$
  Les applications $F_{A_{t}}$ et $F_{A_{s}}$ sont des espérances conditionnelles
  fidèles  de $A$ sur $A_{t}$ (resp. $A_{s}$). Elles conservent $\phi_{a}$ 
  et commutent.
  
  On définit de même $F_{B_{t}}$ et $F_{B_{s}}$ avec des 
  résultats analogues.
   \end{prop}
   \dem

   Les propriétés de $\phi_{a}$ (\ref{Haar}) et celles des éléments de $A_{t}$ 
   (\ref{rem27})
   permettent d'affirmer que pour tout $a$ de $A$, $F_{A_{t}}(a)$ 
   appartient à $A_{t}$ et que l'égalité 
   $F_{A_{t}}(xay) =xF_{A_{t}}(a)y$ est vérifiée pour tous $x$ et $y$ 
   dans $A_{t}$ et $a$ dans $A$.
   Comme $\Delta$ est un homomorphisme d'algèbres involutives, on a 
   $F_{A_{t}}(x^{*})=F_{A_{t}}(x)^{*}$.
   L'identité $\langle p_{b}\triangleright 
   a,p_{b}\rangle=\langle a,p_{b}\rangle$ signifie que $F_{A_{t}}$ 
   conserve $\phi_{a}$ alors $F_{A_{t}}$  est fidèle puisque $\phi_{a}$ est fidèle.
   
   Grâce à la coassociativité du coproduit, on a :
   $$F_{A_{s}}F_{A_{t}}=(\phi_{a}\otimes \id \otimes 
   \phi_{a})(\Delta_{a}\otimes \id)\Delta_{a}=(\phi_{a}\otimes \id \otimes 
   \phi_{a})(\id\otimes \Delta_{a})\Delta_{a}=F_{A_{t}}F_{A_{s}}$$
   \end{proof}
   
   \subsubsection{Automorphisme modulaire de $\phi_{a}$} \label{gs}
 \begin{prop}[BNSz - 4.12 et 4.14]
 Les  éléments $g_{s}=F_{A_{s}}(p_{a})^{1/2}$ et 
 $g_{t}=F_{A_{t}}(p_{a})^{1/2}$ sont inversibles et l'automorphisme 
 modulaire de $\phi_{a}$ est implémenté par $g_{s}g_{t}$ 
 c'est-à-dire que $\phi_{a}(g_{s}^{-1}g_{t}^{-1}.)$ est une trace 
 sur $A$. On posera de même 
 $$\hat{g}_{s}=F_{B_{s}}(p_{b})^{1/2} \quad \text{ et } \quad
 \hat{g}_{t}=F_{B_{t}}(p_{b})^{1/2}.$$
 \end{prop}

 D'après [BNSz - 4.13], on a les formules suivantes :
 \begin{align*}
    &  \hat{g}_{t} =  1_{b} \triangleleft g_{s}  =  1_{b} \triangleleft 
     g_{t} & \qquad  & g_{t}  = 1_{a} \triangleleft \hat{g}_{s}  =  
     1_{a} \triangleleft 
     \hat{g}_{t} \\
   &  \hat{g}_{s}  =  g_{s} \triangleright  1_{b} =  g_{t} \triangleright 1_{b}
      & \qquad  & g_{s}  =  \hat{g}_{s}\triangleright 1_{a}  = \hat{g}_{t} \triangleright 1_{a} 
      \\
     &S(g_t)  =g_s = S^{-1}(g_t)   & \qquad  & S(\hat{g}_t)  =\hat{g}_s  =S^{-1}(\hat{g}_t)
 \end{align*}
 
 Grâce à \ref{cartan} et \ref{Haar}, on en déduit la relation suivante avec les projecteurs de Haar :
 $$g_{s}p_{a}=\varepsilon_t(g_{s}) p_{a}=S(g_{s}) p_{a}=g_{t}p_{a}$$
 et de même on montre : $p_{a}g_{s}=p_{a}g_{t}$.
 
 D'après \ref{escond}, $F_{A_{s}}(g_{t}^{-1})$ appartient à $A_{s}\cap A_{t}$.
  Si A et B sont connexes, $A_{s}\cap A_{t}$ est réduit aux scalaires  
  et $F_{A_{s}}(g_{t}^{-1})$ 
 est un scalaire qui vaut $d^{-1}\phi_{a}(g_{t}^{-1})$ avec 
 $d=\phi_{a}(1_{a})$.
 
  Les mesures de Haar sont  invariantes par les antipodes et coïncident 
 sur les algèbres co-unitales donc on peut écrire les égalités 
 suivantes 
 $$\phi_{a}(g_{s}^{-1})=\phi_{a}(g_{t}^{-1})=
 \phi_{b}(\hat{g}_{s}^{-1})=\phi_{b}(\hat{g}_{s}^{-1})$$
 et on notera $\gamma$ ces scalaires. On en déduit :
 $$F_{A_{s}}(g_{t}^{-1})=F_{A_{t}}(g_{s}^{-1})=F_{B_{s}}(\hat{g}_{t}^{-1})
=F_{B_{t}}(\hat{g}_{s}^{-1})
 = d^{-1}\gamma$$
 
 Si la restriction de $S$ aux algèbres co-unitales est involutive alors 
  la restriction de $\phi_{a}$ à $A_{t}$ est la trace canonique de 
  $A_{t}$ (voir 
  \ref{phiAsBt}) et 
  $d$ est la dimension commune des algèbres co-unitales.
 
 \subsection{Trace sur $A.B$}   
 \subsubsection{Traces sur $A$ et $B$}
 \begin{defn}	Pour $a$ de $A$, on définit en posant  :
 $$tr_{a}(a)= d \gamma^{-2}\phi_{a}(g_{s}^{-1}g_{t}^{-1}a)$$
 une trace normalisée $tr_{a}$ sur $A$. 
 
 On définit de même $tr_{b}$ .
 \end{defn}
 
 \subsubsection{Espérances conditionnelles} \label{escond2} On notera
  $E_{A_{s}}$, $E_{B_{t}}$ etc \ldots les espérances conditionnelles 
 définies par les traces $tr_{a}$ et $tr_{b}$. Elles sont reliées aux espérances 
 définies par $\phi_{a}$ et $\phi_{b}$ par les formules :
 \begin{eqnarray*}
 E_{A_{s}}(a)	 & =d\gamma^{-1}F_{A_{s}}(ag_{t}^{-1}) &\quad (a \in A)  \\
 E_{B_{t}}(b)	 & =d\gamma^{-1}F_{B_{t}}(\hat{g}_{s}^{-1}b) &\quad (b \in B)
 \end{eqnarray*}
 
  \subsubsection{Lemme} \label{lem}
   {\sl Pour tout $b$ de $B$ et tout $y$ de $B_{s}$, on a :
   \begin{itemize}
  \item[(i)]  $b_{(1)}\otimes E_{B_{t}}(yb_{(2)}) 
   =\Delta_{b}(E_{B_{t}}(yb))$
   \item [(ii)] $b^{*}_{(1)}\otimes E_{B_{t}}(yb_{(2)}^{*}) 
   =\Delta_{b}(E_{B_{t}}(yb^{*}))$
   \item [(iii)] $S^{-1}_{b}E_{B_{t}}(b) \triangleright 1_{a}
   	=E_{B_{t}} (b)\triangleright 1_{a}$
   \end{itemize}}
      \dem

   La définition de $E_{B_{t}}$ et les propriétés de $g_{s}$ 
   permettent
   d'obtenir facilement les propriétés de $E_{B_{t}}$ à partir de 
   celles de $F_{B_{t}}$.
   \begin{itemize}
   	\item  [(i)]Comme $\Delta_{b}(y)$ égale $(1_{b} \otimes y)\Delta_{b}(1)$ on 
   peut écrire :
   \begin{eqnarray*}
    b_{(1)}\otimes F_{B_{t}}(yb_{(2)})	 &=  &(\id \otimes \id \otimes 
   \phi_{b})(\id \otimes \Delta_{b})\Delta_{b}(yb) \\
   	 & = & (\id \otimes \id \otimes 
   \phi_{b})(\Delta_{b}\otimes \id)\Delta_{b}(yb)    \\
   	 & = & \Delta_{b}(F_{B_{t}}(yb))
   \end{eqnarray*}
   
   	\item [(ii)]  se démontre de même.
   
   	\item [(iii)] A l'aide de \ref{Haar} et  \ref{cartan}, on obtient :
   \begin{align*}
   S_{b}^{-1}F_{B_{t}}(b) \triangleright 1_{a}&=1_{a(1)}\langle 
   1_{a(2)},S_{b}^{-1}b_{(1)}\rangle\langle p_{a},b_{(2)} \rangle 
   &&=1_{a(1)}\langle S_{a}^{-1}(1_{a(2)}),b_{(1)}\rangle\langle p_{a},b_{(2)} 
   \rangle \\
   &=1_{a(1)}\langle S_{a}^{-1}(1_{a(2)})p_{a},b\rangle 
   &&=1_{a(1)}\langle p_{a}1_{a(2)},S_{b}^{-1}b\rangle \\
   &=1_{a(1)}\langle p_{a}\varepsilon_{s}(1_{a(2)}),S_{b}^{-1}b\rangle 
   &&=1_{a(1)}\langle 
   p_{a}\varepsilon_{s}\circ\varepsilon_{t}(1_{a(2)}),S_{b}^{-1}b\rangle \\
   &=1_{a(1)}\langle p_{a}S_{a}\circ \varepsilon_{t}(1_{a(2)}),S_{b}^{-1}b\rangle 
   &&=1_{a(1)}\langle p_{a} S_{a}(1_{a(2)}),S_{b}^{-1}b\rangle \\
   &=1_{a(1)}\langle 1_{a(2)}p_{a} ,b\rangle 
   &&=F_{B_{t}}(b) \triangleright 1_{a}&&
   \end{align*}
 \end{itemize}
\end{proof}
   
   \subsubsection{D'autres espérances conditionnelles} 
   La proposition 4.2  de [N] se généralise ainsi~:
  \begin{prop}\ 

  \begin{enumerate}
   	\item  
    En posant pour tout $[a \otimes b]$ dans $A.B$,
  $$E_{A}([a \otimes b])=a(E_{B_{t}}(b) \triangleright 1_{a})$$
  on définit une espérance conditionnelle fidèle de $A.B$ dans $A$.
  
  On définit de même l'espérance conditionnelle $E_{B}$.
  \item Le carré $\mathcal {C}$ 
  $$\begin{matrix}
A & \raisebox{1.3 ex}{$\begin{matrix} E_{A}\\ \subset \end{matrix}$}  & A.B  \\
 \\
 \cup &
 & \;\;\cup \; E_{B}\\
 \\
A_{s}=B_{t}&  \subset   &B
\end{matrix}$$
  est commutatif 
  et symétrique (voir [JS - 5.3.6]). L'algèbre $A \cap B$ est 
  $A_{s}=B_{t}$.
   \end{enumerate}	
\end{prop}
   \dem
   \begin{enumerate}\item Etudions les propriétés de $E_{A}$ :
   \begin{itemize}
   	\item  
   Grâce à \ref{AsBt} la définition de $E_A$ ne dépend pas du 
   représentant de $[a \otimes b]$.
   	\item  
   Soient $a$ et $\alpha$ dans $A$ et $b$ dans $B$, comme $[\alpha\otimes 
   1_{b}][a \otimes b]$ vaut $[\alpha a \otimes b]$, on a bien :
   $$E_{A}([\alpha \otimes 1_{b}][a \otimes b])=\alpha E_{A}([a \otimes 
   b]).$$
   Calculons maintenant $E_{A}([a \otimes b][\alpha \otimes 1_{b}])$.
   $$E_{A}([a \otimes b][\alpha \otimes 
   1_{b}])=a(b_{(1)}\triangleright \alpha)(E_{B_{t}}(b_{(2)}) \triangleright 
   1_{a})$$
   Grâce à \ref{lem}(i) et à \ref{rem27}, on obtient :
  $$E_{A}([a \otimes b][\alpha \otimes 1_{b}])=
   a (E_{B_{t}}(b)\triangleright 
  \alpha)=a (E_{B_{t}}(b)\triangleright 1_{a}) \alpha= E_{A}([a \otimes 
  b])\alpha$$
      	\item  
   Vérifions l'égalité $E_{A}(x^{*})=E_{A}(x)^{*}$. Grâce à 
   \ref{lem}(ii), on peut écrire :
  $$E_{A}([a \otimes b]^{*})=E_{A}([b_{(1)}^{*}\triangleright a^{*}\otimes 
  b_{(2)}^{*}])
  =(b_{(1)}^{*}\triangleright a^{*}) (E_{B_{t}}(b_{(2)}^{*})\triangleright 
  1_{a})=E_{B_{t}}(b^{*})\triangleright a^{*}$$
  Comme $S_{b}$  envoie $B_{t}$ sur $B_{s}$, 
 on obtient grâce à \ref{rem27} :
  $$E_{A}([a \otimes b]^{*})=(S_{b}^{-1}E_{B_{t}}(b)\triangleright a)^{*}
  =( a (S_{b}^{-1}E_{B_{t}}(b)\triangleright 1_{a}))^{*}$$
   et \ref{lem}(iii)  permet de conclure.
  \item  
  Soit $\{u_{p}, p \in P\}$ une quasi-base de $B$ sur 
  $B_{t}$ (voir \ref{qb}). Alors un élément $[a\otimes b]$ de $A.B$ 
  s'écrit :
  \begin{align*}
  [a\otimes b]&= \sum_{p \in P}[a \otimes 
  E_{B_{t}}(bu_{p})u_{p}^{*}]\\
  &= \sum_{p \in P}[a (E_{B_{t}}(bu_{p} )\triangleright  1_{a})\otimes 
  1_{b}][1_{a} \otimes u_{p}^{*}]\\
  &= \sum_{p \in P}[E_{A}([a\otimes b][1_{a}\otimes u_{p} ]) 
  \otimes 1_{b}][1_{a} \otimes u_{p}^{*}]
  \end{align*}
  Par involution, nous obtenons :
  $$[a\otimes b]= \sum_{p \in P}[1_{a} \otimes u_{p}]
  [E_{A}([1_{a}\otimes  u_{p}^{*}][a\otimes b]) 
  \otimes 1_{b}]$$
  c'est-à-dire pour tout $x$ de $A.B$
  $$x=\sum_{p \in P}[1_{a} \otimes u_{p}]
  [E_{A}([1_{a}\otimes  u_{p}^{*}]x) 
  \otimes 1_{b}]$$
  \item  
  Nous calculons maintenant $E_{A}(x^{*}x)$.
  \begin{align*}
  E_{A}(x^{*}x)&=\sum_{(p ,q)\in P\times P}E_{A}(x^{*}[1_{a}\otimes 
  u_{q}]) 
  (E_{B_{t}}(u_{q}^{*} u_{p})\triangleright 1_{a})
  E_{A}([1_{a}\otimes  u_{p}^{*}]x) \\
  &=\sum_{(p ,q)\in P\times P}E_{A}(x^{*}[1_{a}\otimes 
  u_{q}][1_{a}\otimes E_{B_{t}}(u_{q}^{*} u_{p})]) 
  E_{A}([1_{a}\otimes  u_{p}^{*}]x) \\
  &=\sum_{p \in P}E_{A}(x^{*}[1_{a}\otimes 
  u_{p}]) 
  E_{A}([1_{a}\otimes  u_{p}^{*}]x) \\
  &=\sum_{p \in P}E_{A}([1_{a}\otimes  u_{p}^{*}]x)^{*}
  E_{A}([1_{a}\otimes  u_{p}^{*}]x) 	
  \end{align*}
  On en déduit que $E_{A}$ est positive et fidèle.
  \end{itemize}
  \item Le carré $\mathcal {C}$ est commutatif puisqu'on peut écrire
   $$E_{A}E_{B}=E_{A_{s}}\otimes E_{B_{t}}=E_{B}E_{A}.$$
   On en déduit que $A_{s}=B_{t}=A \cap B$
  Il est symétrique par définition puisque $A.B$ est l'espace 
  vectoriel engendré par les produits $i_{a}(a)i_{b}(b)$ ($a\in A, b\in B$).
  \end{enumerate}
  \end{proof}
  
  \begin{cor}
    	Pour tout $[a\otimes  b]$ de $A.B$, on a l'égalité :
 $$tr_{a}(E_{A}([a\otimes  b]))=tr_{b}(E_{B}([a\otimes  b]))$$
    \end{cor}
    
    \dem

  Soit $[a \otimes b] \in A.B$.
 $$tr_{a}(E_{A}([a\otimes  b]))=tr_{a}(a(E_{B_{t}}(b)\triangleright 
 1_{a}))=tr_{a}(E_{A_{s}}(a)(E_{B_{t}}(b)\triangleright 1_{a}))$$
 Par un calcul analogue on trouve :
$$ tr_{b}(E_{B}([a\otimes  b]))=tr_{b}
((1_{b}\triangleleft E_{A_{s}}(a)) E_{B_{t}}(b))$$
Par \ref{AsBt}, \ref{gs} et \ref{phiAsBt}, $tr_{a}$ et $tr_{b}$ 
coïncident sur $A_{s}=B_{t}$. On conclut à l'égalité :
$$tr_{a}(E_{A}([a\otimes  b]))=tr_{b}(E_{B}([a\otimes b]))$$ \end{proof}  
  
  \subsubsection{Prolongement des traces $tr_{a}$ et $tr_{b}$ à $A.B$}\ 
 Le corollaire précédent permet de prolonger $tr_{a}$ et $tr_{b}$ à 
 $A.B$. Etudions ce prolongement.
 \begin{prop}
 La formule  
 $$tr([a\otimes b])=tr_{a}(E_{A}([a\otimes  b]))=tr_{b}(E_{B}([a\otimes  
 b]))$$
 définit une trace normalisée fidèle sur $A.B$. 
 Par construction, les 
 espérances $E_{A}$ et $E_{B}$ conservent cette trace.
 On a aussi : $$tr([a\otimes b])=tr_{a}(E_{B_{t}}(b) \triangleright a)=
 tr_{b}(b\triangleleft E_{A_{s}}(a))$$
  \end{prop}
 \dem
 
 On vérifie facilement l'égalité : $tr([1_{a}\otimes 1_{b}])=1$. De 
 plus $tr$ est une forme linéaire positive et fidèle puisque $E_{A}$ est 
 linéaire, positive et fidèle.
 
 Comme $[x \otimes y]$ vaut $[x\otimes 1_{b}][1_{a}\otimes y]$ pour 
 montrer que $tr$ est une trace, il suffit de montrer, pour $a$ et $x$ 
 dans $A$ et $b$ et $y$ dans $B$, les deux 
 identités :
 \begin{enumerate}
\item  $tr([a \otimes b][x\otimes 1_{b}])=tr([x\otimes 1_{b}][a \otimes b])$
 \item  $tr([a \otimes b][1_{a}\otimes y])=tr([1_{a}\otimes y][a \otimes b])$
 \end{enumerate}
 Montrons la première, la seconde se démontre de manière analogue.
 Grâce à \ref{lem}(i), on a :
 \begin{eqnarray*}
 tr([a \otimes b][x\otimes 1_{b}])& = & tr([a(b_{(1)}\triangleright x) 
 \otimes b_{(2)}])  \\
 	 & = & tr_a\big (a (b_{(1)}\triangleright x) 
 (E_{B_{t}}(b_{(2)}) \triangleright 1_{a} )\big)   \\
 & = & tr_a\big (a (E_{B_{t}}(b) \triangleright x)\big)   
 \end{eqnarray*}
 Comme $tr_a$ est une trace, à l'aide des formules \ref{rem27} 
 on peut écrire :
 \begin{eqnarray*}
tr([a \otimes b][x\otimes 1_{b}])
&=&tr_a(a(E_{B_{t}}(b) \triangleright 1_{a}) x)\\
&=&tr_a(xa (E_{B_{t}}(b) \triangleright 1_{a}))\\
&=&tr([x\otimes 1_{b}][a \otimes b])\end{eqnarray*}
De plus,  pour tout $[a\otimes b]$ de $A.B$, on a l'égalité :
 $$tr([a\otimes b])
 =tr_a((E_{B_{t}}(b) \triangleright 1_{a})a) =tr_{a}(E_{B_{t}}(b) \triangleright a).$$
 \end{proof}

\subsubsection{Représentation standard de 
$A.B$ sur $L^2(A,tr)$}\label{GNSa}
\begin{prop}[BSz2 - 4.2]
L'algèbre $A.B$ admet une représentation fidèle $\pi_{\phi}$ sur $L^2(A,\phi_{a})$ qui 
prolonge la représentation de $A$ par multiplication à gauche. En 
particulier, $\pi_{\phi}$ vérifie pour $b$ dans $B$ et $a$ dans $A$, 
$$\pi_{\phi}(b)\Lambda_{\phi}(a)=\Lambda_{\phi}(b \triangleright 
a).$$

L'algèbre $A.B$ est l'extension de Jones de $A_{t} \subset A$ 
représentée sur $L^2(A,\phi_{a})$. Elle 
est donc engendrée par $A$ et $p_{b}$, projecteur de Jones de 
l'inclusion. Plus précisément, $p_{b}$ vérifie :
$\pi_{\phi}(p_{b})\Lambda_{\phi}(a)=\Lambda_{\phi}(F_{A_{t}}(a))$.
\end{prop}

Considérons l'isométrie $U$ de $L^2(A,\phi_{a})$ sur $L^2(A,tr)$ définie 
par 
$$U\Lambda_{\phi}(a)=d^{-1/2}\gamma\Lambda_{tr}(ag_{s}^{1/2}g_{t}^{1/2}).$$
On vérifie facilement que la représentation  $\pi = U\pi_{\phi}U^{-1}$ de $A.B$ prolonge la 
représentation standard de $A$ sur $L^2(A,tr)$ et le projecteur de 
Jones de l'inclusion $A_{t} \subset A$ représentée sur $L^2(A,tr)$ est alors 
$$f_{b}=d\gamma^{-1} 
\hat{g}_{t}^{-1/2}p_{b} \hat{g}_{t}^{-1/2}.$$
en effet à l'aide de \ref{gs} et \ref{rem27}, on obtient pour $a$ 
dans $A$:
\begin{eqnarray*}
\pi(f_{b})\Lambda_{tr}(a) & = 
& d^{3/2}\gamma^{-2}U\Lambda_{\phi}(\hat{g}_{t}^{-1/2}\triangleright 
F_{A_{t}}(\hat{g}_{t}^{-1/2}\triangleright ag_{s}^{-1/2}g_{t}^{-1/2}))  \\
&  =& d^{3/2}\gamma^{-2}U\Lambda_{\phi}(g_{s}^{-1/2}
F_{A_{t}}(g_{s}^{-1/2} ag_{s}^{-1/2}g_{t}^{-1/2}))    \\
	 & = &d\gamma^{-1}\Lambda_{tr}(g_{s}^{-1/2}
F_{A_{t}}(g_{s}^{-1/2} 
ag_{s}^{-1/2}g_{t}^{-1/2})g_{s}^{1/2}g_{t}^{1/2})\\
&=&d\gamma^{-1}\Lambda_{tr}(
F_{A_{t}}(g_{s}^{-1/2} 
ag_{s}^{-1/2}))\\
&=&\Lambda_{tr}(E_{A_{t}}(a))
\end{eqnarray*}

On vérifie facilement que $E_{B_{t}}(f_{b})$ vaut $d^2\gamma^{-2}$.

\begin{prop}
L'algèbre $A.B$ est l'extension de Jones de $A_{t} \subset A$ 
représentée sur $L^2(A,tr_{a})$. Le projecteur de 
Jones est $f_{b}=d\gamma^{-1} \hat{g}_{t}^{-1/2}p_{b} 
\hat{g}_{t}^{-1/2}$ avec  $E_{B_{t}}(f_{b}) = d^2\gamma^{-2}$.	
\end{prop}

\subsubsection{Trace de Markov}
\begin{prop}
 La trace $tr$ est la trace de Markov normalisée de l'inclusion 
 $A_{t}\; 
 \raisebox{1.3 ex}{$\begin{matrix} {\scriptstyle E_{{\scriptscriptstyle A_{t}}}}
 \\ \subset \end{matrix}$}\;A$
  dont l'indice est $d^{-2}\gamma^2$. C'est aussi la trace de Markov 
  de l'inclusion $A\;
 \raisebox{1.3 ex}{$\begin{matrix} {\scriptstyle E_{{\scriptscriptstyle 
 A}}}\\ \subset \end{matrix}$}\;A.B$
\end{prop}

\dem

Comme l'inclusion $A_{t} \subset A$ est connexe, il existe une unique trace de Markov 
dont le module est l'indice de l'inclusion ([GHJ - 2.7.3]).
D'après \ref{escond2} et \ref{GNSa}, on a pour tout $x$ de $A_{t}$ :
$$tr([x\otimes f_{b}]=tr(E_{B_{t}}(f_{b}) \triangleright x)=d^2\gamma^{-2} 
tr([x\otimes 1_{b}] )$$ 
On en déduit que  $tr$ est la trace de Markov de module 
$\beta=d^{-2}\gamma^2$ ([GHJ - 2.7.1]). D'après [GHJ - 2.7.4], $tr$ est aussi la trace de de l'inclusion $A\;
 \raisebox{1.3 ex}{$\begin{matrix} {\scriptstyle E_{{\scriptscriptstyle 
 A}}}\\ \subset \end{matrix}$}\;A.B$.
\end{proof}

\subsection{L'inclusion $M_1 \subset M_2$}

\subsubsection{Construction de $M_1\subset M_2$}\label{R0}
Comme le carré $\mathcal {C}$ est un carré commutatif symétrique
(\ref{escond2}) pour les espérances conditionnelles associées à la  trace
de Markov de l'inclusion $A\subset A.B$,
il  vérifie le corollaire 5.3.4 de [JS] : La trace $tr$ est aussi la
trace de Markov
     des inclusions $A_{s} \subset A$, $B_{t}\subset B$ et $B \subset
A.B$ et on obtient  par  construction de
base une échelle périodique de carrés commutatifs (voir aussi [JS
5.3.5]) qu'on peut préciser grâce à \ref{415}.

$$  \begin{matrix}
A_{0}=A &  \subset   & A_{1}=A\rtimes B &\subset   & A_{2}=A\rtimes
B \rtimes  A &\subset   & A_{3}=A\rtimes B \rtimes A \rtimes B&\dots\\

   \cup &
   & \cup &&\cup &
   & \cup \\

B_{0}=B_{t}&  \subset  &B_{1}=B&\subset   & B_{2}=
B \rtimes  A &\subset   & B_{3}= B \rtimes A \rtimes B&\dots
\end{matrix}$$
     A la limite, on obtient une paire $M_1\subset M_2$ de facteurs
     hyperfinis de type $II_{1}$ : $M_{2}$  (resp. $M_1$) est la fermeture
     faible de $ \cup_{n\in \mathbb{N}} A_{n}$ (resp. $\cup _{n\in
\mathbb{N}} B_{n}$) dans la construction GNS par rapport à la
     trace $tr$ qui se prolonge.

     \subsubsection{} \label{At} De plus d'après [JS - 5.7.1],
     le commutant relatif $M_1' \cap M_2$ est $B' \cap A$. On en déduit la
     proposition suivante :
     \begin{prop} Le commutant relatif $M'_1 \cap M_2$ est l'algèbre $A_{t}$
     contenue dans $A=A_{0}$.
     \end{prop}

     \dem
     
     Le calcul suivant montre que les éléments de $A_{t}$ commutent à
     $B$. Soient $x\in A_{t}$ et $b\in B$, alors on a :
    $$[1_{a} \otimes b][x\otimes 1_{b}]=[ x_{(1)} \otimes b \triangleleft
    x_{(2)}]$$
    Grâce à \ref{cartan}, on en déduit  :
    $$[1_{a} \otimes b][x\otimes 1_{b}]
     =[x 1_{a(1)} \otimes b
     \triangleleft 1_{a(2)}]
     =[x  \otimes (1_{b} \triangleleft 1_{a(1)})(b
     \triangleleft 1_{a(2)})]=[x\otimes b
     \triangleleft 1_{a}] = [x\otimes 1_{b}][1_{a} \otimes b]$$
        D'autre part, si un élément $x$ de $A$ commute à $B$, il commute à
     $f_{b}$ et on a : $$xf_{b}=E_{A_{t}}(x)f_{b}.$$
     Or d'après [GHJ - 2.6.7 (iii)], pour tout $y$ de $A.B$, il existe
     un unique $z$ dans $A$ tel que $ yf_{b}=zf_{b}$. On obtient donc
      l'égalité $x=E_{A_{t}}(x)$ et $x$ appartient à $A_t$.
    \end{proof}

\subsection{Action de $A$ sur $M_1$}
\subsubsection{} Comme en \ref{AB} et [N - 5.6], on peut transformer
les produits croisés à gauche en produits croisés à droite et obtenir
l'échelle de carrés commutatifs suivante isomorphe à celle considérée
en \ref{R0} :
$$  \begin{matrix}
A &  \subset   & A\ltimes B &\subset   & A\ltimes
B \ltimes  A &\subset   & A\ltimes B \ltimes A \ltimes B&\dots&
\subset & M_2\\
   \cup &
   & \cup &&\cup &
   & \cup&&\cup \\

B_{t}&  \subset  &B&\subset   &
B \ltimes  A &\subset   &  B \ltimes A \ltimes B&\dots&\subset&M_1
\end{matrix}$$

\subsubsection{} Nous généralisons maintenant la proposition 5.7 de [N].
\begin{prop}
Soient $i_{a} : a\mapsto[a\otimes 1_{b} \otimes 1_{a}\otimes\dots]$ l'inclusion
de $A$ dans $M_2$ et $E_{M_1}$ l'espérance conditionnelle de $M_2$ sur $M_1$
conservant la trace $tr$. Posons
$$f_{a}=d\gamma^{-1} g_{s}^{-1/2}p_{a} g_{s}^{-1/2}.$$
L'application de $M_1 \times A$ dans $M_1$
$$x \triangleleft a = d^{-2}\gamma^2
E_{M_1}(i_{a}(g_{s}^{-1/2}f_{a}g_{s}^{1/2}) xi_{a}(a))\quad(x\in M_1,
a\in A)$$
définit une action extérieure à droite de $A$ sur $M_1$ telle que $M_2=A\ltimes M_1$.
\end{prop}
   \dem
   
    Notons $x=[b \otimes z]$ un élément  de $B_{n}$
   avec $b \in B$ et abusivement $z \in A_{n-1}$. L'élément $x$ vu dans
   $A_{n} \subset M_2$ s'écrit $[1_{a} \otimes b \otimes z]$ et $[1_{a} \otimes
   1_{b} \otimes z]$ commute
   avec $A_{0}$.
   Observons l'action de $A$ sur $x$ :
   \begin{align*}
   x \triangleleft a & =  d^{-1}\gamma
E_{M_1}(i_{a}(g_{s}^{-1}p_{a})xi_{a}(a))  \\
   & =  d^{-1}\gamma
E_{M_1}([g_{s}^{-1}p_{a}\otimes b \otimes 1][a\otimes 1_{b}\otimes  z])\\
    & =  d^{-1}\gamma
E_{M_1}([g_{s}^{-1}p_{a}a_{(1)}\otimes b \triangleleft a_{(2)}\otimes  z])\\
& =  d^{-1}\gamma
E_{M_1}([g_{s}^{-1}p_{a}\varepsilon_{s}(a_{(1)})\otimes b \triangleleft a_{(2)}
\otimes  z])&\text{(d'après \ref{Haar})}\\
   & =  d^{-1}\gamma
E_{M_1}([g_{s}^{-1}p_{a}\otimes (1_{b}\triangleleft\varepsilon_{s}(a_{(1)}))
(b \triangleleft a_{(2)})\otimes  z])\\
& =  d^{-1}\gamma
E_{M_1}([g_{s}^{-1}p_{a}\otimes (1_{b}\triangleleft a_{(1)})
(b \triangleleft a_{(2)})\otimes  z])\\
& =  d^{-1}\gamma
E_{M_1}([g_{s}^{-1}p_{a}\otimes b \triangleleft a\otimes  z])
\end{align*}
Le carré
$  \begin{matrix}
   A_{n} &\subset   & M_2\\
   \cup && \cup &\\
   B_{n}&  \subset  &M_1
\end{matrix}$ muni des espérances conditionnelles conservant $tr$ est
commutatif,  on en déduit :
\begin{align*}
x \triangleleft a &=  d^{-1}\gamma
E_{B_{n}}([g_{s}^{-1}p_{a}\otimes b \triangleleft a\otimes  z])\\
&=  d^{-1}\gamma
[(1_{b}\triangleleft E_{A_{s}}(g_{s}^{-1}p_{a})) (b \triangleleft a)\otimes
z]\\&=[b \triangleleft a\otimes  z]
	\end{align*}
En effet d'après \ref{gs} et \ref{escond2}, on a :
$$E_{A_{s}}(g_{s}^{-1}p_{a})=d\gamma^{-1}g_{s}^{-1}F_{A_{s}}(p_{a}g_{t}^{-1})
=d\gamma^{-1}g_{s}^{-1}F_{A_{s}}(p_{a}g_{s}^{-1})
=d\gamma^{-1}g_{s}^{-1}F_{A_{s}}(p_{a})g_{s}^{-1}=d\gamma^{-1}$$
L'application $x \mapsto x \triangleleft a$ prolonge donc l'action
duale $[b\otimes z] \mapsto [b \triangleleft a\otimes  z]$ de $A$ sur
$B_{n}=B \ltimes A_{n-1}$. Elle définit une action à droite
faiblement continue de $A$ sur $M_1$. De plus $A\ltimes M_1 = i_{a}(A)M_1$
est le facteur $M_2$.
L'action est extérieure d'après \ref{ext} et \ref{At}.
   \end{proof}
   
   %%%%%%%%%%%%%%%%%%%%%%%%%%%%%%%%%%%

\subsection{Conclusion} \label{conclusion}
On a donc montré que le C*-groupoïde quantique connexe   fini
 $A$ agit exté\-rieurement sur le facteur hyperfini de type
$\mathrm{II}_{1}$. On peut voir $A$ comme le commutant relatif d'une 
inclusion en complétant l'échelle de carrés commutatifs :

$$  \begin{matrix}
B &  \subset   & B\ltimes A &\subset   & B\ltimes
A \ltimes  B &\subset   & B\ltimes A \ltimes B \ltimes A&\dots&
\subset & M_1\\
   \cup &
   & \cup &&\cup &
   & \cup&&&\cup \\

B_{s}=A_{t}&  \subset  &A&\subset   &
A \ltimes  B &\subset   &  A \ltimes B \ltimes A&\dots&\subset&M_0
\end{matrix}$$
par construction de base verticale car le projecteur de Jones est 
$e_{1}=f_{a} \otimes 1 \otimes 1 \otimes \dots$ pour chaque colonne :

$$ \begin{matrix}
C_{0}=A \ltimes B &  \subset   & C_{1}=A \ltimes  B\ltimes A &\subset   & 
C_{2}=A \ltimes B\ltimes
A \ltimes  B &\subset   & \dots&
\subset  &M_2\\
   \cup &
   & \cup &&\cup &
   & &&\cup \\
B &  \subset   & B\ltimes A &\subset   & B\ltimes
A \ltimes  B &\subset   &\dots&
\subset  &M_1\\
   \cup &
   & \cup &&\cup &
   & &&\cup \\
A_{t}&  \subset  &A&\subset   &
A \ltimes  B &\subset   &  &\subset &M_0
\end{matrix}$$
Le facteur $M_{2}=\cup_{n \in \mathbb{N}} C_{n}$ est engendré comme algèbre de von 
Neumann par $M_{1}$ et $e_{1}$, de plus on vérifie 
facilement l'égalité $e_{1}xe_{1}=E_{M_{0}}(x)e_{1}$ pour tout $x$ de 
$M_{1}$ puisqu'elle est vraie pour tout $x$ de $C_{n}$ et tout entier 
$n$. Alors d'après [Pipo2], la tour
 $M_{0} \,\subset \, M_1 \, \raisebox{1.3 ex}{$\begin{matrix} e_{1}\\
\subset \end{matrix}$} \,M_2$
est standard et la grille de carrés commutatifs construit vérifie les hypothèses de la proposition 5.7.5 . 

On peut construire de même le facteur $M_{3}$ avec 
$e_{2}=f_{b}\otimes 1\otimes 1\otimes 1\dots$.
On démontre alors  à l'aide du théorème 5.7.6 de [JS]
que $A$ (resp. $B$) est le commutant relatif $M'_{0} \cap M_{2}$ (resp. 
$M'_{1} \cap M_{3}$). En effet, on a :
$$M'_{0} \cap M_{2}=(1_{a}\otimes 1_{b} \otimes A)' \cap A \ltimes 
B \ltimes 1_{a}$$
 et on conclut comme en \ref{At}.
 
On peut donc  munir $A$ et $B$ en tant que commutants relatifs de 
structures duales de C*-groupoïde quantique fini régulier. Quel 
rapport entre ces nouvelles structures et celles de départ qui 
n'étaient pas nécessairement régulières ? Le calcul de la 
dualité héritée de l'inclusion $M_{0} \subset M_{1}$ en fonction de 
la dualité originelle donne la formule qui sert pour 
la déformation 
explicitée dans la partie suivante. On en déduit qu'on peut 
déformer un C*-groupoïde quantique fini en un C*-groupoïde quantique fini régulier sans modifier la structure de C*-algèbre.

%%%%%%%%%%%%%%%%%%%%%%%%%%%%%%%%%%%%%%%%%%%%%%%%
\section{Déformation régulière d'un C*-groupoïde quantique fini}

On considère deux C*-groupoïdes quantiques finis $A$ et $B$ duaux. Le 
but de cette partie est de déformer la structure de co-algèbre de 
$A$ (on peut bien sûr faire de même dans le même temps pour $B$) pour 
faire de $A$ un C*-groupoïde quantique fini régulier.
Comme il ne sera question que de la structure de 
co-algèbre de $A$, on omettra l'indice~$A$. L'étude de la structure de 
co-algèbre de $B$ est analogue.

\subsection{Eléments séparateurs} \label{separateur} D'après [NV3 - 2.3.4], $(S\otimes 
\id)(\Delta(1))$ est un élément séparateur de $A_{t}\otimes A_{t}$. 
Plus précisément, les relations utiles sont, pour $z$ dans $A_{t}$ :
\begin{eqnarray*}
\Delta(1)(S^{-1}(z)\otimes 1) & = & \Delta(1)(1\otimes z)  \\
(S^{-1}(z)\otimes 1)	\Delta(1) & = & (1\otimes z)\Delta(1)
\end{eqnarray*}

Si $\{\lambda_{i},i \in I\}$ est une famille d'unités
matricielles  de $A_{t}=\otimes_{j\in J} M_{\nu_{j}}(\mathbb
C)q_{j} $ telle que l'élément
$\lambda_{i}$ appartienne au facteur $M_{\nu_{j_{i}}}(\mathbb
C)q_{j_{i}}$, posons $q=\sum_{i\in I}\frac{1}{n_{j_{i}}} \lambda_{i}^* \otimes 
\lambda_{i}$. Alors $(S^{-1}\otimes\id)(q)$ vérifie aussi, pour $z$ dans $A_{t}$ :
\begin{eqnarray*}
(S^{-1}\otimes\id)(q)(S^{-1}(z)\otimes 1) & = & (S^{-1}\otimes\id)(q)(1\otimes z)  \\
(S^{-1}(z)\otimes 1)	(S^{-1}\otimes\id)(q) & = & (1\otimes z)(S^{-1}\otimes\id)(q)
\end{eqnarray*}

\subsection{Proposition} \label{k}
{\it
	Il existe un élément $k$ positif inversible de $A_{t}$ tel que 
\begin{enumerate}
	\item  $k^2=1_{(2)}S(1_{(1)})$

	\item  $k^2$ soit la dérivée de Radon-Nikodym de la trace canonique 
	de $A_{t}$ par rapport à la restriction de la co-unité à $A_{t}$.
	
	\item $\Delta(1) = \Delta(1) (S^{-1}\otimes\id)(q)$ 

	\item  $\Delta(1)=(1 \otimes k^2)(S^{-1}\otimes\id)(q)$
	
	\item La restriction de $S^2$ à $A_{s}A_{t}$ est 
	$\Ad(k^{-2}S(k^2))$. En particulier $S^2(k)$ vaut $k$.
\end{enumerate}
}
\dem 

Le début de l'énoncé est le lemme 4.6 de [BSz2], la 
dernière se trouve dans des notes manuscrites de K. Szlachanyi.

D'après \ref{separateur}, on peut écrire :
$$\Delta(1) (S^{-1}\otimes\id)(q)=
	\Delta(1)(1\otimes \sum_{i\in I}\frac{1}{n_{j_{i}}} 
	\lambda_{i}^*\lambda_{i})=	\Delta(1) $$
On obtient donc (3) et pour (4), on écrit
$$\Delta(1)=\Delta(1) (S^{-1}\otimes\id)(q)=(1 \otimes 
1_{(2)}S(1_{(1)}))(S^{-1}\otimes\id)(q)$$
On a aussi : $\Delta(1)= (S^{-1}(k^{2}) \otimes 
1)(S^{-1}\otimes\id)(q)$ et de $\Delta(1)=\Delta(1)^{*}$, on déduit l'égalité, pour tout $i$ de $I$ :
 $$S^{-1}(\lambda_{i}k^2)=S^{-1}(\lambda_{i}^{*}k^2)^{*}$$
 Comme $S \circ *$ est involutive et que les $\lambda_{i}$ engendrent 
 $A_{t}$ , pour tout $z$ de $A_{t}$ on a :
 $$S^{-1}(zk^2)=S(k^2z)$$
 On en déduit les deux égalités :
 $$S^2(z)=k^{-2}zk^2 \qquad S^2(S(z))=S(k^2)S(z)S(k^{-2})$$
Comme $A_{t}$ et $A_{s}$ commutent, on peut traduire les deux 
expressions de $S^2$ par (5).
\end{proof}

On trouve une proposition analogue dans [V2].

\subsection{Nouvelle dualité et déformation des co-algèbres}
On considère une nouvelle dualité entre $A$ et $B$ :
$$[a,b]=\langle kaS(k),b \rangle\qquad (a\in A,b \in B)$$
et la structure de co-algèbre définie sur $A$ par cette dualité :
 \begin{align*}
 	\tilde{\Delta}(a)&=(1\otimes k^{-1})\Delta(a)(1\otimes k^{-1})\\
 	\tilde{\varepsilon}(a)&=\varepsilon(kaS(k))=\varepsilon(S(k)ak)\\
 	\tilde{S}(a)&=S(k^{-1})kS(a)k^{-1}S(k)
 \end{align*}
 De plus, le projecteur $\tilde{\Delta}(1)$ vaut 
 $(\tilde{S}^{-1}\otimes\id)(q)$.
  
 \begin{thm}
 	L'algèbre $A$ munie de sa structure d'algèbre originelle et de 
 	cette nouvelle structure de co-algèbre est un C*-groupoïde quantique 
 	fini régulier.
 \end{thm}
 \dem
 
 Les nouveaux co-produit, co-unité et antipode étant définis par 
 dualité, il suffit de vérifier les propriétés propres aux C*-groupoïdes 
 quantiques.

  Vérifions que le nouveau co-produit est un homomorphisme d'algèbres 
 involutives. D'après \ref{k}, pour tous $x$ et $y$ dans $A$, on a :
  \begin{align*}
 \tilde{\Delta}(x)\tilde{\Delta}(y)
 &=(1\otimes k^{-1})\Delta(x)(1\otimes k^{-2})\Delta(y)(1\otimes k^{-1})\\
 &=(1\otimes k^{-1})\Delta(x)\Delta(1)(1\otimes 
 k^{-2})\Delta(1)\Delta(y)(1\otimes k^{-1})\\
 &=(1\otimes k^{-1})\Delta(x)\Delta(1)(S^{-1}\otimes\id)(q)\Delta(y)(1\otimes 
 k^{-1}) \\
 &=\tilde{\Delta}(xy)	
 \end{align*}
 L'égalité $\tilde{\Delta}(x)^{*}=\tilde{\Delta}(x^{*})$ est 
 évidente.
 
 Le projecteur $\tilde{\Delta(1)}$ vérifie l'égalité :
 $$(\tilde{\Delta}\otimes 1)\tilde{\Delta}(1)=(1\otimes 
 \tilde{\Delta}(1))(\tilde{\Delta}(1)\otimes 1)$$
 En effet, comme $k$ commute à $A_{s}$, on a :
 \begin{align*}
 	(1\otimes \tilde{\Delta}(1))(\tilde{\Delta}(1)\otimes1)
 	&=(1\otimes 1 \otimes k^{-1})(1\otimes \Delta(1))(1\otimes k^{-1} \otimes k^{-1})( \Delta(1)
 	\otimes1)(1\otimes k^{-1} \otimes 1)\\
 &=(1\otimes k^{-1} \otimes k^{-1})(1\otimes \Delta(1))( \Delta(1)
 	\otimes1)(1\otimes k^{-1} \otimes k^{-1})\\
 &=(1\otimes k^{-1} \otimes k^{-1})(\Delta \otimes \id)\Delta(1)(1\otimes 
 k^{-1} \otimes k^{-1})\\
 &= (\tilde{\Delta}\otimes 1)\tilde{\Delta}(1)
 \end{align*}

 On pose $\tilde{\Delta}(y)=\tilde{y}_{(1)} \otimes \tilde{y}_{(2)}$ 
 et on démontre sans problème la relation :
 $$\tilde{\varepsilon}(xyz)
=\tilde{\varepsilon}(x\tilde{y}_{(1)})\tilde{\varepsilon}(\tilde{y}_{(2)}z)\qquad ((x,y,z)\in A^3)$$
 
Vérifions maintenant la relation entre $\tilde{\Delta}$, $\tilde{\varepsilon}$ et
$\tilde{S}$  pour $a$ dans $A$ :
\begin{align*}
(\tilde{\varepsilon}\otimes \id)(\tilde{\Delta}(1)(a\otimes 1))
&=(\varepsilon\otimes \id)((S(k)\otimes k^{-1})\Delta(1)(ak\otimes k^{-1}))\\ 	
&=(\varepsilon\otimes \id)(\Delta(1)(ak\otimes 1)k^{-1}\quad (\text{d'après \ref{separateur})}\\
&=m(\id \otimes S)(\Delta(ak))k^{-1}\\
&=a_{(1)}kS(a_{(2)})k^{-1}\\
&=\tilde{a}_{(1)}S(k^2)\tilde{S}(\tilde{a}_{(2)})
\end{align*}

On a donc montré pour $a$ dans $A$ :
$$(\tilde{\varepsilon}\otimes \id)(\tilde{\Delta}(1)(a\otimes 
1)=\tilde{a}_{(1)}S(k^2)\tilde{S}(\tilde{a}_{(2)})$$

On écrit cette relation pour l'unité :
$$1=(\tilde{\varepsilon}\otimes \id)(\tilde{\Delta}(1))=\tilde{1}_{(1)}S(k^2)\tilde{S}(\tilde{1}_{(2)})$$

Or $\tilde{a}_{(1)}S(k^2)\tilde{S}(\tilde{a}_{(2)})$ vaut 
$\tilde{a}_{(1)}\tilde{1}_{(1)}S(k^2)\tilde{S}(\tilde{1}_{(2)})\tilde{S}(\tilde{a}_{(2)})$, 
on obtient donc
$$(\tilde{\varepsilon}\otimes \id)(\tilde{\Delta}(1)(a\otimes 1))
=m(\id \otimes \tilde{S})(\tilde{\Delta}(a))$$

 L'involutivité de $\tilde{S}$ sur les sous-algèbres co-unitales est 
 évidente.
 \end{proof}
 
 On peut donc déformer toute paire de C*-groupoïdes quantiques finis 
 en une paire de C*-groupoïdes quantiques finis réguliers sans modifier la structure de C*-algèbre.
%\newpage
\appendix
\section{} 
 Cet appendice est un ajout à la version publiée. Il consiste à donner une version légèrement modifiée de la construction de facteurs hyperfinis de type $II_1$ sur lesquels agissent deux C*-groupoïdes quantiques finis réguliers en dualité de façon à les retrouver, avec leurs structures initiales, comme commutants relatifs de la tour de Jones obtenue. N. Thiéry et l'auteur ont utilisé cette possibilité dans [DT]. Cette nouvelle version consiste en gros à échanger les lignes et les colonnes dans la grille de carrés commutatifs. 
 
 Nous utilisons les notations de 2.6. Nous considérons  deux C*-groupoïdes quantiques finis réguliers (voir 2.5) $A$ et $B$ en dualité.
 
\subsection{Représentation standard de $B.A$ sur $L^2(A,tr)$}\label{GNSabis}
Utilisant la remarque 6.2, nous notons $B.A$ les produits croisés isomorphes $B\rtimes A$ et $B\ltimes A$. Dans \ref{GNSa}, nous considérons la représentation $\pi$ de $A.B$ sur $L^2(A,tr)$ comme extension de Jones de $A_{t} \subset A$. Dans [BSz2 - p. 183-184], G. Böhm et K. Szlach\'anyi définissent aussi une représentation $\pi'_{\phi}$ de $B.A$ sur $L^2(A,\phi_{a})$ qui nous permet d'obtenir une représentation $\pi'$ de $B.A$ sur $L^2(A,tr)$ comme extension de Jones de $A_{s} \subset A$. La démarche est analogue à \ref{GNSa} et 6.4.7, nous résumons les résultats dans la proposition suivante :
\begin{prop} (1)
L'algèbre $B.A$ admet une représentation fidèle $\pi'_{\phi}$ sur $L^2(A,\phi_{a})$ qui 
prolonge la représentation de $A$ par multiplication à gauche. En 
particulier, $\pi'_{\phi}$ vérifie pour $b$ dans $B$ et $a$ dans $A$, 
$$\pi'_{\phi}(b)\Lambda_{\phi}(a)=\Lambda_{\phi}(a \triangleleft S_b^{-1}(b)).$$

L'algèbre $B.A$ est l'extension de Jones de $A_{s} \subset A$ 
représentée sur $L^2(A,\phi_{a})$. Elle 
est donc engendrée par $A$ et $p_{b}$, projecteur de Jones de 
l'inclusion. Plus précisément, $p_{b}$ vérifie :
$\pi_{\phi}(p_{b})\Lambda_{\phi}(a)=\Lambda_{\phi}(F_{A_{s}}(a))$.

(2) Soit $U$ l'isométrie de $L^2(A,\phi_{a})$ sur $L^2(A,tr)$ définie 
par 
$$U\Lambda_{\phi}(a)=d^{-1/2}\gamma\Lambda_{tr}(ag_{s}^{1/2}g_{t}^{1/2}).$$
La représentation  $\pi'= U\pi'_{\phi}U^{-1}$ de $B.A$ prolonge la 
représentation standard de $A$ sur $L^2(A,tr)$ et 
 $B.A$ est l'extension de Jones de $A_{s} \subset A$ représentée sur $L^2(A,tr)$. Le projecteur de 
Jones est $f_{b}=d\gamma^{-1} \hat{g}_{t}^{-1/2}p_{b} 
\hat{g}_{t}^{-1/2}$ avec  $E_{B_{s}}(f_{b}) = d^2\gamma^{-2}$.

(3) La trace $tr$ (définie sur $B.A$ comme en 6.4.5) est la trace de Markov normalisée de l'inclusion 
 $A_{s}\; 
 \raisebox{1.3 ex}{$\begin{matrix} {\scriptstyle E_{{\scriptscriptstyle A_{s}}}}
 \\ \subset \end{matrix}$}\;A$
  dont l'indice est $d^{-2}\gamma^2$. C'est aussi la trace de Markov 
  de l'inclusion $A\;
 \raisebox{1.3 ex}{$\begin{matrix} {\scriptstyle E_{{\scriptscriptstyle 
 A}}}\\ \subset \end{matrix}$}\;B.A$
\end{prop}

\subsection{Construction des facteurs}\label{grille} 
Comme en 6.5.1 et 6.7, à partir des carrés commutatifs munis de la trace de Markov prolongeant $tr_a$ et $tr_b$ :
$$\begin{matrix}
A & \raisebox{1.3 ex}{$\begin{matrix} E_{A}\\ \subset \end{matrix}$}  & A.B  \\
 \\
 \cup &
 & \;\;\cup \; E_{B}\\
 \\
A_{s}=B_{t}&  \subset   &B
\end{matrix}\qquad \qquad \text{et} \qquad \qquad \begin{matrix}
A & \raisebox{1.3 ex}{$\begin{matrix} E_{A}\\ \subset \end{matrix}$}  & B.A  \\
 \\
 \cup &
 & \;\;\cup \; E_{B}\\
 \\
A_{t}=B_{s}&  \subset   &B
\end{matrix}$$
on construit  des échelles périodiques de carrés commutatifs dont  les lignes sont des constructions de bases :
$$\begin{array}{ccccccccccc}B_1 & \subset &N^3_0= A_0.B_1 & \subset & N^3_1=B_0.A_0.B_1 & \subset & N_2^3=A.B_0.A_0.B_1 & \subset & \cdots & \subset & M_3 \\\cup &  & \cup &  & \cup &  & \cup &  &  &  & \cup \\A_s & \subset &  N^2_0=A_0 & \subset & N^2_1=B_0.A_0 & \subset & N^2_2=A.B_0.A_0 & \subset & \cdots & \subset & M_2 \\ &  & \cup &  & \cup &  & \cup &  &  &  & \cup \\ &  & B_s=A_t & \subset & N^1_1=B_0 & \subset & N^1_2=A.B_0 & \subset & \cdots & \subset & M_1 \\ &  &  &  & \cup &  & \cup &  &  &  & \cup \\ &  &  &  &N^0_1= B_t=A_s & \subset & N^0_2=A & \subset & \cdots & \subset & M_0\end{array}$$

On a noté $A_0$ et $B_0$ les copies de départ des C*-groupoïdes quantiques finis réguliers et leurs inclusions dans les algèbres construites, $N^r_t$ l'algèbre de la $r$-ème ligne et $t$-ème colonne (les lignes sont comptées à partir de celle de $M_0$ et la colonne  $-1$ est la plus à gauche). Les facteurs $M_i$ sont les facteurs engendrés par les lignes d'algèbres.  Les colonnes sont aussi des constructions de base ainsi que
 \begin{displaymath}
  M_0 \,\subset \, M_1 \, \stackrel{e_1}{\subset} \, M_2 \,\stackrel{e_2}{\subset} \, M_3 \end{displaymath} 
  dont les projecteurs sont définis dans  $N_1^3$ par $e_1=[1_b\otimes f_a \otimes 1_b]$ et $e_2=[1_b\otimes 1_a \otimes f_b]$ avec $f_{a}=d\gamma^{-1} g_{s}^{-1/2}p_{a} g_{s}^{-1/2}$.

\subsection{Commutants relatifs}\label{commutant}
On applique alors le théorème 5.7.6 de [JS] à cette grille et on en déduit  :
\begin{align*}
M'_{0} \cap M_{1}&=(N_2^0)' \cup N^1_1=(A \otimes 1_{b} )' \cap (1_{a} \otimes B_0) \subset A.B_0\\
M'_{1} \cap M_{2}&=(N_1^1)' \cup N^2_{0}=(B_0 \otimes 1_{a})' \cap (1_{b}  \otimes  A_0) \subset B_0.A_0\\
M'_{0} \cap M_{2}&=(N_2^0)' \cup N^2_1=(A \otimes 1_{b} \otimes 1_{a})' \cap (1_{a} \otimes B_0.A_0) \subset A.B_0.A_0\\
M'_{1} \cap M_{3}&=(N_1^1)' \cup N^3_{0}=(B_0 \otimes 1_{a} \otimes 1_{b})' \cap ( A_0.B_1  \otimes 1_{b}) \subset B_0.A_0.B_1\\
M'_{0} \cap M_{3}&=(N_2^0)' \cup N^3_1=(A \otimes 1_{b} \otimes 1_{a} \otimes 1_{b})' \cap (1_{a} \otimes B_0.A_0.B_1) \subset A.B_0.A_0.B_1
\end{align*}
Comme en \ref{At}, on conclut que le carré commutatif $\;\;\begin{matrix}
M'_{1} \cap M_{3} &  \subset   & M'_{0} \cap M_{3}  \\
 \cup & & \cup \\
 M'_{1} \cap M_{2}&  \subset   &M'_{0} \cap M_{2}
\end{matrix}$
est exactement le carré commutatif $\;\;\begin{matrix}
B_1 & \subset   & A_0.B_1  \\
\cup & & \cup \\
A_{s} & \subset   &A_0
\end{matrix}$ au coin supérieur gauche de la grille. De plus $M'_{0} \cap M_{1}$ est l'algèbre $A_t$, donc $M'_{0} \cap M_{1}\subset M'_{0} \cap M_{2}\subset M'_{0} \cap M_{3}$ est la construction de base (voir 6.4.6) et les inclusions construites sont d'indice fini et de profondeur $2$.
\subsection{Actions de $A$ sur $M_1$ et de $B$ sur $M_2$}\label{actionbis}
Dans cette nouvelle construction, le commutant relatif $A_0=M'_{0} \cap M_{2}$ agit à gauche sur $M_1$ par son action à gauche sur $B_0$ contenu dans $M_1$.
\begin{prop}
Soient $i_{a} : a\mapsto[\cdots\otimes 1_{a} \otimes 1_{b}\otimes a]$ l'inclusion
de $A_0$ dans $M_2$ et $E_{M_1}$ l'espérance conditionnelle de $M_2$ sur $M_1$
conservant la trace $tr$. 
L'application de $A_0 \times M_1$ dans $M_1$
$$a \triangleright x = d^{-2}\gamma^2
E_{M_1}(i_{a}(a)xi_{a}(g_{s}^{1/2}f_{a}g_{s}^{-1/2}) )\quad(x\in M_1,
a\in A)$$
définit une action extérieure à gauche de $A_0$ sur $M_1$ telle que $M_2= M_1\rtimes A_0$.

Soient $i_{b} : b \mapsto[\cdots\otimes 1_{b} \otimes 1_{a}\otimes b]$ l'inclusion
de $B_1$ dans $M_3$ et $E_{M_2}$ l'espérance conditionnelle de $M_3$ sur $M_2$
conservant la trace $tr$. 
L'application de $B_1 \times M_2$ dans $M_2$
$$b \triangleright y = d^{-2}\gamma^2
E_{M_2}(i_{b}(b)yi_{b}(\hat{g}_{s}^{1/2}f_{b}\hat{g}_{s}^{-1/2}) )\quad(y\in M_2,
b\in B)$$
définit une action extérieure à gauche de $B_1$ sur $M_2$ telle que $M_3= M_2\rtimes B_1$.

De plus on a : $b \triangleright i_a(a) =i_a(b \triangleright a)$.
\end{prop}
   \dem
   Commençons par un lemme :
   
 \begin{lem}
Dans $B.A$, pour $a\in A$ et $b\in B$, on a :
$$[1_b\otimes a][b\otimes p_ag_s^{-1}]=[(a\triangleright b) \otimes p_ag_s^{-1}]$$
\end{lem}
 \dem
  \begin{align*}
  [1_b\otimes a][b\otimes p_ag_s^{-1}]&=[a_{(1)}\triangleright b\otimes a_{(2)}p_ag_s^{-1}]\\
&=[a_{(1)}\triangleright b\otimes \varepsilon_t(a_{(2)})p_ag_s^{-1}]
& \text{(d'après \ref{Haar})}\\
&=[(a_{(1)}\triangleright b)( \varepsilon_t(a_{(2)})\triangleright 1_b)\otimes p_ag_s^{-1}]\\
&=[(a_{(1)}\triangleright b)( a_{(2)}\triangleright 1_b)\otimes p_ag_s^{-1}]\\
&=[(a\triangleright b) \otimes p_ag_s^{-1}]
\end{align*}\end{proof}
 Notons $x=[z \otimes b]$ un élément  de $N^1_r$  avec $b \in B_0$ et  $z \in N_{r}^0$. D'après \ref{commutant}, comme $z$ appartient à $M_0$, $z$ commute à $A_0$. L'élément $i_{a}(a)xi_{a}(g_{s}^{1/2}f_{a}g_{s}^{-1/2})$ considéré dans $N^2_r$ s'écrit donc 
 $[z\otimes 1_{b}\otimes a][1_{N_{r}^0} \otimes b\otimes g_{s}^{1/2}f_{a}g_{s}^{-1/2}]$.
Le lemme permet de simplifier l'expression de l'action pour $x=[z \otimes b]$ :

   $$   a \triangleright x  =  d^{-2}\gamma^2
E_{M_1}(i_{a}(a)xi_{a}(g_{s}^{1/2}f_{a}g_{s}^{-1/2}) )  =  d^{-1}\gamma
E_{M_1}([z\otimes(a\triangleright b) \otimes p_ag_s^{-1}])$$
Le carré
$  \begin{matrix}
   N^2_r &\subset   & M_2\\
   \cup && \cup &\\
   N^1_r&  \subset  &M_1
\end{matrix}$ muni des espérances conditionnelles conservant $tr$ est
commutatif,  on en déduit : $a \triangleright x   =  d^{-1}\gamma
E_{N^1_r}([z\otimes(a\triangleright b) \otimes p_ag_s^{-1}])$
puis en appliquant la proposition 6.4.4, $a \triangleright x   =  d^{-1}\gamma
[z \otimes (a\triangleright b) E_{A_t}(p_ag_s^{-1})\triangleright 1_b)]$.
De plus, d'après \ref{gs} et \ref{escond2}, on a : 
$$E_{A_{t}}(p_{a}g_{s}^{-1})=d\gamma^{-1}F_{A_{t}}(g_{s}^{-1}p_{a})g_{t}^{-1}
=d\gamma^{-1}F_{A_{t}}(g_{t}^{-1}p_{a})g_{t}^{-1}
=d\gamma^{-1}g_{t}^{-1}F_{A_{t}}(p_{a})g_{t}^{-1}=d\gamma^{-1}$$
On obtient donc $a \triangleright x   = [z \otimes (a\triangleright b) ]$. L'application $x \mapsto a \triangleright x$ étend l'action
duale $[z\otimes b] \mapsto [z\otimes a \triangleright b]$ de $A_0$ sur
$N^1_r=N_{r}^0.B$. Elle définit une action à gauche
faiblement continue de $A_0$ sur $M_1$. De plus $ M_1 \rtimes A_0= M_1i_{a}(A_0)$
est le facteur $M_2$.
L'action est extérieure d'après \ref{ext} et \ref{commutant}.

L'étude de l'action de $B_1$ sur $M_2$ est analogue.
   \end{proof}

\subsection{Structures des commutants relatifs}
Dans $M'_0 \cap M_3$ regardons la dualité entre $A_0= M'_{0} \cap M_{2}$ et $B_1=M'_{1} \cap M_{3}$ définie pour ces commutants relatifs comme en \ref{dualite}.
Nous allons utiliser ici, pour la première fois, le fait que $A$ et $B$ sont réguliers.
\subsubsection{L'opérateur $h$}\label{hbis} Nous avons déjà mentionné en \ref {phiAsBt}, \ref{gs} et \ref{k} que si $A$ est régulier, la restriction de $\phi_a$ à $A_s$ est la trace canonique de $A_s$. D'autre part, la définition 6.4.1 relie la trace de Markov de l'inclusion  à $\phi_a$ par la formule : 
$tr(a)=d\gamma^{-2}\phi_a(g_{s}^{-1}g_{t}^{-1}a)$ pour $a \in A$. 
Grâce à \ref{gs}, on en déduit, pour $a\in A_s$, la relation $tr(a)=\gamma^{-1}\phi_a(g_{s}^{-1}a)$ entre la trace de Markov et la trace canonique sur $A_s$ donc $h^2$ est égal à $\gamma g_{s}$.
\subsubsection{Action de $B$ sur $A_0 \subset M_2$}\label{actiondeB1}
On peut alors reprendre les formules \ref{actionbis} de l'action d'un élément $b$ de $B_1$  sur un élément $i_a(a)$ de $A_0 \subset M_2$ en écrivant :
$$i_a(b \triangleright a)= b \triangleright i_a(a)=d^{-2}\gamma^2 E_{M_2}(i_{b}(b)[i_a(ah)\otimes f_{b}(h^{-1}\triangleright 1_b)])$$
\subsubsection{Dualité des commutants relatifs}
 
\begin{prop}
La dualité entre $A_0$ et $B_1$ résultant de la tour des facteurs (définie en\ref{dualite}) coïncide avec la dualité donnée initialement entre $A$ et $B$. Les structures de C*-groupoïdes quantiques résultant de la tour sont les structures initiales de $A$ et $B$.
\end{prop}
\begin{proof}
Notons $[a,b]$ la dualité entre $a\in A_0$ et $b\in B_1$ définie en\ref{dualite}.
Grâce à \ref{grille},\ref{commutant}, \ref{actionbis}, on obtient :
\begin{align*}
[a,b]&= d^{-4}\gamma^4 tr([i_a(ah)\otimes f_b][i_a(f_ah)\otimes 1_b]i_b(b))\\
&= d^{-4}\gamma^4 tr(i_b(b)[i_a(ah)\otimes f_b][i_a(f_ah)\otimes 1_b])\\
&= d^{-4}\gamma^4 tr(E_{M_2}(i_b(b)[i_a(ah)\otimes f_b(h^{-1}\triangleright 1_b)])[i_a(hf_ah)\otimes 1_b])\\
&=d^{-2}\gamma^2 tr([i_a(b\triangleright a) \otimes 1_b][i_a(hf_ah)\otimes 1_b])\quad \text{d'après \ref{actiondeB1}}\\
&=d^{-2}\gamma^2 tr_a((b\triangleright a)hf_ah) \\
&=d^{-1}\gamma^2 tr_a((b\triangleright a)p_a) \quad \text{d'après \ref{hbis}}\\
&=\phi_a(g_{s}^{-1}g_{t}^{-1}p_a(b\triangleright a))\\
&=\phi_a(g_{s}^{-2}p_a(b\triangleright a))\\
&=\langle g_{s}^{-2}p_a(b\triangleright a),p_b \rangle\\
&=\langle p_a(b\triangleright a),p_b\triangleleft g_{s}^{-2} \rangle \quad \text{d'après \ref{dual}}\\
&=\langle p_a(b\triangleright a),p_b(1_b \triangleleft g_{s}^{-2}) \rangle\quad \text{d'après 2.6.4}\\
\\
&=\langle p_a(b\triangleright a),p_b \hat{g}_{t}^{-2} \rangle \quad \text{d'après \ref{gs}}\\
\\
&=\langle (b\triangleright a),(p_b \hat{g}_{t}^{-2})\triangleleft p_a \rangle\\
\end{align*}
Or, par \ref{escond} et \ref{gs},  on a :
$$(p_b \hat{g}_{t}^{-2})\triangleleft p_a=F_{B_s}(p_b \hat{g}_{s}^{-2})=F_{B_s}(p_b)F_{B_s}(p_b)^{-1}=1$$

Donc, par \ref{dual}, les deux dualités entre $A_0$ et $B_1$ coïncident et comme les structures de C*- algèbres sont conservées, les structures de C*-groupoïdes quantiques aussi.
\end{proof}
\subsection{Conclusion} 
Nous réunissons dans le théorème suivant les résultats de cet annexe.
\begin{thm}
Etant donnés $A$ et $B$ deux C*-groupoïdes quantiques finis réguliers en dualité, on peut construire une tour de Jones (de profondeur $2$) de facteurs hyperfinis de type $\rm{II}_1$ :
 \begin{displaymath}
  M_0 \,\subset \, M_1 \, \stackrel{e_1}{\subset} \, M_2 \,\stackrel{e_2}{\subset} \, M_3 \end{displaymath} 
  telle que 
  \begin{enumerate}
\item le carré commutatif $\;\;\begin{matrix}
M'_{1} \cap M_{3} &  \subset   & M'_{0} \cap M_{3}  \\
 \cup & & \cup \\
 M'_{1} \cap M_{2}&  \subset   &M'_{0} \cap M_{2}
\end{matrix}$
soit isomorphe au carré commutatif $\;\;\begin{matrix}
B & \subset   & A.B  \\
\cup & & \cup \\
A_{s} & \subset   &A
\end{matrix}$ avec $e_1=[d\gamma^{-1} g_{s}^{-1/2}p_{a} g_{s}^{-1/2}\otimes 1_b]$ et $e_2=[1_a \otimes d\gamma^{-1} \hat{g}_{s}^{-1/2}p_{b} \hat{g}_{s}^{-1/2}]$
\item les structures de C*-groupoïdes quantiques sur les commutants relatifs obtenues à partir de la tour de facteurs grâce à la dualité \ref{dualite} coïncident avec les structures initiales de $A$ et $B$.
\item $A$ (resp. $B$) agisse à gauche sur $M_1$ (resp. $M_2$) de telle sorte que $M_2$ (resp. $M_3$) soit isomorphe à $M_1\rtimes A$ (resp.  $M_2\rtimes B$).
\end{enumerate}

\end{thm}

%\newpage
\section*{Références}

[Bi] D.Bisch : Bimodules, higher relative commutants and the fusion
algebra associated to a subfactor {\it The Fields Institutes for Research
in Mathmatical Sciences Communications Series 13 (1997), 13-63.}

[BNSz] G. Böhm, F. Nill et K. Szlach\'anyi : Weak Hopf algebras I.
Integral theory and C*-structure.{\it  J. Algebra 221 (1999) 385-438.}

[BSz1] G. Böhm et K. Szlachanyi :  A Coassociative C*-Quantum Group
with Non-Integral Dimensions.
{\it Lett. in Math. Phys., 35 (1996),437-456.}

[BSz2] G. Böhm et K. Szlach\'anyi : Weak Hopf algebras II.
Representation Theory, Dimension  and the Markov Trace.
{\it  J. Algebra 233 (2000) 156-212.}

   [Da1] M.-C. David : Paragroupe d'Adrian Ocneanu et algèbre de
Kac.{\it Pacific Journal of
   mathematics, Vol 172, No2, 1996.}
   
    [Da2] M.-C. David : Couple assorti de systèmes de Kac et inclusions de facteurs de type $II_1$.{\it Journal of functional analysis 159, 1-42 (1998).}
    
    [DT] M.-C. David et N. M. Thiéry :Exploration of finite dimensional Kac algebras
  and lattices of intermediate subfactors of irreducible inclusions (arXiv:math.QA/0812.3044v2)

   [GHJ]  F. M. Goodman, P. de la Harpe et V. F. R. Jones : Coxeter
Graphs and Towers of
   algebras. {\it MSRI Publications number 14.}

   [J]	V. Jones : Index for subfactors. {\it Invent. Math. 72 1-25  (1983).}
   
   [JS] V. Jones et V.S.Sunder : Introduction to Subfactors.{\it London 
Mathematical Society. Lecture Notes Series 234. Cambridge university 
press}

   [L] R.Longo : A duality for Hopf algebras and subfactors I. {\it
   Comm. Math. Phys. 159 (1994), 133-150}

   [N]  D. Nikshych : Duality for action of weak Kac algebras and
   crossed product inclusions of $\mathrm{II}_{1}$ factors.{\it journal
of Operator Theory 46 (2001) $n°$ 3 suppl. 635-655}

[NSzW] F. Nill, K. Szlach\'anyi  et H.-W. Wiesbrock : Weak Hopf
algebras and reducible Jones inclusions of depth 2, I : From crossed
products to Jones towers. {\it prépublication math.QA/9806130 (1998).}

[NV1] D. Nikshych et L. Vainerman :  A characterisation of depth 2
subfactors of  $\mathrm{II}_{1}$factors. {\it J. Func. Analysis 171
(2000) no. 2,
278-307.}

[NV2] D. Nikshych et L. Vainerman :  A Galois correspondence for
$\mathrm{II}_{1}$ factors and quantum groupoids. {\it J. Func.
Analysis 178 (2000)
113-142.}

[NV3] D. Nikshych et L. Vainerman :  Finite quantum groupoids and
their applications. {\it
"New Directions in Hopf Algebras",
Editors S. Montgomery and H.-J. Scheneider, MSRI Publications Vol. 43,
Cambridge University Press (2002), pp. 211 - 262.}

   [PiPo 1]	M. Pimsner et S. Popa :  Entropy and index for
   subfactors.{\it Ann.Scient.ENS 19 (1986)}
   p. 57-106

   [PiPo 2]	M. Pimsner et S. Popa :  Iterating the basic construction.
   {\it Trans.A.M.S. 310 (1988) No 1
   p.127-134.}

[Szy]  W. Szymanski : Finite index subfactors and Hopf algebras
crossed products. Proc.Amer.
   Math. Soc. 120 (1994) 519-528.
   
[V1] J.-M. Vallin. Groupoïdes quantiques finis.{\it Journal of Algebra
{\bf 239},215-261 (2001)}

[V2] J.-M. Vallin. Deformation of finite dimensional C*-Quantum 
Groupoids. math.QA/0310265.

[W] Y.Watatani. Index for C*-subalgebras.{\it Memoirs of the AMS 424
(1990)}.

\end{document}